\documentclass{article}[12pt]
\usepackage{latexsym,amssymb,amsthm,amsfonts,amsmath,epsfig,psfrag
}
\newtheorem{lemma}{Lemma}[section]
\newtheorem{theorem}{Theorem}
\newtheorem{cor}{Corollary}[section]
\newtheorem{prop}{Proposition}[section]
\theoremstyle{definition}

\def \N{{I\!\!N}}
\def \Z{{\mathbb Z}}
\def \R{{\mathbb R}}
\def \H{{\mathbb H}}

\def \S{{\mathbb S}}
\def \[{\langle }
\def \]{\rangle }
\def \M{{\cal M }}
\def \L{{\cal L }}

\newcommand{\conv}{\mathop{\rm conv}\nolimits}
\renewcommand{\dim}{\mathop{\rm dim}\nolimits}

\def \t{\tilde}

\def \o{\overline}
\def \w{\widetilde}
\def \lf{\left\lfloor}
\def \rf{\right\rfloor}
 
\begin{document}

\begin{center}
{\Large  \bf
}

\vspace{15pt}
{\Large \bf On compact hyperbolic Coxeter $d$-polytopes\\
with $d+4$ facets.}

\vspace{15pt}
{\large Anna Felikson\footnote{Partially supported by grants NS-5666.2006.1
and INTAS YSF-06-10000014-5916.},
Pavel Tumarkin\footnote{Partially supported by grants
MK-6290.2006.1, NS-5666.2006.1 and INTAS YSF-06-10000014-5766.} }

\end{center}

\vspace{11pt}

\begin{center}
\parbox{11cm}%
{\small 
{\it Abstract.}
We show that 
there is no compact hyperbolic Coxeter $d$-polytope with $d+4$ facets
for $d\ge 8$. This bound is sharp: examples of such polytopes up to
dimension $7$ were found by Bugaenko~\cite{Bu1}.
We also show that in dimension $d=7$ the polytope with 11 facets
is unique.

}
\end{center}

\medskip

\tableofcontents

\vspace{8pt}

\section*{Introduction}

A Coxeter polytope in the spherical, hyperbolic or
Euclidean space is a polytope whose dihedral angles are all 
integer submultiples of $\pi$.
These polytopes are very important among the other
acute-angled polytopes since 
a group generated by reflections with respect to the facets of a Coxeter
polytope is discrete.
On the other hand, a fundamental chamber of any discrete reflection 
group is a Coxeter polytope.
 
It is well known~\cite{C} that any spherical Coxeter polytope (containing no
pair of opposite points of the sphere) is a simplex
and any compact Euclidean Coxeter polytope is either a simplex or a 
direct product of simplices. 
See~\cite{C} for the description of these polytopes.

At the same time, hyperbolic Coxeter polytopes are still far
from being classified.
It is proved that no compact hyperbolic Coxeter polytope 
exists in dimensions $d\ge 30$~\cite{V_abs},
and no finite volume hyperbolic Coxeter polytope 
exists in dimensions $d\ge 996$~\cite{996}.
These bounds do not look sharp:
examples of compact polytopes are known up to dimension
$ 8$ only~\cite{Bu1},~\cite{Bu2}, and examples of finite volume 
polytopes are known up to dimension $21$ only~\cite{V_units}, 
~\cite{VK},~\cite{Bo}.

We will focus on compact Coxeter polytopes.
Besides the restriction on the dimension and some series of
examples~\cite{M},~\cite{ImH},~\cite{Al},~\cite{R}, 
there exists a classification of compact hyperbolic Coxeter polytopes of
certain combinatorial types.
More precisely, simplices are classified in~\cite{L},
$d$-polytopes with $d+2$ facets were listed in~\cite{K} and~\cite{Ess_eng},
and the classification of $d$-polytopes with $d+3$ facets
is contained in~\cite{n3}. 

\bigskip

In this paper we examine a class of the next complexity, 
namely, compact hyperbolic Coxeter $d$-polytopes with
$d+4$ facets. We show that no such a polytope exists in hyperbolic
space of dimension $d\ge 8$.
We also stress that in dimensions $2\le d \le 7$ these polytopes do  
exist~\cite{Bu1}.

The paper is organized as follows.
Section~\ref{preliminaries} is preparatory:
we recall basic notions concerning Coxeter diagrams
and combinatorics of simple polytopes.
We recall also some facts connecting the combinatorial
(metrical) properties of a face of a polytope
to the combinatorial
(metrical) properties of the polytope itself.
Section~\ref{section without small} is devoted to Coxeter diagrams
containing no Lann\'er diagram of order less than 5.
In particular, we prove that a Coxeter diagram of any 
compact hyperbolic Coxeter polytope contains a Lann\'er subdiagram of
order less than 5.
In Section~\ref{section_lift} we develop a theory of liftings 
to connect the combinatorics of a face of a Coxeter polytope to 
a subdiagram of the Coxeter diagram of this polytope.
Finally, in Sections~\ref{new} and~\ref{d8,n+4}  we use all these tools
to show the absence of the polytopes 
in dimensions $d\ge 8$.
 
The paper was partially written in the Max Plank Institute for Mathematics in Bonn.
The authors are grateful to the Institute for hospitality.

\section{Preliminaries}
\label{preliminaries}
In this section we list the essential facts about Coxeter diagrams,
Gale diagrams and diagrams of missing faces.
Concerning Coxeter diagrams we follow mainly~\cite{V_Refl_gp} and~\cite{29}.
For the details about Gale diagrams see~\cite{G}.
See also~\cite{Ess2} for an
overview about Coxeter polytopes, Gale diagrams  
and diagrams of missing faces.    
At the end of the section we recall a recent result of
Allcock~\cite{Al} which states that Coxeter polytopes often have
some Coxeter faces and describes the Coxeter diagrams of
that faces.

\subsection*{Coxeter diagrams}
\label{cox_def}

{\bf 1.}
{\it An abstract Coxeter diagram} $\Sigma$ is a finite 1-dimensional 
simplicial complex with weighted edges, where weights $w_{ij}$
are positive, and if $w_{ij}<1$ then $w_{ij}=\cos \frac{\pi}{m_{ij}}$ 
for some integer $m_{ij}\ge 3$.  
%
A {\it subdiagram} of $\Sigma$ is a subcomplex with the same weights
as in $\Sigma$.
%
The {\it order}  $|\Sigma|$ of the diagram  $\Sigma$
is the number of nodes of $\Sigma$.

If $\Sigma_1$ and $\Sigma_2$ are subdiagrams of an abstract
Coxeter diagram $\Sigma$, 
we denote by $\[\Sigma_1,\Sigma_2\]$ a subdiagram of $\Sigma$ 
spanned by all nodes of $\Sigma_1$ and $\Sigma_2$. 
Denote by $\Sigma_1\setminus v$ and $\Sigma_1\setminus \Sigma_2$
the subdiagrams of the diagram $\Sigma_1$
spanned by all nodes of $\Sigma_1$ except either $v$ 
or nodes of $\Sigma_2$ respectively.

Given an abstract Coxeter diagram $\Sigma$ with nodes
$v_1,\dots,v_n $ and weights $w_{ij}$, 
we
construct a symmetric $n\times n$ matrix $G(\Sigma)=(g_{ij})$, where 
$g_{ii}=1$, $g_{ij}= -w_{ij}$ if  $v_i$ and $v_j$ are  adjacent, and
$g_{ij}=0$ otherwise.
%
By $\det(\Sigma)$ and a signature of $\Sigma$ we mean the determinant 
and the signature of $G(\Sigma)$.


We can draw edges of Coxeter diagram in the following way:
if the weight $w_{ij}$ equals $\cos(\frac{\pi}{m_{ij}})$, 
$v_i$ and $v_j$ are joined by an $(m_{ij}-2)$-fold edge
or a simple edge labeled by $m_{ij}$;
if $w_{ij}=1$, $v_i$ and $v_j$ are joined by a bold edge;
if $w_{ij}>1$, $v_i$ and $v_j$ are joined by a dotted edge 
labeled by $w_{ij}$ (or without any label).

We write $[v_i,v_j]=m_{ij}$ if  $w_{ij}=\cos(\frac{\pi}{m_{ij}})$,
and  $[v_i,v_j]=\infty$ if $v_iv_j$ is a dotted edge.
We write $[v_i,v_j]=2$ if $v_i$ and $v_j$ are not joined.

\vspace{7pt}


An abstract Coxeter diagram $\Sigma$ is {\it elliptic} if
$G(\Sigma)$ is positive definite;
$\Sigma$ is {\it parabolic} if 
any indecomposable component of
$G(\Sigma)$ is degenerate and positive semidefinite;
a connected diagram $\Sigma$ is a {\it Lann\'er} diagram if 
$\Sigma$ is neither elliptic nor parabolic but
any proper subdiagram of $\Sigma$ is elliptic;
$\Sigma$  is {\it hyperbolic} if
$\Sigma$ is indefinite with negative inertia index equal to 1; 
$\Sigma$  is {\it superhyperbolic} if 
its negative inertia index is greater than 1;
$\Sigma$ is {\it admissible} if
$\Sigma$ contains no parabolic subdiagrams and
$\Sigma$ is 
not superhyperbolic.


\vspace{7pt}

Table~\ref{el-par} contains
the list of elliptic and connected parabolic Coxeter diagrams with 
their standard notation.
Lann\'er diagrams are listed in
~\cite[Table~3]{29}.
In particular, there are finitely many Lann\'er diagrams of order greater
than 3, and the maximal order of a Lann\'er diagram is 5.
We represent the list of Lann\'er diagrams of order 4 and 5
in Table~\ref{Lan}.

\begin{table}
\caption{Connected elliptic and parabolic Coxeter diagrams are listed
  in the left and in the right column respectively} 
\label{el-par}
\smallskip
\begin{center}
\begin{tabular}{|cc@{\quad}|cc|}
\hline
\raisebox{0pt}{${ A_n}$ $(n\ge 1)$}  & 
\raisebox{0pt}{\epsfig{file=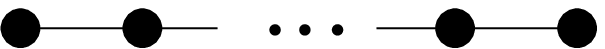,width=0.2\linewidth}}&
\multicolumn{2}{c|}{
\begin{tabular}{cc}
\raisebox{0pt}[20pt][5pt]{${ \widetilde A_1}$} & 
\raisebox{0pt}[20pt][5pt]{\epsfig{file=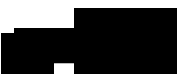,width=0.052\linewidth}}\\
\raisebox{3pt}[15pt][14pt]{${ \widetilde A_n}$ $(n\ge 2)$}  & 
\raisebox{-8pt}[25pt][7pt]{\epsfig{file=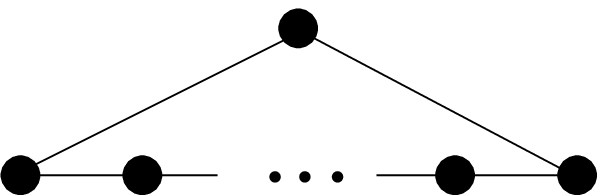,width=0.2\linewidth}}
\end{tabular}
}\\
\hline
\raisebox{-1pt}[23pt][7pt]{${ B_n=C_n}$} & 
\raisebox{-7pt}[23pt][7pt]{\epsfig{file=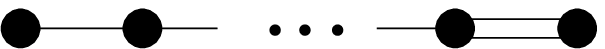,width=0.2\linewidth}}&
\raisebox{7pt}[23pt][7pt]{${ \widetilde B_n}$ $(n\ge 3)$}  & 
\raisebox{-0pt}[30pt][7pt]{\epsfig{file=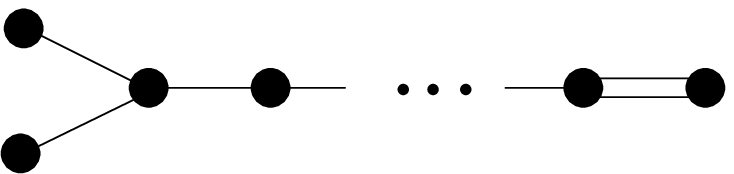,width=0.2\linewidth}}\\
\cline{3-4}
\raisebox{7pt}[15pt][7pt]{ $(n\ge 2)$} & 
&
\raisebox{-0pt}[15pt][7pt]{${ \widetilde C_n}$ $(n\ge 2)$}  & 
\raisebox{-0pt}[15pt][7pt]{\epsfig{file=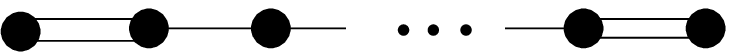,width=0.2\linewidth}}\\
\hline
\raisebox{7pt}[23pt][7pt]{${ D_n}$ $(n\ge 4)$} & 
\raisebox{0pt}[30pt][7pt]{\epsfig{file=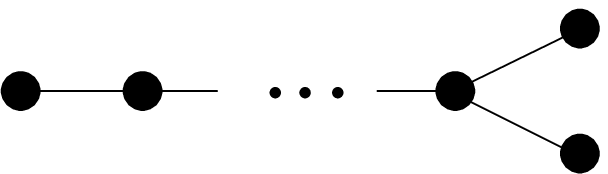,width=0.2\linewidth}}&
\raisebox{7pt}[23pt][7pt]{${ \widetilde D_n}$ $(n\ge 4)$} & 
\raisebox{0pt}[30pt][7pt]{\epsfig{file=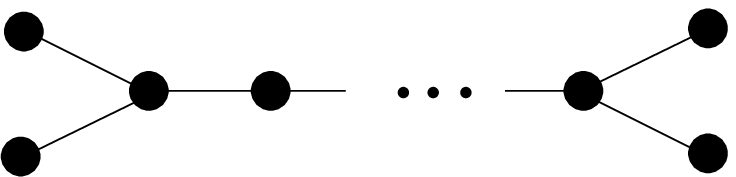,width=0.2\linewidth}}\\
\hline
\raisebox{0pt}[15pt][7pt]{${ G_2^{(m)}}$}  & 
\psfrag{m}{{\scriptsize $m$}}
\raisebox{0pt}[15pt][7pt]{\epsfig{file=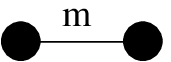,width=0.052\linewidth}}& 
\raisebox{0pt}[15pt][7pt]{${ \widetilde G_2}$} & 
\raisebox{0pt}[15pt][7pt]{\epsfig{file=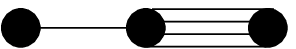,width=0.082\linewidth}}\\
\hline
\raisebox{0pt}[15pt][7pt]{${ F_4}$}  & 
\raisebox{0pt}[15pt][7pt]{\epsfig{file=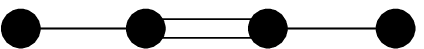,width=0.13\linewidth}}&
\raisebox{-0pt}[15pt][7pt]{${ \widetilde F_4}$} & 
\raisebox{-0pt}[15pt][7pt]{\epsfig{file=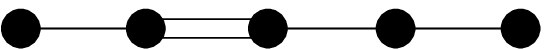,width=0.16\linewidth}}\\
\hline
\raisebox{8pt}[15pt][7pt]{${ E_6}$}  & 
\raisebox{0pt}[30pt][7pt]{\epsfig{file=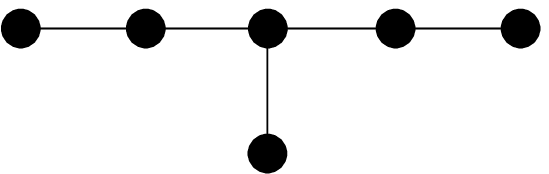,width=0.16\linewidth}}&
\raisebox{8pt}[35pt][7pt]{${ \widetilde E_6}$} & 
\raisebox{-8pt}[40pt][17pt]{\epsfig{file=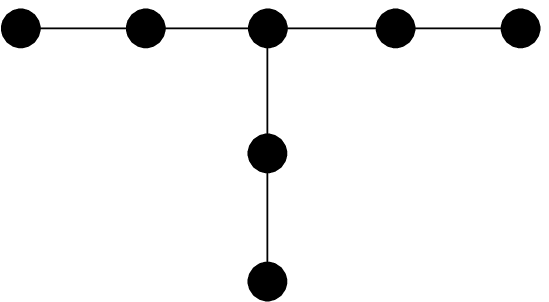,width=0.16\linewidth}}\\
\hline
\raisebox{5pt}[25pt][7pt]{${ E_7}$}  & 
\raisebox{-0pt}[30pt][7pt]{\epsfig{file=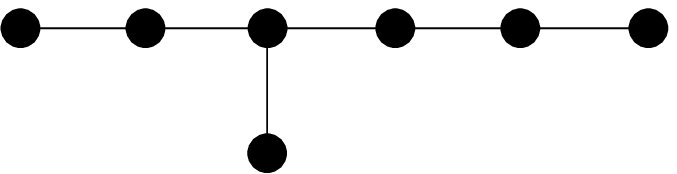,width=0.2\linewidth}}&
\raisebox{5pt}[25pt][7pt]{${ \widetilde E_7}$} & 
\raisebox{-0pt}[25pt][7pt]{\epsfig{file=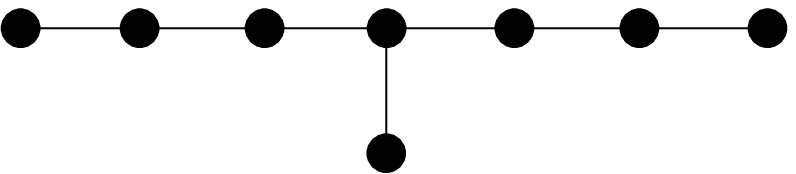,width=0.22\linewidth}}\\
\hline
\raisebox{5pt}[25pt][7pt]{${ E_8}$}  & 
\raisebox{-0pt}[30pt][7pt]{\epsfig{file=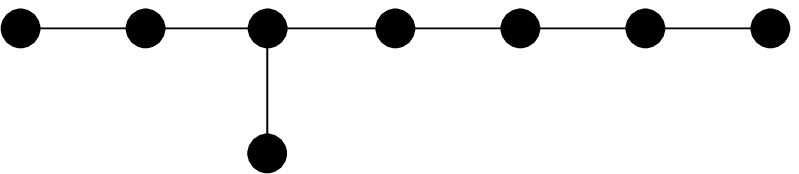,width=0.24\linewidth}}&
\raisebox{5pt}[25pt][7pt]{${ \widetilde E_8}$} & 
\raisebox{-0pt}[25pt][7pt]{\epsfig{file=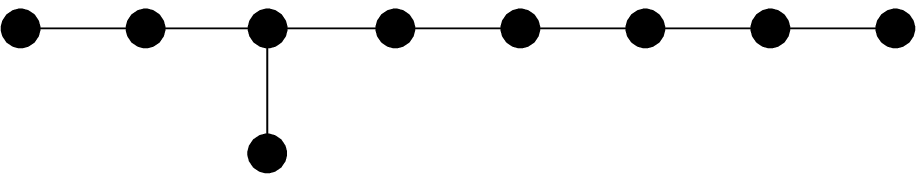,width=0.25\linewidth}}\\
\hline
\raisebox{0pt}[15pt][7pt]{${ H_3}$}  & 
\raisebox{0pt}[15pt][7pt]{\epsfig{file=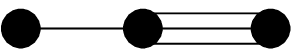,width=0.1\linewidth}}&
& 
\\
\cline{1-2}
\raisebox{0pt}[15pt][7pt]{${ H_4}$}  & 
\raisebox{0pt}[15pt][7pt]{\epsfig{file=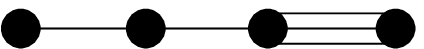,width=0.13\linewidth}}&
& 
\\
\hline

\end{tabular}
\end{center}

\end{table}


\bigskip

{\bf 2.}
It is convenient to describe Coxeter polytopes by their Coxeter diagrams.
Let $P$ be a Coxeter polytope with facets $f_1,\dots,f_r$.
The Coxeter diagram $\Sigma(P)$ of the polytope $P$ 
is a diagram with nodes $v_1,\dots,v_r$; 
two nodes $v_i$ and $v_j$ are not joined if $f_i$ is orthogonal
to $f_j$;
 $v_i$ and $v_j$ are joined by an edge with weight
$$w_{ij}=
\begin{cases}  
\cos \frac{\pi}{k},&\text{ if  $f_i$ and $f_j$ form a  dihedral angle 
$\frac{\pi}{k}$;}\\
1,&\text{ if $f_i$ is parallel to $f_j$;}  \\
\cosh \rho,&\text{ if $f_i$ and $f_j$ diverge and $\rho$ is the distance 
from $f_i$ to $f_j$.}
\end{cases}
$$

\noindent
If $\Sigma=\Sigma(P)$ then $M(\Sigma)$ coincides with the Gram matrix
of unit normal vectors to the facets of $P$.

It is shown in~\cite{V_Refl_gp} that if $\Sigma=\Sigma(P)$ 
is a Coxeter diagram 
of a compact $d$-dimensional
hyperbolic polytope $P$, then $\Sigma$ is an admissible connected hyperbolic 
 diagram with positive inertia index equal to $d$.
In particular, $\Sigma$ contains no bold edges and other 
parabolic subdiagrams.
Elliptic subdiagrams of $\Sigma$ are in one-to-one correspondence
with faces of $P$: a $k$-face $F$ corresponds to an elliptic
subdiagram $\Sigma_F$ of order $d-k$ whose nodes correspond to
the $(d-1)$-faces containing $F$.

\bigskip

{\bf 3.}
Given a Coxeter diagram $\Sigma$, it is easy to check
if $\Sigma$ is superhyperbolic or not.
However, this requires the computation of a signature of a rather big
matrix if the order of $\Sigma$ is not small.
In the case when $\Sigma$ is a union of two subdiagrams either joined by
only one edge or having a unique node in common, there exists a  more effective way
to state that $\Sigma$ is superhyperbolic~\cite{V_abs}.

Let $T$ be a subdiagram of $\Sigma$ such that 
$\det(\Sigma\setminus T)\ne 0$.
A {\it local determinant} of $\Sigma$ on a subdiagram $T$ is
$\det(\Sigma,T)=\frac{\det (\Sigma)}{\det(\Sigma\setminus T)}$.

\begin{prop}[\cite{V_abs}, Prop.~12]
\label{super1}
If a Coxeter diagram $\Sigma$ consists of two subdiagrams $\Sigma_1$
and $\Sigma_2$ having a unique vertex $v$ in common, and nodes of
$\Sigma_1\setminus v$ do not attach to $\Sigma_2\setminus v$,  then
$$
\det(\Sigma,v)=\det(\Sigma_1,v)+\det(\Sigma_2,v)-1.
$$  

\end{prop}

\begin{prop}[\cite{V_abs}, Prop.~13]
\label{super2}
If a Coxeter diagram $\Sigma$ is generated by disjoint subdiagrams 
$\Sigma_1$ and $\Sigma_2$ joined by a single edge
$v_1v_2$, then
$$
\det(\Sigma,\[v_1,v_2\])=\det(\Sigma_1,v_1)\det(\Sigma_2,v_2)-w_{12}^2,
$$  
where $w_{12}$ is the weight of the edge $v_1v_2$.

\end{prop}

\begin{prop}[\cite{V_abs}, Prop.~15]
\label{super3}
Suppose that a Coxeter diagram $\Sigma$ is a union of two disjoint 
hyperbolic subdiagrams $\Sigma_1$ and $\Sigma_2$
joined by a unique edge $v_1v_2$
and that the subdiagrams $\Sigma_1 \setminus v_1$ and 
$\Sigma_2 \setminus v_2$ are elliptic.
Assume that in addition  one of the following conditions holds:\\
1) $v_1v_2$ is a simple edge and \
$\det(\Sigma_1,v_1)\det(\Sigma_2,v_2)>\frac 1 4$;\\
2) $v_1v_2$ is a double edge and \
$\det(\Sigma_1,v_1)\det(\Sigma_2,v_2)>\frac 1 2$.

Then the diagram $\Sigma$ is superhyperbolic.

\end{prop}

In~\cite[Table~2]{V_abs} Vinberg listed some useful local determinants.
We use Propositions~\ref{super1}--~\ref{super3} together with 
Table~2 of~\cite{V_abs} 
throughout this paper not referring to them every time.
We also use the fact that all local determinants of Lann\'er diagrams
shown in Table~\ref{Lan} (of the present paper) on their open vertices 
(see the definition below) are greater than 0.95, this fact
can be checked by a  direct computation.   
In particular, we will use the following corollary of 
Prop.~\ref{super3}:

\begin{prop}
\label{2lan}
Suppose that a diagram $\Sigma$ consists of two disjoint Lann\'er
diagrams $L_1$ and $L_2$ of order $5$ each joined by a unique non-dotted 
edge $v_1v_2$, such that $L_i\setminus v_i$ (where
 $v_i\in L_i$, $i=1,2$) is of the 
type $H_4$ or $F_4$. Then $\Sigma$ is superhyperbolic. 

\end{prop}

\bigskip

{\bf 4.}
Let $\Sigma$ be an abstract Coxeter diagram and $v$ be a node of $\Sigma$.
Suppose that it is possible to construct an abstract Coxeter
diagram $\[\Sigma,x\]$, $x\notin \Sigma$, such that 
$x$ is joined with $v$ in $\[\Sigma,x\]$ and 
for any non-elliptic subdiagram $\Sigma'\subset \[\Sigma,x\]$
the intersection $\Sigma'\cap \Sigma=\Sigma'\setminus x$ 
is non-elliptic, too. 
Then $v$ is called an {\it open vertex} of $\Sigma$.

The vertex $v$ is called {\it doubly open}
if it is possible to construct a connected Coxeter
diagram $\[\Sigma,x_1\]$, $x_1\notin \Sigma$, 
such that  $x_1$ is joined with $v$ in $\[\Sigma,x_1\]$ 
and $x_1$ is an open vertex of  $\[\Sigma,x_1\]$.  

\vspace{7pt}

If $v\in Sigma$ is not an open vertex of a Lann\'er diagram $\Sigma$
and $x\notin \Sigma$ is a node joined with $v$,
then the diagram $\[\Sigma, x\]$ contains a non-elliptic subdiagram
$M$ such that $x,v\in M$.
Indeed, if $x$ is not connected with $\Sigma \setminus v$,
the statement is obvious. If we add some edges connecting $x$ and 
$\Sigma \setminus v$ then the non-elliptic subdiagrams of $\[\Sigma,x\]$
remain non-elliptic.


In~\cite[Tabelle~2]{Ess2} Esselmann listed all Lann\'er diagrams
of order 4 and 5 containing open vertices.
In Table~\ref{Lan} we list Lann\'er diagrams
of order 4 and 5 and introduce the notation.
The open vertices are marked black,
the doubly open vertices are encircled.

\begin{table}[htb]
\begin{center}
\caption{Lann\'er diagrams of order 4 and 5.
Open vertices are marked black,
doubly open vertices are marked by additional circles.
The upper index equals the order of the diagram.
}
\label{Lan}
\psfrag{41}{${\cal L}^4_1:$}
\psfrag{42}{${\cal L}^4_2:$}
\psfrag{43}{${\cal L}^4_3:$}
\psfrag{44}{${\cal L}^4_4:$}
\psfrag{45}{${\cal L}^4_5:$}
\psfrag{46}{${\cal L}^4_6:$}
\psfrag{47}{${\cal L}^4_7:$}
\psfrag{48}{${\cal L}^4_8:$}
\psfrag{49}{${\cal L}^4_9:$}
\psfrag{51}{${\cal L}^5_1:$}
\psfrag{52}{${\cal L}^5_2:$}
\psfrag{53}{${\cal L}^5_3:$}
\psfrag{54}{${\cal L}^5_4:$}
\psfrag{55}{${\cal L}^5_5:$}
\epsfig{file=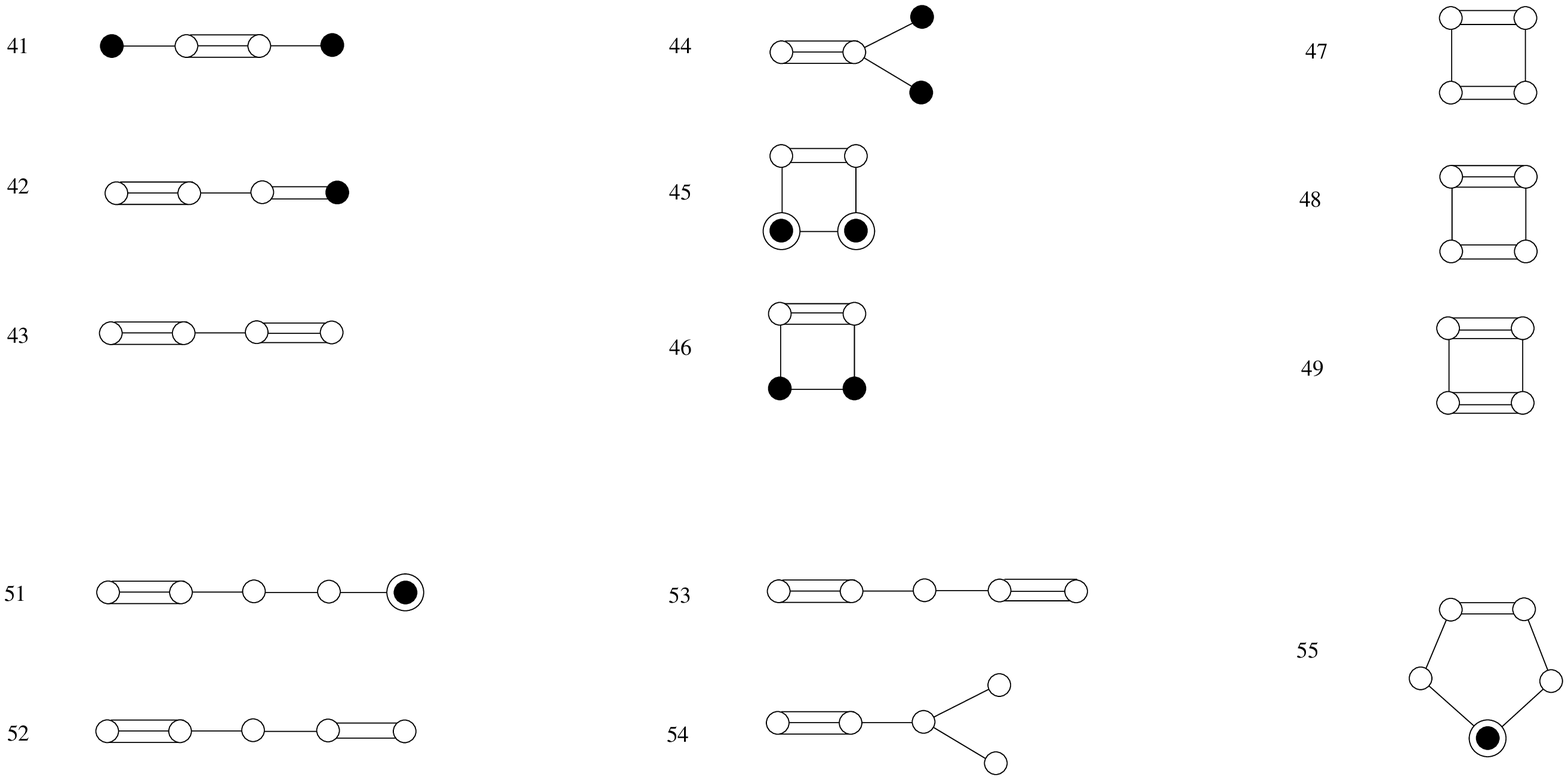,width=0.75\linewidth}
\end{center}
\end{table}

\begin{prop}
\label{2subd_of_L4}
Any Lann\'er diagram of order $4$ contains at least two subdiagrams of the
types $H_3$ or $B_3$.

\end{prop}

This follows directly from the classification of Lann\'er diagrams.

\begin{prop}[\cite{V_abs}, Prop.~2]
\label{joined}
Any two Lann\'er diagrams of $\Sigma$ are joined by at least one edge.

\end{prop}

The statement is obvious: if it is false then $\Sigma$ is a
superhyperbolic diagram.

\subsection*{Gale diagrams and diagrams of missing faces}

{\bf 1.}
As it is shown in~\cite[cor. of Th.~3.1]{V_Refl_gp},
any compact hyperbolic Coxeter polytope is simple
(i.e. a $k$-dimensional face of a $d$-polytope is an intersection of
exactly $d-k$ facets).
In this section we discuss the combinatorics of simple polytopes.
We abuse notation by writing ``$d$-polytope'' instead of
``$d$-dimensional polytope''. 

Every combinatorial type of simple $d$-polytope with $d+k$ facets can be
represented by its {\it Gale diagram}. This consists of $d+k$ points
$a_1,\dots,a_{d+k}$ of $(k-2)$-dimensional unit sphere 
$\S^{k-2}\subset \R^{k-1}$
centered  at the origin.
Each point $a_i$ corresponds to a facet
$f_i$ of $P$.
The combinatorial type of a simple convex polytope can be read off from
the Gale diagram in the following way:
%
for any $J\subset \{1,\dots,d+k\}$ 
the intersection of facets
$\{f_j | j\in J \}$
is a face of $P$ if and only if the origin is contained in the interior of
$\conv\{a_{j} | j\notin J\}$ 
(where $\conv X$ is a convex hull of the set $X$).



The set of points  $\{a_1,\dots,a_{d+k}\}\subset \S^{k-2}$
is a Gale diagram of some $d$-dimensional polytope $P$ with $d+k$ facets
if and only if every open halfspace $H^+$ in $\R^{k-1}$
bounded by a hyperplane through the origin
contains at least two of  points $a_1,\dots,a_{d+k}$.


Two Gale diagrams are called {\it isomorphic} if the corresponding polytopes
are combinatorially equivalent.
 
\medskip

Let $P$ be a simple polytope.
The facets $f_1,\dots,f_m$ of $P$ compose a {\it missing face} of $P$
if $\bigcap\limits_{i=1}^m f_i=\emptyset$ but
any proper subset of $\{f_1,\dots,f_m\}$ has a non-empty intersection.
Clearly, if a  set of facets have no  common intersection then this set of
facets contains at least one missing face.
This leads to the following lemma:


\begin{lemma}
\label{min_in_a_halfspace}
Let $P$ be a simple $d$-polytope with $d+k$ facets and let
$\{a_1,\dots,a_{d+k}\}\subset \S^{k-2}$ be a Gale diagram of $P$.
Let  $H^+$ be any open halfspace bounded  by a hyperplane through the
origin. Then $H^+$ contains a set $I\subset \{a_1,\dots,a_{d+k}\}$
which corresponds to a missing face of $P$.

\end{lemma}

%

\bigskip

{\bf 2.}
If  $k=2$ then the Gale diagram of $P$ is one-dimensional, i.e. vertices
$a_i$ of the diagram lie on the 0-dimensional unit sphere.
In other words, each of the points $a_i$ coincides either with the
point $-1$ or with the point $1$. This leads to the following

\begin{prop}[\cite{G}]
\label{prod}
A simple $d$-polytope with $d+2$ facets is a direct product of two simplices
$\Delta^{n-k} \times \Delta^{k}$ (where $0\le k\le [n/2]$ and
$\Delta^m$ stays for an $m$-dimensional simplex).

\end{prop}

As it is shown in~\cite{Ess_eng}, 
any compact Coxeter $d$-polytope with 
$d+2$ facets is either a simplicial prism or a product of two
triangles. Coxeter prisms are listed in~\cite{K}, 
and the remaining Coxeter polytopes of this type are listed in~\cite{Ess_eng}
(there are 7 of them). We call these  7 polytopes 
{\it Esselmann polytopes}.

If  $k=3$ then the Gale diagram of $P$ is two-dimensional, i.e. vertices
$a_i$ of the diagram lie on the unit circle.

A {\it standard Gale diagram} of simple $d$-polytope with $d+3$ facets
consists of vertices
$v_1,\dots,v_m$ of regular $m$-gon ($m$ is odd) in $\R^2$
centered  at the
origin, which are labeled according to the following rules:

1) Each label is a positive integer, the sum of labels equals $d+3$.

2) The vertices that lie in any open halfspace bounded by a line
 through the origin have labels whose sum is at least two.

It is easy to check (see, for example,~\cite{G}) that any
two-dimensional Gale diagram is isomorphic to some standard
diagram. Two simple $d$-polytopes with $d+3$ facets are
combinatorially equivalent if and only if their standard Gale diagrams
are congruent.

\bigskip


For $k>3$ there is no notion of standard Gale diagram known to us.
We will use another diagram to encode the combinatorics of 
a simple polytope.

\medskip

A {\it diagram of missing faces} is a finite set
$D$ with a specified collection
  $\M_D$ of subsets of $D$ satisfying the following
condition: if $M, M'\in \M_D$ then 
$M'\not\subset M$.
The order $|M|$ of $M$ is the number of elements in $M$. 
Elements of $\M_D$ are  called {\it missing faces} of $D$. 

A diagram $D_1\subset D$ of missing faces is a {\it subdiagram}
of $D$ if for any $M\subset D_1$ the following holds: 
$M\in \M_{D_1}$ if and only if $M\in \M_{D}$.


\medskip

To represent a diagram of missing faces
we draw a point (vertex) for each element of $D$
and encircle a set of points corresponding to $M$  
if and only if $M\in \M_D$.

\bigskip
{\bf 3.}
For any simple polytope $P$ we construct a diagram of missing faces 
$D(P)$ in the following way:
the elements correspond to the facets of $P$,
a set of elements is a missing face of $D(P)$ if and only if
the corresponding facets compose a missing face of $P$.


The combinatorics of $P$ can be recovered from $D(P)$:
the set of facets have a non-empty intersection if and only if
the corresponding set in $D(P)$ contains no missing faces.

\begin{lemma}
\label{disjoint}
Let $P$ be a simple polytope and $D(P)$ be a diagram of missing faces
of $P$. 
For any $M\in \M_{D(P)}$ there exists $M'\in \M_{D(P)}$ 
such that $M\cap M'=\emptyset$.

\end{lemma}

\begin{proof}
Consider a Gale diagram $G(P)$ of $P$.
It consists of several points on a 
sphere $\S^d$ for some $d$.
Let $\overline M$ be the points of $G(P)$ corresponding to the vertices
of $M$. Since $M$ is a missing face, 
$\conv (G(P)\setminus \overline M)$ does not contain 
the origin. Equivalently, there exists a hyperplane $H$ through 0
separating a halfspace $H^+$ containing $G(P)\setminus \overline M$.
Denote by $H^-$ another halfspace with respect to $H$.
By Lemma~\ref{min_in_a_halfspace},  $H^-$
contains at least one subset corresponding to some
missing face. Since this subset belongs to $\overline M$,
and no missing face contains another missing face, we obtain that
$H$ separates  $G(P)\setminus \overline M$ from $\overline M$,
i.e.  $H^-\cap G(P)=\overline M$. 
By Lemma~\ref{min_in_a_halfspace},  
$H^+$ contains at least one subset corresponding to some missing face $M'$.
Clearly, $M$ and $M'$ are disjoint and the lemma is proved.

\end{proof}

We also need the following two Propositions:

\begin{prop}[\cite{Ess2}, Lemma~1.6]
\label{face}
Let $P$ be a simple polytope and $f$ be a facet of $P$.
Let $\{f_1,\dots,f_k\}$ be the set of all facets of $P$ such that
$f_i\cap f\ne 0$ and $f_i\ne f$.
Denote $f'_i=f_i\cap f$ for $i=1,\dots, k$, 
and for any subset $F\subseteq \{f_1,\dots,f_k\}$
denote $F'=\{f'_i \, |\, f_i\in F\}$. 

The set of faces $F'$ is a missing face of $f$ if and only if
either
\begin{itemize}
\item[1)] $\{f\}\cup F$ is a missing face of $P$
or 
\item[2)] $F$ is a missing face of $P$ and $F$ contains no proper subset
$F_0$ such that $\{f\}\cup F_0$ is a missing face of $P$. 

\end{itemize}
\end{prop}

\begin{prop}[\cite{Ess2}, Lemma~1.9]
\label{lanner_for_any_facet}
Let $P$ be a simple polytope, and let $f$ be any facet of $P$.
There exists a missing face of $P$ containing $f$.

\end{prop}

\subsection*{Lann\'er diagrams and missing faces}

Let $P$ be a compact Coxeter polytope in $\H^n$ and
$\Sigma(P)$ be the Coxeter diagram of $P$.
Let $L$ be any Lann\'er subdiagram of $\Sigma(P)$.
By definition of Lann\'er diagrams,
 the facets corresponding to $L$ 
compose a missing face of $P$
(and any missing face of $P$ corresponds to some Lann\'er diagram
in $\Sigma(P)$).
Thus, a diagram of missing faces $D(P)$ can be easily reconstructed
by $\Sigma(P)$: to obtain $D(P)$ one should take $\Sigma(P)$,
encircle all Lann\'er diagrams and
delete all edges. 

In the same way we can construct a diagram $D(\Sigma)$ of missing faces
for any admissible Coxeter diagram $\Sigma$.

The correspondence ``Lann\'er diagrams $\longleftrightarrow$ 
missing faces'' shows in particular
that if $P$ is a compact hyperbolic Coxeter
polytope,
then $P$ has no missing faces of order greater than 5.


\subsection*{Faces of Coxeter polytopes}

Let $P$ be a compact hyperbolic Coxeter $d$-polytope, and denote by $\Sigma$
its Coxeter diagram. Let $S_0$ be an elliptic subdiagram of $\Sigma$. 
By~\cite[Th.~3.1]{V_Refl_gp}, $S_0$ corresponds to a face of $P$ of dimension
$d-|S_0|$. Denote this face by $P(S_0)$. $P(S_0)$ itself is an
acute-angled polytope~\cite{A}, but it 
may not be a Coxeter polytope. Borcherds obtained the following
sufficient condition for $P(S_0)$ to be a Coxeter polytope.

\smallskip

\begin{prop}[\cite{Bo}, Example 5.6] 
\label{bor}
Suppose $P$ is a Coxeter polytope with diagram $\Sigma$, and 
$S_0\subset \Sigma$
is an elliptic subdiagram that has no
$A_n$ or $D_5$ component. Then $P(S_0)$ itself is a Coxeter polytope.

\end{prop}

\smallskip

Facets of $P(S_0)$ correspond to those nodes of $\Sigma$ that compose an elliptic
diagram together with $S_0$. The following result of Allcock shows how
to compute dihedral angles of $P(S_0)$. 

We say that a node of $\Sigma$ {\it attaches} to 
$S_0$ if it is joined with any node of $S_0$ by an edge of any type.
Let $a$ and $b$ be the facets of $P(S_0)$ coming from facets $A$ and $B$ of
$P$, i.e. $a=A\cap P(S_0)$ and $b=B\cap P(S_0)$. 
Denote by $v_A$ and $v_B$ the nodes of $\Sigma$ corresponding to the
facets $A$ and $B$.  Then the
angles of $P(S_0)$ can be computed in the following way.

\smallskip

\begin{prop}[\cite{Al}, Th. 2.2.]
\label{al}
Under the hypotheses of Prop.~\ref{bor},
\begin{itemize}
\item[(1)]
If neither $v_A$ nor $v_B$ attaches to $S_0$, then $\angle ab=\angle AB$. 
\item[(2)]
If just one of $v_A$ and $v_B$ attaches to $S_0$, say to the component 
$S_0^i$, then 
\begin{itemize}
\item[(a)]
if $A\perp B$ then $a\perp b$;
\item[(b)] 
if $v_A$ and $v_B$ are joined and adjoining $v_A$ and $v_B$ 
to $S_0^i$
yields a diagram $B_k$ (resp. $D_k$, $E_8$ or $H_4$) then 
$\angle ab=\pi/4$ (resp. $\pi/4$, $\pi/6$ or $\pi/10$);
\item[(c)] 
otherwise, $a$ and $b$ do not meet.
\end{itemize} 
\item[(3)] 
If $v_A$ and $v_B$ attach to different components of $S_0$, then
\begin{itemize}
\item[(a)] 
if $A\perp B$ then $a\perp b$;
\item[(b)] 
otherwise, $a$ and $b$ do not meet.
\end{itemize} 
\item[(4)] 
If $v_A$ and $v_B$ attach to the same component of $S_0$, say $S_0^i$,
then 
\begin{itemize}
\item[(a)] 
if $v_A$ and $v_B$ are not joined and $S_0^i\cup \{A,B\}$ is a diagram
$E_6$ (resp. $E_8$ or $F_4$) then 
$\angle ab=\pi/3$ (resp. $\pi/4$ or $\pi/4$);
\item[(b)] 
otherwise, $a$ and $b$ do not meet.
\end{itemize} 
\end{itemize}

\end{prop}

\smallskip

Let $w\in \Sigma$ be a {\it neighbor} of $S_0$,
so that $w$ attaches to $S_0$ by some edges.
We call $w$ {\it good} if
$\[S_0,w\]$ is an 
elliptic diagram, and {\it bad} otherwise. We denote by $\overline S_0$ the
subdiagram of $\Sigma$ consisting of nodes corresponding to facets of
$P(S_0)$. The diagram $\overline S_0$ is spanned by good neighbors of
$S_0$ and by all nodes not joined to $S_0$.
If $P(S_0)$ is a Coxeter
polytope, denote its Coxeter diagram by $\Sigma_{S_0}$.   
By ordinary edge we mean non-dotted edges.
By simple edges we mean 1-fold edges.
By empty edge we mean two nodes which are not joined.

For any node $v\in \Sigma$ which is not a bad neighbor of
 $S_0$ (i.e. $v\in \o S_0$) we denote by
$\t v$ the corresponding node of $\Sigma_{S_0}$.   

If $\Sigma_{S_0}$ does not differ from $\o S_0$,
we consider the diagram $\Sigma_{S_0}$ as a subdiagram of 
$\Sigma$, and we do not differ $v$ and $\t v$.  

\bigskip
\bigskip

\begin{cor}
\label{dif}
Under the hypotheses of Prop.~\ref{bor},
\begin{itemize}
\item[(a)]
If $S_0$ is of the type $H_4$, $F_4$, or $G_2^{(m)}$ for $m\ge 6$,
or any other diagram having no good neighbors, then
$\overline S_0=\Sigma_{S_0}$. 
\item[(b)]
If $S_0$ is of the type $H_3$, then $\overline S_0$ may be obtained from $\Sigma_{S_0}$ by
replacing some dotted edges by ordinary edges.
\item[(c)]
If $S_0$ is of the type $G_2^{(5)}$, then $\overline S_0$ may be obtained from $\Sigma_{S_0}$ by
replacing some edges labeled by 10 by simple edges, and some dotted
edges by ordinary edges.
\item[(d)]
If $S_0$ is of the type $B_n$, $n\ge 3$, then  $\overline S_0$ 
may be obtained
from $\Sigma_{S_0}$ by replacing some double edges by simple edges, and some dotted
edges by ordinary edges.
\item[(e)]
If $S_0$ is of the type $B_2=G_2^{(4)}$, then  $\overline S_0$ may be obtained
from $\Sigma_{S_0}$ by replacing some double edges by simple edges, and some dotted
edges by ordinary or empty edges.
\end{itemize}
\end{cor}

The corollary immediately follows from Prop.~\ref{al}. One should 
only notice that all neighbors of diagrams listed in item (a) are bad.

\begin{cor}
\label{dif2}
Under the hypotheses of Prop.~\ref{bor},
let $S_1\subset \Sigma_{S_0}$ be a subdiagram of the type 
$G_2^{(m)}$, where $m\ne
4, 10$, and let $S_1'$ be the corresponding subdiagram of
$\o S_0$.
Then no node of $S_1'\subset \Sigma$ is a good neighbor
of $S_0$.


In particular, any subdiagram
$S_2\subset \Sigma_{S_0}$ of the type 
 $F_4$, $H_4$, $H_3$ or $G_2^{(m)}$, where $m\ne 4, 10$,
corresponds to a subdiagram of the same type in $\o S_0$.

Any Lann\'er subdiagram $L\subset \Sigma_{S_0}$ of order $5$
corresponds to a Lann\'er subdiagram of $\o S_0$.

\end{cor}

\begin{lemma}
\label{bad}
Suppose that $S_0\subset\Sigma$ is an elliptic subdiagram and $|S_0|<d$. 
Then $S_0$ has at most 
$|\Sigma|-d-1$ bad
neighbors. In particular, if $P$ has $d+4$ facets, then any elliptic
subdiagram of $\Sigma$ of order less than $d$ has at most $3$ bad neighbors. 

\end{lemma}

\begin{proof}
The lemma follows from the fact that a
$d$-polytope has at least $d+1$ facets. 

\end{proof}

\begin{lemma}
\label{bad&Lanner}
Let $S\subset \Sigma$ be an elliptic subdiagram containing no
component of the types $A_n$ and $D_5$,
and let $v$ 
be a bad neighbor of $S$. Then $v$ is joined with each Lann\'er
subdiagram of $\o S$.

\end{lemma}

The statement immediately follows from Proposition~\ref{joined}.

\section{Admissible Coxeter diagrams without small Lann\'er diagrams}
\label{section without small}

A Lann\'er diagram (a missing face) $L$ is {\it small} if $|L|<5$.
A missing face $M$ is {\it large } if $|M|>5$.

\begin{figure}[htb]
\begin{center}
\psfrag{1}{$\Theta_1$}
\psfrag{2}{$\Theta_2$}
\psfrag{3}{$\Theta_3$}
\epsfig{file=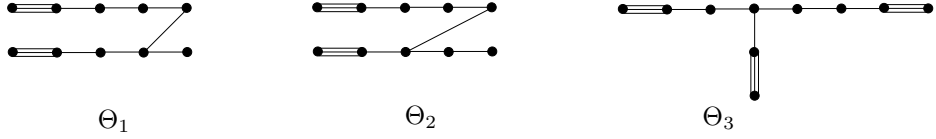,width=0.85\linewidth}
\caption{Diagrams without small Lann\'er  subdiagrams. }
\label{10only5}
\end{center}
\end{figure}

\begin{lemma}
\label{Table5}
Let $\Sigma$ be a connected
admissible Coxeter diagram containing
no small Lann\'er subdiagrams.
Suppose that each node of $\Sigma$ belongs to some Lann\'er subdiagram 
of $\Sigma$. Then $|\Sigma|\le 10$.
If $|\Sigma|=10$ then $\Sigma$ is one of the three diagrams
shown in Fig.~\ref{10only5}.

If in addition $\det (\Sigma)=0$ then $\Sigma=\Theta_1$ 
(see Fig.~\ref{10only5}).

\end{lemma}

\begin{proof}

Since $\Sigma$ is a connected diagram and $\Sigma$ 
contains no small Lann\'er subdiagrams, any edge of $\Sigma$ is either
a simple edge, or a double edge, or a triple edge.
Consider two cases: either $\Sigma$ contains ${\cal L}_5^5$ or it does
not (see Table~\ref{Lan} for the notation).

\bigskip
\noindent
{\bf Case 1.} Suppose that $\Sigma$ contains ${\cal L}_5^5$. Clearly,
only the open vertex of  ${\cal L}_5^5$ can be joined with some other 
nodes of $\Sigma$,
otherwise $\Sigma$ contains a parabolic subdiagram $\widetilde B_4$, 
or $\widetilde B_5$, or a small Lann\'er diagram.
By the similar reason, the diagram $\Sigma \setminus   {\cal L}_5^5$
is a linear diagram without double edges.
Since any node of $\Sigma$ belongs to some Lann\'er diagram,
for  $|\Sigma|>9$ such a diagram is superhyperbolic.

\bigskip
\noindent
{\bf Case 2.} Suppose that $\Sigma$ does not contain ${\cal L}_5^5$.

\begin{itemize}
\item
{\it Any triple edge of $\Sigma$ is incident to one of the leaves of 
$\Sigma$. The other vertex of the triple edge has valency two.
}

\noindent
{\it Proof:} the statement follows immediately from the absence of small Lann\'er
subdiagrams.

\item
{\it $\Sigma$ contains no double edge.}

\noindent
{\it Proof:} Suppose that $v_1v_2$ is a double edge of $\Sigma$.
Assume that $v_1$ is a leaf of $\Sigma$.
By Prop.~\ref{lanner_for_any_facet}, the node $v_1$ belongs to some Lann\'er diagram $L$.
Clearly, $L$ contains $v_2$ and $L$ is a linear Lann\'er diagram.  
Let $v$ be any node of $\Sigma\setminus L$  connected to $L$.
Then $\[L,v\]$ contains either parabolic or small Lann\'er subdiagram.

Therefore, neither $v_1$ nor $v_2$ is a leaf of $\Sigma$.
Then $\Sigma$ contains  edges $v_0v_1$ and $v_2v_3$
($v_0\ne v_3$, otherwise $\[v_1,v_2,v_3\]$ contains either a Lann\'er or
a parabolic subdiagram).
Since $\Sigma$ does not contain neither small Lann\'er diagrams nor
parabolic diagrams,
these edges are simple. Hence,  $\[v_0,v_1,v_2,v_3\]$
is a diagram of the type $F_4$.
If $v_4$ is any node of $\Sigma$ connected to $\[v_0,v_1,v_2,v_3\]$,
then $\[v_0,v_1,v_2,v_3,v_4\]$ either is a diagram of the type  
${\cal L}_5^5$ or it
contains  a small Lann\'er or a parabolic
subdiagram, which contradicts  the assumption.

\item
{\it $\Sigma$ is a tree.}
\noindent

{\it Proof:}
By the previous two statements, 
any minimal cycle in $\Sigma$ is a parabolic diagram $\widetilde A_n$.

\end{itemize}

Remove from $\Sigma$ all leaves that belong to triple edges.
Denote the obtained diagram by $\Sigma'$.
It follows from above that $\Sigma'$ is a simply laced (i.e. it
contained simple edges only) tree.
It follows from the classification of parabolic diagrams that
a simply laced tree without parabolic subdiagrams is an elliptic diagram.
Therefore, $\Sigma'$ is one of $ A_n, D_n, E_6, E_7, E_8$.

Now, append  triple edges to some of the leaves of $\Sigma'$.
If $\Sigma'= A_n $ or $D_{n+1} $ for $n\ge 9$, then
some vertices of $\Sigma$ belong to no Lann\'er diagram.
For the remaining cases we have either $|\Sigma|<10$,
or $\Sigma$ is one of the diagrams listed in the lemma, or
  $\Sigma$  is superhyperbolic.
A direct calculation of determinants shows that $\det \Sigma=0$ 
if and only if $\Sigma=\Theta_1$.

\end{proof}

\begin{cor}
\label{5}
Let $P$ be a compact Coxeter polytope.
Then $\Sigma(P)$ contains a small  Lann\'er diagram.

\end{cor}

\begin{proof}
Suppose that  $\Sigma(P)$ contains no small Lann\'er diagrams.
Then, by Lemma~\ref{disjoint}, $\Sigma(P)$ contains two disjoint
Lann\'er diagrams $L_1$ and $L_2$ of order 5, and hence, $|\Sigma(P)|\ge 10$.
A subdiagram $\[L_1, L_2\]$ is connected, otherwise it is superhyperbolic. 
It follows from Lemma~\ref{Table5}
that $\Sigma(P)$ is one of the diagrams 
listed in Fig.~\ref{10only5}.
In all three cases $\Sigma(P)$
contains a Lann\'er subdiagram which intersects any other Lann\'er subdiagram of $\Sigma(P)$.
By Lemma~\ref{disjoint}, $\Sigma(P)$ cannot be a Coxeter diagram
of a compact Coxeter polytope
and the statement is proved.

\end{proof}

\section{Liftings}
\label{section_lift}

\label{k-lifting}
Let $D$ be a diagram of missing faces
and $\Sigma$ be an admissible Coxeter diagram.
$\Sigma$ is a {\it $0$-lifting} of $D$ if
 there is a bijection  $\phi:D \to \Sigma $ such that
$M\in \M_D$ if and only if $\phi(M)$
is a Lann\'er diagram. $\phi$ is called  a {\it lifting bijection}.


$\Sigma$ is a {\it $k$-lifting} of $D$ ($k\in \N$)
if $\Sigma$
contains a subset $A$ of additional nodes satisfying
\begin{itemize}
\item[1)] $|A|=k$;
\item[2)] There exists an injection   $\phi:D \to \Sigma$ 
taking bijectively
$D$ to $\Sigma \setminus A$.
 $\phi$ is called a {\it lifting injection};
\item[3)] For any  Lann\'er diagram $L$ of $\Sigma$ the set
$\phi^{-1}(L\setminus A)$ contains a missing face of $D$;
\item[4)] For any $M\in \M_D$ there exists a Lann\'er diagram
$L\subset \[\phi(M), A\]$ containing $\phi(M)$.
\item[5)] For any set $\{a_1,\dots,a_r\}\subset A$ the diagram
$\Sigma\setminus \{a_1,\dots,a_r\}$ is not a $(k-r)$-lifting of $D$.

\end{itemize}

\bigskip

In~\cite{Ess2} a 0-lifting of $D$ is called a ``hyperbolic realization'' of $D$.
%
%
A $0$-lifting satisfies the general definition of a $k$-lifting.
When the number $k$ is not important, we write ``lifting'' instead of
``$k$-lifting''.
%

\bigskip

Let $\Sigma$ be a lifting of $D$.
Suppose that $L$ is a Lann\'er diagram of $\Sigma$ and 
for any 
$M\in \M_D$ we have $L\nsubseteq \[\phi(M),A\]$.
We say that $L$ is an {\it additional} Lann\'er diagram of 
the lifting $\Sigma$. 

Let $D$ be an abstract diagram of missing faces,  $\Sigma$ be an 
abstract Coxeter diagram and $\phi$ be an injection 
$\phi: D\to \Sigma$.
Suppose that there exists $\Sigma' \subset \Sigma$ 
which is a non-elliptic subdiagram not containing
$\phi(M)$ for any $M\in\M_D$. 
Then we say that $\Sigma'$ is a
{\it conflicting}  subdiagram (and conclude that $\Sigma$ is not a
lifting of $D$ with given injection).  


Notice that $\phi(M)$ is a subset of the set of nodes of
$\Sigma$. When we mean a Coxeter diagram spanned by
the nodes of  $\phi(M)$, we write
$\[\phi(M)\]$ (compare the notation introduced in Section~\ref{cox_def}).


\begin{lemma}
\label{sublifting}
Let $D$ be a diagram of missing faces and $D_1\subset D$ be a subdiagram.
Let $\Sigma$ be a $k$-lifting of $D$.
Then $\Sigma$ contains a subdiagram $\Sigma_1$ which is a $k_1$-lifting of $D_1$
for some $k_1\le k$.

\end{lemma}

\begin{proof}
Let $A$ be the set of additional nodes of $\Sigma$.
Take  $\phi|_{D_1}$
for the  lifting injection of $D_1$. 
By definition of liftings, $\[\phi(D_1),A\]$ contains some $k_1$-lifting of
$D_1$, and we have $k_1\le k$.

\end{proof}

The proof of the following lemma is based on the multiple usage of
Proposition~\ref{face}.

\begin{lemma}
\label{lift}
Let $P\subset \H^d$ be a simple
Coxeter polytope and let $f$ be an $m$-face of $P$.
Let $D(f)$ be a diagram of missing faces  of $f$.
Then $\Sigma(P)$ contains a subdiagram $\Sigma_0$ which is a
$k$-lifting of $D(f)$ for some $k\le d-m$.

\end{lemma}

\begin{proof}
Let $f_1,...,f_{d-m}$ be the facets of $P$ containing $f$.
Let $F_1=f_1$ and define inductively  $F_i=F_{i-1}\cap f_i$.
Clearly, $F_i$ is a facet of $F_{i-1}$ and $F_{d-m}=f$.

If $m=d$, the lemma is trivial.
Assume by induction that the statement holds for any face of dimension greater than
$m$.
Let us prove for the $m$-face $f$.
By the induction assumption, $\Sigma(P)$ contains a subdiagram 
$\Sigma_{r-1}$
which is a $k$-lifting of $D(F_{d-m-1})$, where $k\le d-(m+1)$.
The following claim completes the proof.

\noindent
{\bf Claim}: either $\Sigma_{r-1}$ contains a $k$-lifting of $f$  
or $\[\Sigma_{r-1},v_{d-m}\]$ contains a $(k+1)$-lifting
of $f$, where $v_{d-m}$ is a node of $\Sigma(P)$ corresponding to the
facet $f_{d-m}$. \\

\noindent
Proof of the claim:

\begin{itemize}
\item[]

Let $\Pi_1,...,\Pi_s$ be the facets of $F_{d-m-1}$.
Denote by $J$ a set of indices such that
 $\Pi_i\cap f_{d-m}\ne \emptyset$ for $i\in J$, 
and denote  $\pi_i=\Pi_i\cap f_{d-m}$ for $i\in J$.
Then $\{ \pi_i\}$ is the set of facets of $f=F_{d-m}$.

Let $\phi_1$ be the lifting injection for $F_{d-m-1}$,
in particular $\phi_1$ takes $\{\Pi_i \}$ to $\Sigma_{r-1}\setminus A$,
where $A$ is the set of additional nodes.
Denote by $\psi$ the map $\pi_i\to \Pi_i$.
Denote $\phi=\phi_1\circ \psi$.
Then $\phi$ is an injection from $\{\pi_i\}$ to $\Sigma\setminus A$. 
Let $\{\Pi_i\,|\,i\in I \}$ be a missing face of $D(F_{d-m})$
(where $I\subset J$ is some index set). 
By Prop.~\ref{face}, either $\{\pi_i\,|\,i\in I\}$ or 
$\{f_{d-m}\}\cup \{\pi_i\,|\,i\in I\}$ is a missing face of
$D(F_{d-m-1})$.
This proves property 4) of the definition of liftings
for either $\Sigma_{r-1}$ or $\[\Sigma_{r-1},v_{d-m}\]$.
By the same proposition,
if $\{\Pi_i\,|\,i\in K\}$ is not a missing face in $D(F_{d-m})$
then neither $\{\pi_i\,|\,i\in K\}$ nor 
$\{f_{d-m}\}\cup\{\pi_i\,|\,i\in K\} $
is a missing face in $D(F_{d-m-1})$.
This proves property 3).

Thus, either $A$ or $A\cup \{f_{d-m}\}$ contains a set of 
additional nodes
for some $k$-lifting of $D(F_{d-m-1})$.
Hence, $k\le d-(m+1)+1=d-m$
and everything is proved.

\end{itemize}

\end{proof}

\begin{cor}
\label{simplex}
Let $P\subset \H^d$ be a compact Coxeter polytope and $f$ be a face of $P$.
Then $D(f)$ contains no large missing faces.
In particular, if $\dim f >4 $ then $f$ is not a simplex.

\end{cor}

\begin{proof}
By the definition of liftings,  
for any missing face $M$ of $D$ there exists a Lann\'er diagram 
$L$ in $\Sigma$
such that $|L|\ge |M|$.
Hence, the statement follows from the fact that a Lann\'er diagram 
contains at most 
5 nodes.
%

\end{proof}

\begin{lemma}
\label{cvetok}
Let $D$ be a diagram of missing faces such that $|M|=5$
for any $M\in \M_D$.
Then any lifting of $D$ is a $0$-lifting.

\end{lemma}

\begin{proof}
Suppose that  $D$ is a $k$-lifting, where $k>0$,
and $\phi$ is a lifting injection.
Remove from $D$  all additional nodes and denote by $D_1$ 
the diagram obtained. Then $D_1$ is a $0$-lifting, 
and $\phi$ is a lifting bijection. Indeed, the conditions
2),3) and 5) of the definition of $k$-lifting hold evidently.

For the condition 4), consider any missing face
$M\in \M_D$.
By definition, $\phi(M)$ belongs to some Lann\'er diagram of
$\Sigma$.
Note, that $\[\phi(M)\]$ is a subdiagram of order 5.
Since the order of a Lann\'er diagram does not exceed 5,
$\[\phi(M)\]$ is a Lann\'er diagram.

%

\end{proof}

Now we will prove several lemmas about liftings
we will use later.\\

\noindent
{\bf Notation.} 
\vspace{-7pt}
\begin{itemize}
\item
Let $D$ be a diagram of missing faces and $N_1,...N_r\subset D$.
We write $D=\bigcup \limits_{i=1}^r N_i$ if for any node of $D$ 
there exists a set $N_i, i\in\{ 1,\dots,r\}$ containing this node.

\item
Let $D=\bigcup \limits_{i=1}^r N_i$ be a diagram of
missing faces, where $N_i\cap N_j =\emptyset$ for any $i\ne j$.
Suppose that there exists $k\in \{1,\dots,r\}$ such that
 $M\in \M_D$ if and only if
$M=\bigcup \limits_{t=i}^{i+k-1} N_t$ for some $i\in\{1,\dots,r-k+1\}$.
Then we write $D=\lf N_1,N_2,\dots,N_r \rf_k$.

\item
When we are interested in combinatorial type of $D$ only rather than in 
concrete subdiagrams $N_i$,
we write $\lf |N_1|,\dots,|N_r| \rf_k$.
For example, $\lf 1,4,1,3 \rf_2$ stays for the diagram

\begin{center}
\epsfig{file=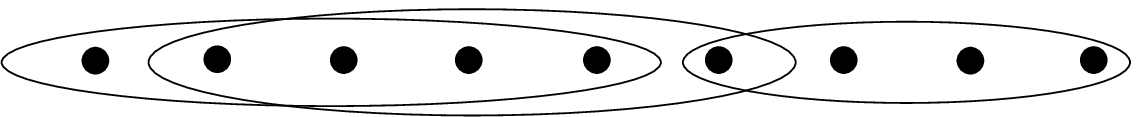,width=0.28\linewidth}
\end{center}

\item
If  there exists $k\in \{1,\dots,r\}$ such that
$M\in \M_D$ if and only if
$M=\bigcup \limits_{t=i}^{i+k-1} N_t$ for some $i\in{1,\dots,r}$,
where $t$ is considered as a number modulo $ r$,
we write  $D=(N_1,N_2,\dots,N_r)_k$
or  $(|N_1|,\dots,|N_r|)_k$.

\item
We write $\Sigma \approx \lf N_1,N_2,\dots,N_r \rf_k$,
if $\Sigma$ contains no parabolic subdiagrams and 
the structure of Lann\'er subdiagrams of $\Sigma$ corresponds
to the diagram of missing faces  $ \lf N_1,N_2,\dots,N_r \rf_k$,
i.e.  $\Sigma$ consists of disjoint subdiagrams
$L_1,\dots,L_r$, $|L_i|=|N_i|$, 
and  $L\subset \Sigma$ is a Lann\'er subdiagram if and only if
$L=\[L_i,L_{i+1}, \dots, L_{i+k-1}\]$
for some
$i\in\{1,\dots,r-k+1\}$.

Similarly we use the notation
$\Sigma \approx \lf |N_1|,\dots,|N_r| \rf_k$, 
$\Sigma \approx  (N_1,N_2,\dots,N_r)_k$
and  $\Sigma \approx (|N_1|,\dots,|N_r|)_k$.

\end{itemize}

\begin{prop}[\cite{Ess2}, Lemma~4.7]
\label{spletayushaya}
Let $D=\lf M,N \rf_1$ be a diagram of missing faces, where
$|M|\ge 3$ and $|N|\ge 4$.
Then $D$ has no $0$-liftings.

\end{prop}

\begin{lemma}
\label{lift54}
The diagram  $\Theta_1$ shown in Fig.~\ref{10only5}
is the only lifting of $D=\lf 4,5 \rf_1$.
The additional node of this lifting is the leaf of $\Theta_1$
which does not belong to a triple edge.

\end{lemma}

\begin{proof}
Let $\Sigma$ be a lifting of $D=\lf M_4,M_5 \rf_1$,
where $|M_5|=5$ and $|M_4|=4$, and let $\phi$ be a lifting injection.
Suppose that $\[\phi(M_4)\]$ is a Lann\'er diagram of $\Sigma$.
Then by Prop.~\ref{spletayushaya} 
$\Sigma\approx \lf \phi(M_4),\phi(M_5) \rf_1$ is a superhyperbolic diagram.
Hence,  $\[\phi(M_4)\]$ is an elliptic diagram and
$\Sigma$ contains at least one additional node $a$.
Then the diagram  $\[a,\phi(M_4),\phi(M_5)\]$ satisfies the properties
1)-4) of the definition of liftings. So, by property 5) $a$ is a
unique additional node of $\Sigma$, and $\Sigma$ is a 1-lifting.
By property 3) of the definition
of liftings,
any Lann\'er subdiagram of\ $\Sigma$ consists of 5 nodes.
Since $|\Sigma|=|\[a,\phi(M_4),\phi(M_5)\]|=10$, 
 Lemma~\ref{Table5} implies that
$\Sigma$ is one of the diagrams $\Theta_1$, $\Theta_2$ and $\Theta_3$
shown in Fig.~\ref{10only5}.
Note that $\Sigma$ contains two disjoint Lann\'er diagrams 
$\[\phi(M_4),a\]$ and $\[\phi(M_5)\]$. 
Hence, $\Sigma=\Theta_1$ or $\Theta_2$.

It is easy to find a unique lifting injection for $\Theta_1$:
$a$ is the only leaf of $\Sigma$ which does not belong to a triple edge,
$\phi(M_4)$ are the remaining points in the lower row of vertices 
(see Fig.~\ref{10only5}).

Lifting injection for $\Theta_2$ does not exist:
$\phi(M_4)$ belongs to any Lann\'er diagram
different from $\phi(M_5)$, on the other hand,
$\Theta_2$ contains no such a quadruple of vertices,
so the lemma is proved.

\end{proof}

\begin{lemma}
\label{flower}
Let $D=N\cup \{x_1,x_2,x_3\}$ be a diagram of missing faces,
where $|N|=4$. 
If $\M_D=\{N\cup x_1,N\cup x_2,N\cup x_3 \}$
then $D$ has no liftings.

\end{lemma}


\begin{proof}
Suppose that $\Sigma$ is a lifting of $D$ and $\phi$ is a lifting injection.
By Lemma~\ref{cvetok}, $\Sigma$ is a 0-lifting.
Now, consider $\Sigma_{ij}=\[\phi(N\cup \{x_i,x_j\})\]$, $i\ne j$.   
Clearly, $\Sigma_{ij}\approx \lf 1,4,1 \rf_2$. 
By~\cite[Lemma~5.3]{Ess2}, $\lf 1,4,1 \rf_2$
is one of the following diagrams:

\begin{center}
\epsfig{file=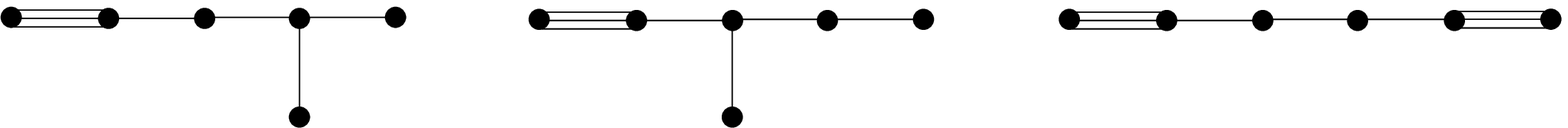,width=0.85\linewidth}
\end{center}

Taking each of these diagrams for $\Sigma_{12}$ and trying to add $x_3$
in order to compose a correct diagram $\Sigma_{13}$, we obtain
a conflicting subdiagram.

\end{proof}

\begin{figure}[htb]
\begin{center}
\epsfig{file=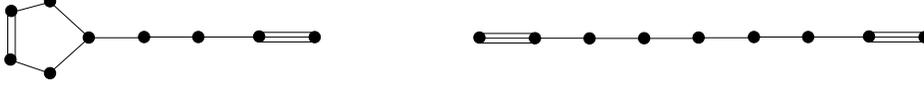,width=0.85\linewidth}
\caption{1-liftings of $\lf 4,4 \rf_1$. }
\label{fig414}
\end{center}
\end{figure}

\begin{lemma}
\label{lift44}
%
%
The diagram $D=\lf 4,4 \rf_1$ has no $0$-liftings, and any $1$-lifting 
of $D$ is one of the two diagrams  shown in Fig.\ref{fig414}.

\end{lemma}

\begin{proof}
By Prop.~\ref{spletayushaya}, there are no 0-liftings of $D$. 
Let $\Sigma$ be a 1-lifting of $D$, let $\phi$ be a lifting injection,
and let $a$ be the additional node.

Denote by $M_1$ and $M_2$ the missing faces of $D$.
We may assume that $\[\phi(M_1)\]$ is an elliptic subdiagram of $\Sigma$
and $\[\phi(M_1),a\]$ is a Lann\'er diagram.
We consider two cases:

\begin{itemize}
\item[Case 1.]
Suppose that $\[\phi(M_2)\]$ is a Lann\'er diagram.
By Prop.~\ref{spletayushaya}, $\Sigma$ contains an additional 
 Lann\'er diagram $L=\[\phi(M_1),x\]$, where $x\in
\phi(M_2)$. By Lemma~\ref{flower}, this additional Lann\'er
diagram $L$ is unique.
Thus, 
$\Sigma\approx \lf a,\[\phi(M_1)\],x,\[\phi(M_2)\setminus x\] \rf_2
= \lf 1,4,1,3 \rf_2$. 
This is impossible by~\cite[Folgerung~5.10]{Ess2}.

\item[Case 2.]
Suppose that both $\[ \phi(M_1)\]$ and $\[\phi(M_2)\]$ are elliptic.
Then $\[\phi(M_1),a\]$ and $\[\phi(M_2),a\]$ are Lann\'er diagrams. 
Note that $a$ is an open vertex of $\[\phi(M_1),a\]$,
otherwise $\Sigma$ has a conflicting Lann\'er diagram.
Similarly, $a$ is an open vertex of $\[\phi(M_2),a\]$.
Hence, each of  $\[\phi(M_1),a\]$ and   $\[\phi(M_2),a\]$
is one of ${\cal L}_1^5$ and ${\cal L}_5^5$.

If $\Sigma$ contains an edge joining 
 $\[\phi(M_1)\]$ with  $\[\phi(M_2)\]$
then $\Sigma$ has a conflicting subdiagram.
Therefore, 
 $\[\phi(M_1)\]$ and  $\[\phi(M_2)\]$
are not connected in $\Sigma$,
and $\Sigma$ is one of the two diagrams  shown in Fig.~\ref{fig414}
(if $\[\phi(M_1),a\]=\[\phi(M_2),a\]={\cal L}_5^5$
then $\Sigma$ contains a parabolic subdiagram).


\end{itemize}

\end{proof}

\begin{figure}[!h] 
\begin{center}
\epsfig{file=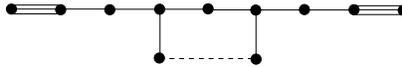,width=0.37\linewidth}
\caption{A unique compact hyperbolic Coxeter 8-polytope 
$P_8$ with 11 facets.}
\label{p8}  
\end{center}
\end{figure}

Denote by $\!P_8$ a unique compact hyperbolic Coxeter 8-polytope with
11 facets (see~\cite{Bu1},\cite{Ess2}), 
and denote by $\Sigma(P_8)$ its Coxeter diagram,
see Fig.\ref{p8}.

\begin{lemma}
\label{442}
Let $D=\[4,4,2\]_1$. 
Let $\Sigma$ be a $0$- or $1$-lifting of $D$
containing no Lann\'er diagram of order $3$.
Suppose that the positive inertia index of  $\Sigma$ does not exceed $8$.
Then $\Sigma=\Sigma(P_8)$.

\end{lemma}


\begin{proof}
By Lemmas~\ref{sublifting} and~\ref{lift44}, $D$ has no 0-liftings.
Let $\Sigma$ be a 1-lifting of $D$, let $\phi $ be a lifting
injection,
and let $a$ be the additional node.
Denote by $M_1, M_2$ and $M_3$  the missing faces of $D$,
where $|M_1|=|M_2|=4$ and $|M_3|=2$. Let $M_3=\{u_1,u_2\}$. 
By Lemma~\ref{lift44}, $\[\phi(M_1),\phi(M_2),a\]$
is one of the two diagrams  shown in Fig.\ref{fig414}.
Since $\Sigma$ has only one additional node,
the absence of  Lann\'er diagram of order 3 implies that
$\[\phi(M_3)\]$ is a Lann\'er diagram.
Therefore, any Lann\'er subdiagram of $\Sigma$
distinct from $\[\phi(M_3)\]$ is a Lann\'er diagram of order 5.

Suppose that the node  $u_1$ does not belong to any additional Lann\'er subdiagram
of $\Sigma$. Then $\[\phi(M_1),\phi(M_2),u_1\]$
is an elliptic subdiagram of $\Sigma$ of order 9 in contradiction to 
the assumption that the positive inertia index of $\Sigma$ 
is at most 8.

Hence, both $u_1$ and $u_2$  belong to some additional Lann\'er diagrams.
Consider $X_1=\[u_1,\phi(M_1),\phi(M_2),a\]$. 
$X_1$ is a subdiagram of order
10 containing no small Lann\'er diagrams. Moreover,
each vertex of $X_1$ belongs to some Lann\'er diagram.
By Lemma~\ref{Table5},  $X_1$ is one of the diagrams $\Theta_1$, $\Theta_2$
and $\Theta_3$ shown in Fig.~\ref{10only5}.
By the assumption of the lemma, the positive inertia index of $\Sigma$ 
is at most 8. Therefore, $\det(X_1)=0$ and 
$X_1=\Theta_1$ (see Lemma~\ref{Table5}).
By Lemma~\ref{lift44}, $\[\phi(M_1),\phi(M_2),a\]$ is a linear
subdiagram of $\Theta_1$.
Hence, $u_1$ is the only leaf of $\Theta_1$ that is not incident to  
a triple edge.
Similarly, for  $X_2=\[u_2,\phi(M_1),\phi(M_2),a\]$ 
we have  $X_2=\Theta_1$.
Moreover, $u_1$ and $u_2$ are joined with different vertices of
$\[\phi(M_1),\phi(M_2),a\]$, otherwise one of the subdiagrams
$\[\phi(M_i),a,u_1,u_2\]$ is superhyperbolic.
Hence, $\Sigma=\Sigma(P_8)$.

\end{proof}

\begin{lemma}
\label{14131}
Let $D=(1,4,1,3,1)_2$.
Let $\Sigma$ be a $0$- or $1$-lifting of $D$
containing no Lann\'er diagram of order $3$.
Suppose that the positive inertia index of  $\Sigma$ does not exceed $8$.
Then $\Sigma=\Sigma(P_8)$.

\end{lemma}

\begin{proof}
Let $\phi$ be a lifting injection.
Denote the subsets of $D$ as follows:
$(1,4,1,3,1)_2=(v,N_4,u,N_3,w)_2$.


Suppose $\[\phi(N_3),\phi(u)\]$ is a Lann\'er diagram.
Then 
$\[\phi(v),\phi(N_4),\phi(u), \phi(N_3)\]\approx \lf 1,4,1,3 \rf_2$, 
which is impossible 
by~\cite[Folgerung~5.10]{Ess2}.  
Therefore, $\[\phi(N_3),\phi(u)\]$ is an elliptic diagram
and $\Sigma$ is not a 0-lifting of $D$.

Now, assume that $\Sigma$ is a 1-lifting and $a$ is a unique
additional node of $\Sigma$.
Then the diagram
$X=\[a,\phi(D\setminus \{w\})\]$ satisfies the conditions of
Lemma~\ref{Table5}.
Since $|X|=10$, $X$ is one of the diagrams $\Theta_1$, $\Theta_2$
and $\Theta_3$ shown in Fig.~\ref{10only5}.
By the assumption, the positive inertia index of $\Sigma$ 
is at most 8. Hence, $\det(X)=0$
and $X=\Theta_1$: 

\begin{center}
\psfrag{a}{$\[a,\phi(N_3)\]$}
\psfrag{b}{$\[\phi(N_4)\]$}
\psfrag{u}{$\phi(u)$}
\psfrag{v}{$\phi(v)$}
\psfrag{1}{\small $x_1$}
\psfrag{2}{\small $x_2$}
\psfrag{3}{\small $x_3$}
\psfrag{4}{\small $x_4$}
\psfrag{5}{\small $x_5$}
\psfrag{6}{\small $x_6$}
\psfrag{7}{\small $x_7$}
\psfrag{8}{\small $x_8$}
\psfrag{9}{\small $x_9$}
\psfrag{X}{\Large $X=$}
\epsfig{file=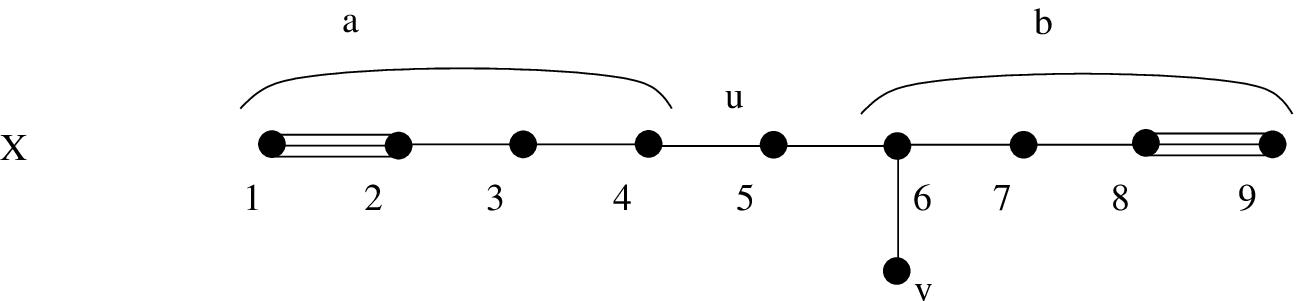,width=0.5\linewidth}
\end{center} 

Consider now $\[X,\phi(w)\]$.
Since $\[\phi(N_4),\phi(w)\]$ should be an elliptic diagram,
 $\phi(w)$  does not attach to $\[x_6,x_7,x_8,x_9\]=\[\phi(N_4)\]$.
Further, $\phi(w)$ is not connected to  $\[x_1,x_2\]$,
otherwise $\[x_1,x_2,x_3,\phi(w)\]$ is a conflicting subdiagram.
So, besides $\phi(v)$, $\phi(w)$ may be attached to $x_3$, $x_4$, or
$x_5$ only.
Any edge joining $w$ with $x_3,x_4,x_5$ is a simple edge,
otherwise one of the diagrams  $\[x_2,x_3,x_4,\phi(w)\]$,
$\[x_3,x_4,x_5,\phi(w)\]$
and $\[x_4,x_5,x_6,\phi(w)\]$ is a conflicting subdiagram. 
In order to obtain no parabolic subdiagrams,
$\phi(w)$ should be  connected to one of the vertices
$x_3,x_4,x_5$ only (except $\phi(v)$).

If $\phi(w)$ attaches to $x_3$ then the subdiagram 
$\Sigma \setminus \phi(v)$ is superhyperbolic.
If $\phi(w)$ is connected to $x_5$ then $\phi(N_3\cup\{w\})$ 
is not contained
in any Lann\'er diagram in contradiction to the definition of liftings.
Thus, $\phi(w)$ attaches to $x_4$. By assumption $\Sigma$ contains no
Lann\'er diagrams of order 3,   so $\[\phi(v),\phi(w)\]$ is a Lann\'er
diagram and
$\Sigma=\Sigma(P_8)$. 

\end{proof}

\begin{lemma}
\label{1323}
The diagram $D=\lf 1,3,2,3 \rf_2$ has no $0$-liftings.

\end{lemma}

\begin{proof}
Suppose that $\Sigma$ is a 0-lifting of 
$D= \lf 1,3,2,3 \rf_2$.
Comparing the 0-liftings of the diagram $\lf 3,2,3 \rf_2$ 
(see~\cite[Lemma~5.12]{Ess2})
with the 0-liftings of
the diagrams $\lf 2,3,1 \rf_2$
(see~\cite[Table~4, item~8]{n3}), we obtain that 
$\Sigma$ is one of the following diagrams:

\begin{center}
\epsfig{file=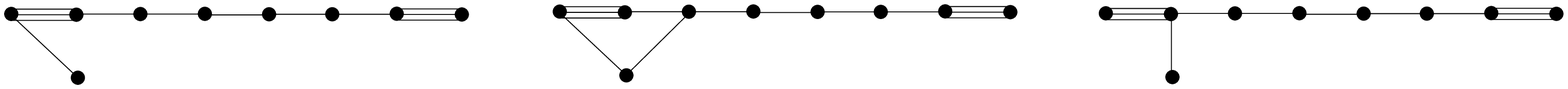,width=0.97\linewidth}
\end{center} 

All of these diagrams  are superhyperbolic
and the lemma is proved.

\end{proof}

\begin{lemma}
\label{32311}
Let $D=(3,2,3,1,1)_2$. 
There is no $0$- or $1$-lifting of $D$ having the positive inertia index
at most $8$ and containing no Lann\'er subdiagram of order $3$.

\end{lemma}

\begin{proof}
Suppose that $\Sigma$ is either  0- or 1-lifting of $D$.
Let $\phi$ be the lifting injection.
Denote the subsets of $D$ as follows:
$(3,2,3,1,1)_2=(N_1,N_2,N_3,v,u)_2$.
Denote $X=\[\phi(N_1), \phi(N_2), \phi(N_3)\]$.

The subdiagram $X$ of $\[\phi(D)\]$ contains no 
small Lann\'er diagrams.
Hence, $X \approx [3,2,3]_2$, 
and by Lemma~5.12 of~\cite{Ess2} we have

\begin{center}
\psfrag{N1}{$\[\phi(N_1)\]$}
\psfrag{N2}{$\[\phi(N_2)\]$}
\psfrag{N3}{$\[\phi(N_3)\]$}
\psfrag{1}{\small $x_1$}
\psfrag{2}{\small $x_2$}
\psfrag{3}{\small $x_3$}
\psfrag{4}{\small $x_4$}
\psfrag{5}{\small $x_5$}
\psfrag{6}{\small $x_6$}
\psfrag{7}{\small $x_7$}
\psfrag{8}{\small $x_8$}
\psfrag{X}{\Large $X=$}
\epsfig{file=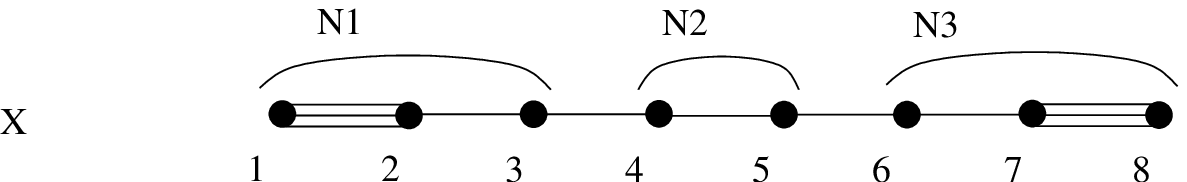,width=0.5\linewidth}
\end{center}

By the  definition of liftings, the diagram
$\[\phi(N_1),\phi(u)\]$ belongs to some
Lann\'er diagram of $\Sigma$.
If $\[\phi(N_1),\phi(u)\]$ is a Lann\'er diagram,
we have $\[X,\phi(u)\]\approx 
\lf \phi(u),\phi(N_1),\phi(N_2),\phi(N_3) \rf_2=\lf 1,3,2,3 \rf_2$,
which is impossible by Lemma~\ref{1323}.


Thus, $\[\phi(N_1),\phi(u)\]$ is an elliptic diagram
and $\Sigma$ contains a unique additional node $a$.
Similarly, $\[\phi(N_3),\phi(v)\]$ is an elliptic diagram. 
Since $\Sigma$ contains no Lann\'er diagram of order 3,
$\[\phi(u), \phi(v)\]$ is a Lann\'er diagram.
The remaining small Lann\'er diagrams are
$\[\phi(N_1),\phi(u)\]$ and  $\[\phi(N_3),\phi(v)\]$.
Therefore, any additional Lann\'er diagram of $\Sigma$ contains either
$\[\phi(N_1),\phi(u)\]$ or  $\[\phi(N_3),\phi(v)\]$,
and the diagram $Y=\[X,a,\phi(u)\]$ satisfies the conditions of
Lemma~\ref{Table5}. Since $|Y|=10$,
$|Y|$ is one of the diagrams $\Theta_1,\Theta_2$ and
$\Theta_3$ shown in Fig.~\ref{10only5}.
The diagram $\Theta_1$ contains no subdiagram 
$\Sigma_0\approx \lf 3,2,3 \rf_2$.
The diagrams $\Theta_2$ and $\Theta_3$ have positive inertia index 9
in contradiction to the assumption  of the lemma.

\end{proof}

\begin{lemma}
\label{1lift53}
The diagrams $\lf 5,3 \rf_1$ and $\lf 5,3,2 \rf_1$  have neither 
$0$-lifting nor $1$-lifting.

\end{lemma}

\begin{proof}
We will prove the Lemma for $D=\lf 5,3 \rf_1$,
the statement for $\lf 5,3,2 \rf_1$ follows immediately.

By Prop.~\ref{spletayushaya}, $D$ has no 0-lifting.
Suppose that $\Sigma$ is a 1-lifting of $D$, $\phi$ is a lifting
injection and $a$ is the additional node.  
Denote by $M_5$ and $M_3$ the missing faces of $D$ 
of orders 5 and 3 respectively.
Then $\[\phi(M_3),a\]$ is a Lann\'er diagram.
By Prop.~\ref{spletayushaya}, $\Sigma$ has at least one additional 
Lann\'er diagram. By the definition of liftings, any additional Lann\'er
diagram
of $\Sigma$ contains $\phi(M_3)$.

Suppose that $L=\[\phi(M_3),x\]$ is an additional Lann\'er diagram of
$\Sigma$ of order 4, $x\in \phi(M_5)$.
Then $x$ is an open vertex of $L$, otherwise $\Sigma$ contains a
conflicting subdiagram. Moreover, it is easy to see that
$x$ is a doubly open vertex of $L$,
that implies $L={\cal L}^4_5$.
On the other hand, $\[L,a\]\approx \lf x,M_3,a \rf_2=\lf 1,3,1 \rf_2$.
By~\cite[Table~4, item~8]{n3}, no diagram $\Sigma_0$ such that
$\Sigma_0\approx \lf 1,3,1 \rf_2$ may contain ${\cal L}^4_5$.
This means that $\Sigma$ contains no additional Lann\'er diagram of order 4.

Therefore, $\Sigma$ contains at least one additional Lann\'er diagram 
of order 5, denote it by
$L=\[\phi(M_3),x_1,x_2\]$, $x_1, x_2\in \phi(M_5)$.
Since $L$ is connected, we may assume that $x_1$ is joined with
some $y\in \phi(M_3).$
Then $x_1$ is an open vertex of $\[\phi(M_5)\]$, and 
$y$ is an open vertex of $\[\phi(M_3),a\]$.
Since none of Lann\'er diagrams of order 5 has more than one open vertex,
we conclude that $x_2$ does not attach to $\[\phi(M_3)\]$.
Hence, $x_2$ is joined with $x_1$.
By the same reason, any other additional Lann\'er diagram $L'$ 
contains $x_1$, and if $x_k\in L'$, $x_k\ne x_1$, then $x_k$ attaches to
$x_1$ and does not attach to $\[\phi(M_3)\]$.

If $L$ is the only additional Lann\'er diagram of $\Sigma$
then 
\begin{center}
$\Sigma\approx
\lf \[\phi(M_5)\setminus \{x_1,x_2\}\],\[x_1,x_2\],\[\phi(M_3)\],a \rf_2=
 \lf 3,2,3,1 \rf_2$, 
\end{center}
which is impossible by Lemma~\ref{1323}.
Hence, $\Sigma$ contains at least one more additional Lann\'er diagram
$L'$. 
Suppose that $L'=\[\phi(M_3),x_1,x_3\]$, $x_3\in \phi(M_5)$,
$x_3\ne x_2$. 
As it was shown above, the open vertex $x_1$ of $\phi(M_5)$
is connected with two other
nodes, $x_2$ and $x_3$,
so $\[\phi(M_5)\]={\cal L}_5^5$. 
In particular, the edges $x_1x_2$ and $x_1x_3$ are simple.
Now, taking in mind that 
$\[L,a\] \approx \lf a,\[\phi(M_3)\],\[x_1,x_2\] \rf_2=\lf 1,3,2 \rf_2$
and
$\[L',a\] \approx \lf a,\[\phi(M_3)\],\[x_1,x_3\] \rf_2=\lf 1,3,2 \rf_2$,
we check all possibilities for $\lf 1,3,2 \rf_2$ such that 
$x_1x_2$ and $x_1x_3$
are simple edges (see~\cite[Table~4, item~8]{n3}). 
In this way we obtain that either $\Sigma$ contains some conflicting
Lann\'er diagram or $\Sigma$ is one of the following diagrams:

\begin{center}
\psfrag{1}{$x_1$}
\psfrag{2}{$x_2$}
\psfrag{3}{$x_3$}
\psfrag{a}{$a$}
\psfrag{y}{$y$}
\epsfig{file=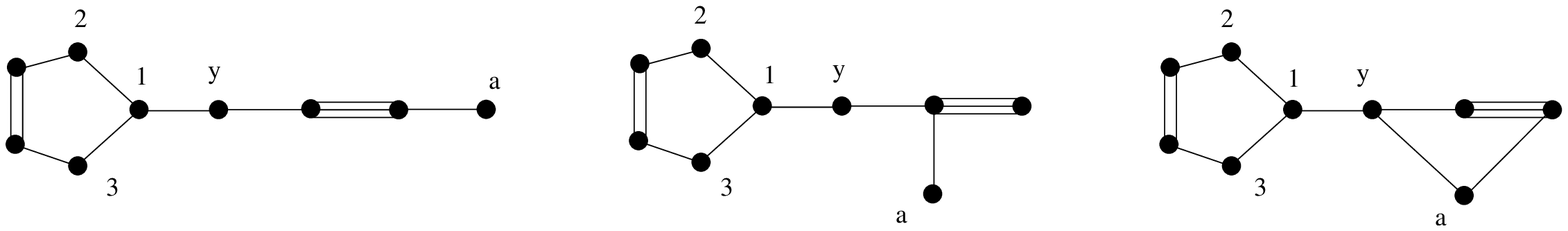,width=0.8\linewidth}
\end{center}

All these diagrams are superhyperbolic.
We come to a contradiction which proves the lemma.

\end{proof}

%
%
%
%

\section{Absence of polytopes in dimensions $\ge$9}
\label{new}

Suppose that there exists a compact hyperbolic Coxeter $d$-polytope
$P$ with $d+4$ facets. Denote by $\Sigma$ its Coxeter diagram.

We will refer to the classification of compact hyperbolic Coxeter 
$d$-polytopes with $d+k$ facets for $k\le 3$ 
(see~\cite{L},~\cite{Ess_eng},~\cite{K} and~\cite{n3}).
In particular, we remind that $P_8$ is a unique Coxeter 8-polytope with 11
facets, and $\Sigma(P_8)$ is its Coxeter diagram (see Fig.~\ref{p8}).

From now on by polytope we mean a compact hyperbolic Coxeter polytope.
By ``a diagram of a polytope'' we  mean  a Coxeter diagram of the polytope.

\begin{lemma}
\label{nok}
Suppose that $d\ge 9$. Then \\
(1) any node of $\Sigma$ is incident to at most one non-simple edge;\\
(2) $\Sigma$ contains no Lann\'er diagram of order $3$;\\
(3) $\Sigma$ contains no edge of multiplicity $\ge 4$.

\end{lemma}

\begin{proof}
Suppose the lemma is broken. Then there exists either a node incident
to two dotted edges, or a subdiagram
$S_0\subset \Sigma$ of the type $G_2^m$ for some $m\ge 4$ having at least
one bad neighbor. 

Suppose that  there exists a node $v$ incident
to two dotted edges.
Then $v$ corresponds to a
$(d-1)$-dimensional (non-Coxeter) face $f$ with at most 
$d+1=(d-1)+2$ facets,
so $f$ is either a simplex or a product of two simplices.
If $f$ is a simplex then $P$
has a large missing face, which is impossible by 
Cor.~\ref{simplex}. If $f$ is a product of two simplices and it has no
large missing faces then $f$  is a product of two 4-dimensional
simplices and $D(f)=\lf 5,5 \rf_1$ (since $d\ge 9$). 
By Prop.~\ref{spletayushaya} 
the diagram $\lf 5,5 \rf_1$ has no liftings, so we contradict   
Lemma~\ref{lift}.

Suppose that $\Sigma$ contains
 a subdiagram
$S_0\subset \Sigma$ of the type $G_2^m$ for some $m\ge 4$ having at least
one bad neighbor.
Then $P(S_0)$ is a Coxeter $(d-2)$-polytope with
at most $(d-2)+3$ facets, so $(d-2)$ is either equal to 8 or less
than 7 (see~\cite{n3}). If $d-2<7$ we obtain that $d\le 8$. 
If $d-2=8$ we have $d=10$, and $P(S_0)$ is a unique Coxeter 8-polytope 
with 11 facets (see Fig.~\ref{p8}).
The diagram $\Sigma_{S_0}$ 
contains a subdiagram of the type $H_4$ with 2
neighbors (the neighbors are attached to this diagram by simple edges). 
Cor.~\ref{dif2} implies that this diagram corresponds to
a subdiagram  
$S_1\subset \Sigma$ of the type $H_4$ with 2 neighbors.
So, $P(S_1)$ is a
Coxeter  6-polytope with at most 8 facets, which is impossible
(see~\cite{Ess_eng},~\cite{K}).     

\end{proof}

\begin{lemma}
\label{no5-9}
Suppose that $d\ge 9$. Then $\Sigma$ does not contain any Lann\'er diagram
of order $5$. 

\end{lemma}

\begin{proof}
Suppose there exists a diagram $L_0\subset \Sigma$ of order 5. 
Clearly, it
contains a subdiagram $S_0\subset L_0$, such that $S_0=H_4$ or
$F_4$. Then $P(S_0)$ is a $(d-4)$-dimensional Coxeter polytope, and its 
Coxeter diagram $\Sigma_{S_0}$ is a subdiagram of $\Sigma$ 
(Cor.~\ref{dif}).
The diagram $S_0$ has at least one neighbor, the node
$L_0\setminus S_0$. Thus, $P(S_0)$ has either $d-3$, or $d-2$, or $d-1$
facets. 

\bigskip

\noindent
{\bf Case 1:} $P(S_0)$ has $d-3$ facets. Then  $P(S_0)$ is a 
Coxeter simplex, and we obtain $d-4\le 4$, i.e. $d\le 8$.

\bigskip

\noindent
{\bf Case 2:} $P(S_0)$ has $d-2$ facets. Then $P(S_0)$ is a 
Coxeter $(d-4)$-polytope with $(d-4)+2$
facets. Results of~\cite{Ess2} and~\cite{K} imply that $d-4\le 5$,
i.e. $d\le 9$. If $d=9$, $P(S_0)$ is a 5-dimensional prism.
Since $\Sigma_{S_0}=\o S_0$,
condition (1)  of Lemma~\ref{nok} and~\cite{K} imply that there are only
four possibilities for $\Sigma_{S_0}$: all of
them are presented on Fig.~\ref{prisms}.

\begin{figure}[!h] 
\begin{center}
\psfrag{a}{a)}
\psfrag{b}{b)}
\psfrag{c}{c)}
\psfrag{d}{d)}
\psfrag{u}{$u$}
\psfrag{u1}{$u_1$}
\psfrag{u2}{$u_2$}
\psfrag{u3}{$u_3$}
\psfrag{u4}{$u_4$}
\epsfig{file=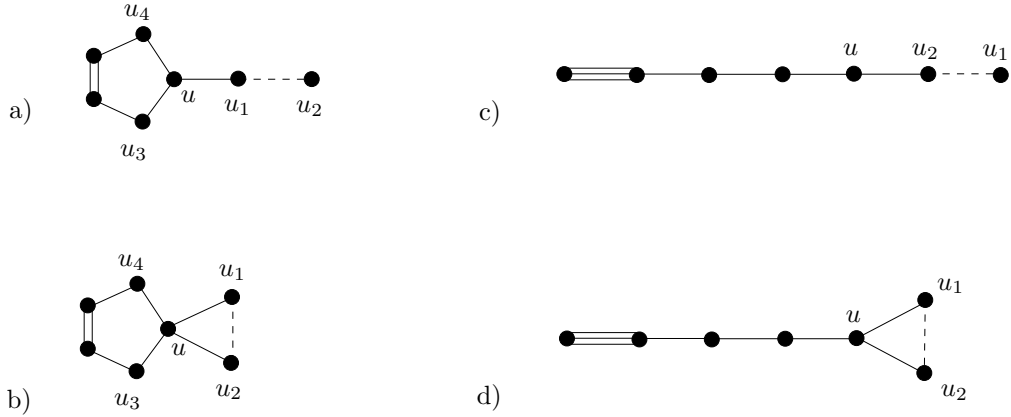,width=0.92\linewidth}
\caption{5-prisms satisfying condition (1) of
Lemma~\ref{nok}.}
\label{prisms}  
\end{center}
\end{figure}

Let $L_1$ be a unique  Lann\'er subdiagram of
$\o S_0$ of order 5, and denote by $S_1$ a unique subdiagram of
$L_1$  
of the type $H_4$ or $F_4$.
Denote $u=L_1\setminus S_1$.
The diagram consists of the following parts:
\begin{center}
\psfrag{L1}{$L_1$}
\psfrag{S1}{$S_1$ {\small = $H_4$ or $F_4$}}
\psfrag{S0}{$S_0$ {\small = $H_4$ or $F_4$}}
\psfrag{S'}{$\overline S_0=\Sigma_{S_0}$ 
{\small = one of the 5-prisms from Fig.~\ref{prisms}}}
\psfrag{u1}{$u_1$}
\psfrag{u2}{$u_2$}
\psfrag{u}{$u$}
\psfrag{v}{$a_1$}
\psfrag{w}{$a_2$}
\psfrag{a}{$a$}
\psfrag{b}{$b$}
\epsfig{file=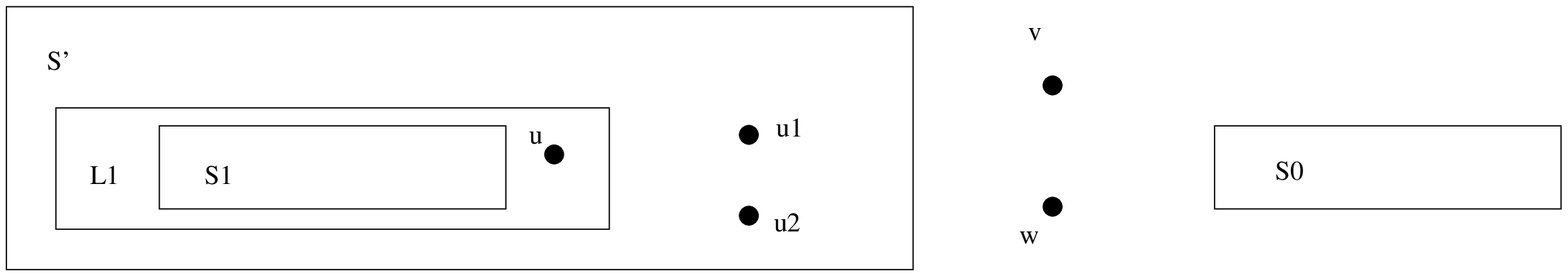,width=0.9\linewidth}
\end{center}
where $a_1$ and $a_2$ attach to $S_0$, and $S_0$ is not connected to
$\overline S_0$. 
 Since any two indefinite subdiagrams of
$\Sigma$ are joined in $\Sigma$ by at 
least one edge, each of $a_1$ and $a_2$ attaches to $L_1$. If they both
attach to $S_1$, then $S_1$ has at least 3 neighbors and we refer to
Case~1. So, one of $a_1$ and $a_2$, say $a_1$, attaches to $u=L_1\setminus
S_1$  (see Fig.~\ref{prisms}) and does not attach to $S_1$.
%
%
%

Consider the diagrams (a) and (b) with nodes indexed as shown in   
Fig.~\ref{prisms}.  Let $S_2=\o {S_0}\setminus\[u_2,u_4\]$. 
$S_2$ is a diagram of type $B_5$ with at
least 3 bad neighbors ($u_4$, $u_2$ and $a_1$).  Thus, $a_2$ does not
attach to $S_2$. Considering a subdiagram
$S_3=\o S_0\setminus\[u_2,u_3\]$ instead of $S_2$, we obtain
that  $a_2$
does not attach to $u_4$ either. Thus, $a_2$ does not attach to $L_1$, and
an indefinite subdiagram $\[S_0,a_2\]$ is not joined with $L_1$. This
means that $\Sigma$ is superhyperbolic, which is impossible.

Now let $\o S_0$ be one of the diagrams (c) and (d) with nodes
indexed as shown in Fig.~\ref{prisms}. 
We have two possibilities: $a_2$ either attaches to $S_1$ or not.
Suppose that $a_2$ does not attach to $S_1$. Then $P(S_1)$ is a
Coxeter 5-polytope with 8 facets. However, the Coxeter diagrams of
such polytopes do not satisfy Lemma~\ref{nok} (see~\cite{n3}).
Thus, $a_2$ is joined with $S_1$ and the diagram $S_1$ has exactly two
neighbors (otherwise we refer to Case~1), so the polytope $P(S_1)$ is
a Coxeter 5-prism. 
By the reasoning of the previous paragraph 
we may assume that the Coxeter diagram $\Sigma_{S_1}$ of $P(S_1)$ is of the
type (c) or (d) (see Fig.~\ref{prisms}). Notice that the diagram
$\o S_1=\Sigma_{S_1}$ contains $S_0$. Therefore, $S_0$ is of the type $H_4$.

Further, the diagram $\o S_1$ contains $u_1$, $u_2$ and $a_1$.
Since $P(S_1)$ is a 5-prism and $\o {S_1}=\[S_0,u_1,u_2,a_1\]$, 
$\o S_1$
contains a Lann\'er diagram of order 5. The nodes $u_1$ and $u_2$ do
not attach to $S_0$, so $\[S_0,a_1\]$ is a Lann\'er
diagram. It follows that it should be joined with $L_1$, and by
Lemma~\ref{2lan} $a_1$ attaches to $u$  by a dotted edge. Therefore, the
diagram $\[S_0,a_1,u, S_1\]$ is
\smallskip
\begin{center}
\psfrag{u}{$u$}
\psfrag{v}{$a_1$}
\psfrag{S0}{$S_0$}
\psfrag{S1}{$S_1$}
\epsfig{file=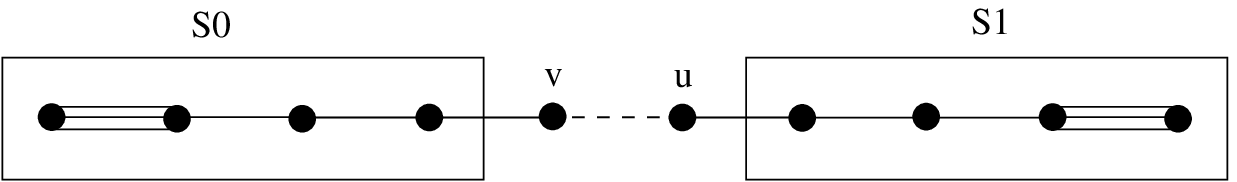,width=0.4\linewidth}
\end{center}
Recall that the nodes $u_1$ and $u_2$ are joined by a dotted
edge. Since $u_1,u_2\in \o S_0$ and  $u_1,u_2\in \o S_1$,
they may be 
joined in $\Sigma\setminus a_2$ with $u$ and $a_1$ only. We 
obtain the following four possibilities for the diagram 
$\Sigma\setminus a_2$:  
\begin{center}
\medskip
\psfrag{p1}{$v_1$}
\psfrag{p2}{$v_2$}
\psfrag{p3}{$v_3$}
\psfrag{p4}{$v_4$}
\psfrag{p5}{$a_1$}
\psfrag{p6}{$u$}
\psfrag{p7}{$v_5$}
\psfrag{p8}{$v_6$}
\psfrag{p9}{$v_7$}
\psfrag{p10}{$v_8$}
\psfrag{p11}{$u_1$}
\psfrag{p12}{$u_2$}
\epsfig{file=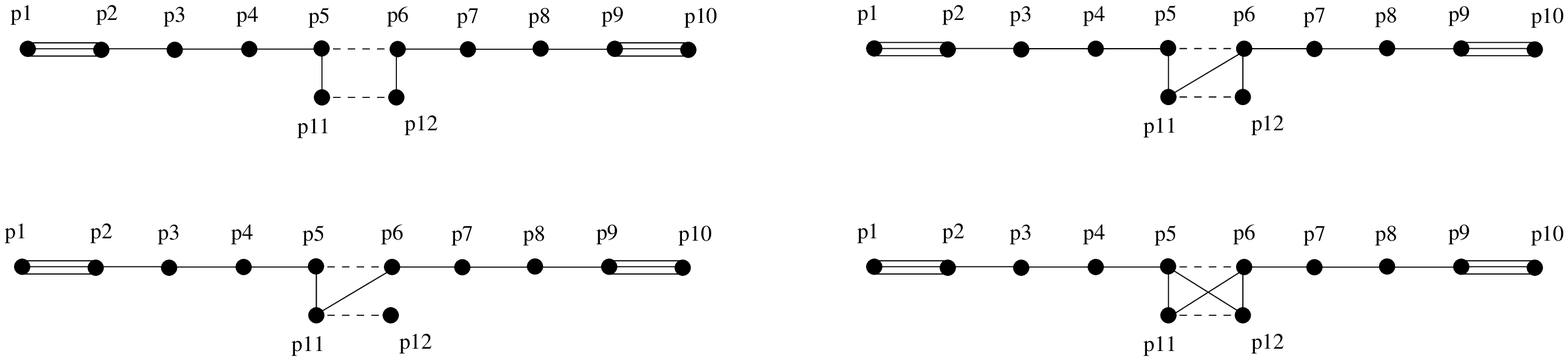,width=0.92\linewidth}
\end{center}

By the assumption, the node $a_2$ attaches to  $S_1=\[v_5,v_6,v_7,v_8\]$.
Since $a_2\notin \Sigma_{S_0}$, $a_2$ also attaches to   
$S_0=\[v_1,v_2,v_3,v_4\]$. Further, $a_2$ attaches to
$\[u_{1},u_{2}\]$, otherwise an indefinite diagram
$\[S_0,a_2\]$ is not joined with the Lann\'er diagram $\[u_{1},u_{2}\]$.
 Since the diagrams $\[v_2,v_3,v_4,a_1,u_1\]$,
$\[v_7,v_6,v_5,u,u_{1}\]$ and
$\[v_7,v_6,v_5,u,u_{2}\]$
already have 3 bad neighbors each, no dotted edge ends in
$a_2$ (see (1) of Lemma~\ref{nok}). 
Carefully examining possible 
multiplicities of edges, it is easy
to see that we always obtain either a Lann\'er subdiagram of order 3, or
a parabolic subdiagram, or a subdiagram of the type $H_4$ with at
least 3 neighbors, and we refer to Case 1.   

\bigskip

\noindent
{\bf Case 3:} $P(S_0)$ has $d-1$ facets. Then  $P(S_0)$ is a Coxeter
$(d-4)$-polytope with $(d-4)+3$ 
facets. Since $d\ge9$, we have $d-4\ge 5$. According to~\cite{n3},
Coxeter diagrams of 
all such polytopes either do not satisfy the conditions of
Lemma~\ref{nok}, or contain a subdiagram of the type $H_4$ with at
least two neighbors, so we refer to the previous cases.     

\end{proof}

\begin{lemma}
\label{no4-9}
Suppose that one of the following holds:\\
1) $d\ge 9$, or\\ 
2) $d=8$, $\Sigma$ satisfies conditions (1)--(3) of Lemma~\ref{nok}
and $\Sigma$ 
does not contain any Lann\'er diagram of order $5$.\\
Then $\Sigma$ does not contain any Lann\'er diagram
of order $4$. 

\end{lemma}

\begin{proof}
Suppose the contrary. 
Let $L_0$ be a Lann\'er subdiagram  of $\Sigma$ of order 4. 
Denote by $S_0$ a subdiagram of $L_0$ of the type $H_3$ (if any) or
$B_3$ (otherwise). Then $P(S_0)$ is a Coxeter 
$(d-3)$-polytope with at most $(d-3)+3$ facets, and its Coxeter
diagram $\Sigma_{S_0}$ does not contain
Lann\'er diagrams of order 5 (Lemmas~\ref{no5-9} and~\ref{lift}). 
If $P(S_0)$ is a simplex or a product of two
simplices, then $d-3\le 3$ or  $d-3\le 4$ respectively, so in this
case $d\le 7$. 

Suppose that $P(S_0)$ has $(d-3)+3$ facets. According to~\cite{n3}, for
$d\ge 8$ the Coxeter diagrams of almost all Coxeter 
$(d-3)$-polytopes with $d$
facets contain Lann\'er diagrams of order 5. There is only one
exclusion: a 5-polytope with 8 facets having the following Coxeter
diagram $\Sigma_{S_0}$: 
\medskip
\begin{center}
\psfrag{e1}{\small $e_1$}
\psfrag{e2}{\small $e_2$}
\psfrag{b1}{\small $b_1$}
\psfrag{b2}{\small $b_2$}
\psfrag{c1}{\small $c_1$}
\psfrag{c2}{\small $c_2$}
\psfrag{d1}{\small $d_1$}
\psfrag{d2}{\small $d_2$}
\epsfig{file=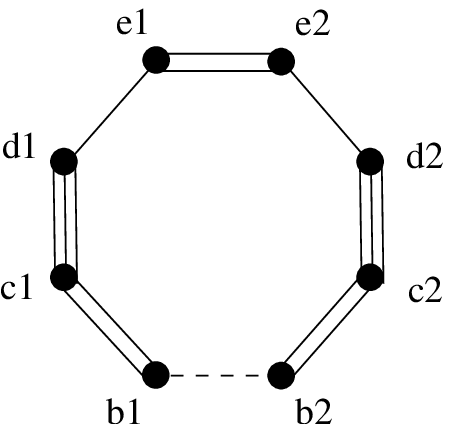,width=0.2\linewidth}
\end{center}
%

By Cor.~\ref{dif2}, no node of $\o S_0\setminus \[b_1,b_2\]$
is a good neighbor of $S_0$.
The node $b_1$ is a good neighbor of $S_0$,
otherwise $\Sigma$ contains an order 3 Lann\'er diagram  
$\[b_1,c_1,d_1\]$, which contradicts  the assumption. 
Thus, the edge  $b_1c_1$ is a simple edge of $\Sigma$, 
and the subdiagram 
$\[d_1,c_1,b_1,x_1,x_2\]\subset \Sigma$ is a Lann\'er diagram of order 5
(where $x_1$ and $x_2$ are the ends of the simple edge of $S_0$).
The contradiction with the assumption of the lemma completes the proof.

\end{proof}

\begin{theorem}
There is no compact hyperbolic Coxeter $d$-polytope with $d+4$ facets
for $d\ge 9$.

\end{theorem}

\begin{proof}
Lemmas~\ref{nok},~\ref{no5-9} and~\ref{no4-9} imply that all Lann\'er
diagrams in $\Sigma$ have order $2$. By~\cite[Satz~6.9]{Ess2}, in this case
$d\le 4$.  

\end{proof}

\section{Absence of polytopes in dimension 8}
\label{d8,n+4}

In this section, we prove that no compact Coxeter polytope with 12
facets exists in $\H^8$. 
Suppose there exists an 8-polytope $P$ with 12 facets. 
First, we show some properties of its Coxeter diagram $\Sigma$. Most
of them repeat ones that  
hold for polytopes in large dimensions. However, the proofs in
eight-dimensional case are much more complicated.

\begin{lemma}[~\cite{nodots}, Lemma~1]
\label{same}
Let $\Sigma(P)$ be a Coxeter diagram of a Coxeter $d$-polytope
$P$. Then no proper subdiagram of $\Sigma(P)$ is a Coxeter 
diagram of a finite volume Coxeter $d$-polytope.

\end{lemma}

In particular,
$\Sigma$ does not contain $\Sigma(P_8)$ as a proper subdiagram.

\begin{lemma}
\label{nok8-2}
Suppose that nodes $v_1$ and $v_2$ of $\Sigma$ are joined by an
edge of multiplicity $k\ge 2$. Then\\
(1) the subdiagram $\[v_1,v_2\]$ has at most one bad neighbor;\\
(2) $k\le 3$;\\
(3) if the subdiagram $\[v_1,v_2\]$ has a bad neighbor, then it has
at least one good neighbor, too. 

\end{lemma} 

\begin{proof}
Denote by $S_0$ the diagram $\[v_1,v_2\]$. If $S_0$ has at least two
bad neighbors, $P(S_0)$ is a Coxeter 6-polytope with at most 8
facets, which is impossible~\cite{K},~\cite{Ess2}.

Suppose $k\ge 4$ or $S_0$ has a unique neighbor which is bad. 
In particular, all neighbors of $S_0$ are bad. In this case $P(S_0)$
is a Coxeter 6-polytope with 
at most 9 facets, and its Coxeter diagram $\Sigma_{S_0}$ is a subdiagram
of $\Sigma$. 
By~\cite{n3}, there are only three such polytopes, their diagrams are 
shown in Fig.~\ref{d6n9}.
A Coxeter diagram of each of these three polytopes
contains a multiple edge with at least two bad neighbors
(see Fig.~\ref{d6n9}).

\begin{figure}[htb]
\begin{center}
\psfrag{1}{$P_6^1$}
\psfrag{2}{$P_6^2$}
\psfrag{3}{$P_6^3$}
\psfrag{10}{\scriptsize $10$}
\epsfig{file=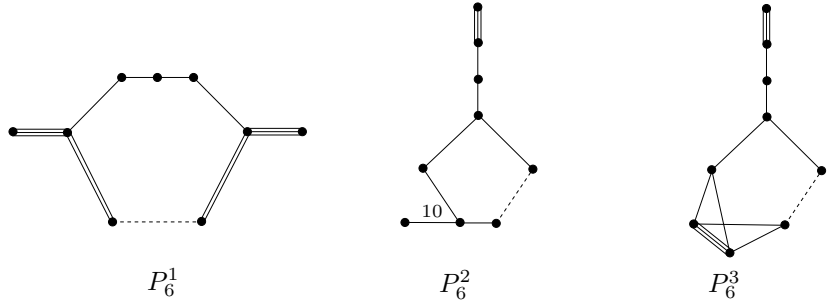,width=0.75\linewidth}
\caption{Compact hyperbolic Coxeter 6-polytopes with 9 facets. }
\label{d6n9}
\end{center}
\end{figure}

\end{proof}

\begin{lemma}
\label{no3-8}
(1) $\Sigma$ contains no Lann\'er diagram of order $3$.\\
(2) No end of a dotted edge is incident to a triple edge.\\
(3) If $\Sigma$ contains a unique dotted edge, then no end of that
edge is incident to a double edge.
\end{lemma}

\begin{proof}
Suppose that the lemma is broken. Denote by $L_0$ a Lann\'er diagram
of order 3 or a subdiagram consisting of adjoint dotted and multiple
edges. $L_0$ contains an edge of multiplicity 2 or 3. We consider the
following two cases:  

\bigskip

\noindent
{\bf Case 1:}  $L_0$ contains an edge of multiplicity 3.\\
Denote by $S_0$ the subdiagram consisting of this edge. $P(S_0)$ is a
Coxeter
6-polytope with at most 9 facets, so we may assume that $P(S_0)$ has
exactly 9 facets. 
Then $\Sigma_{S_0}$ is one of the diagrams shown in
Fig.~\ref{d6n9}. By Cor.~\ref{dif}, the diagram $\overline
S_0\subset\Sigma$ can be obtained
from $\Sigma_{S_0}$ by replacing some dotted edges by ordinary edges,
and some edges labeled by 10 by simple edges. Applying (or not) such a
procedure to the Coxeter diagrams of the polytopes $P_6^1$ and $P_6^3$
we see that the resulting diagrams contain either a parabolic
subdiagram or a multiple edge with at least two bad neighbors.

Now suppose that $P(S_0)=P_6^2$. Denote by $u$ a unique bad neighbor
of $S_0$ (it is unique by Lemma~\ref{nok8-2}). 
By Cor.~\ref{dif}, $\Sigma_{S_0}$   can by transformed to $\o S_0$ 
by the change of some edges labeled by 10 into simple edges and some 
dotted edges into ordinary edges.
By Cor.~\ref{dif2}, the dotted edge of $\Sigma_{S_0}$ remains dotted in
$\o S_0$.
By Lemma~\ref{nok8-2}, the edge of $\Sigma_{S_0}$ labeled by 10
is a simple edge of $\o S_0$. Therefore, the leaf
$a\in \Sigma_{S_0}$ that is incident to the edge labeled by 10
is a good neighbor of 
$S_0$ in $\Sigma$. The remaining nodes of $\Sigma_{S_0}$
cannot be good neighbors of $S_0$ in view of Cor.~\ref{dif2}.
So, the subdiagram $\Sigma\setminus u$ is a Coxeter diagram 
$\Sigma(P_8)$, which is impossible by Lemma~\ref{same}.

\bigskip

\noindent
{\bf Case 2:}  $L_0$ contains no edge of multiplicity 3.\\ 
Denote by $S_0$ the subdiagram consisting of a double edge contained
in $L_0$. $P(S_0)$ is again a 6-polytope 
with 9 facets. By Prop.~\ref{al}, the diagram $\overline
S_0\subset\Sigma$  can be obtained from  
$\Sigma_{S_0}$ by replacing some dotted edges by ordinary (or empty) edges,
and some double edges by simple or empty edges. Applying (or not) such a
procedure to the Coxeter diagrams of the polytopes $P_6^3$ and $P_6^2$
we see that the resulting diagrams contain a multiple edge (a triple one or
labeled by 10 respectively)  with at least two bad neighbors. 

Suppose that $P(S_0)=P_6^1$. Both double edges from $\Sigma_{S_0}$ become
simple or empty edges in $\overline S_0$, otherwise 
$\Sigma$ contains a Lann\'er diagram with a triple edge and
we refer to Case 1. 
Furthermore, by Cor.~\ref{dif2}, only one end of each double edge can
be a good neighbor of $S_0$. Hence, by Prop.~\ref{al} 
both double edges of $\Sigma_{S_0}$ 
are simple edges in $\o S_0$.
This implies that a
unique dotted edge in $\Sigma_{S_0}$ remains a dotted edge 
in $\overline S_0$ 
(otherwise the diagram $\o S_0$ is superhyperbolic). 
So, in this case we may assume that $L_0$ contains no dotted edge. By
Prop.~\ref{al}, there exist only two good neighbors of $S_0$ --- the
ends of the dotted edge. They attach to $S_0$ by simple edges. 
Consider the diagram $S_1$ of the type $H_4$ containing a triple edge, a
simple edge that was double in $\Sigma_{S_0}$ and a simple edge joining
$S_0$ with its good neighbor.
Since
$L_0$ is a Lann\'er diagram without edges of multiplicity greater than
two, $L_0$ contains 3 edges, and we obtain that the  diagram 
$S_1$ of the type $H_4$ has at least 4 bad neighbors (two neighbors in
$\o S_0$ and two ones in $S_0$). This contradicts lemma~\ref{bad}, which
completes the proof of the lemma. 
  
%
%

\end{proof}

The following lemma is a corollary of the main result of~\cite{nodots}. 

\begin{lemma}[\cite{nodots}, cor. of Th.~A]
\label{e2}
If $\Sigma(P)$ is a Coxeter diagram of a compact Coxeter
polytope $P\subset \H^d$, where $d>4$, then
$\Sigma(P)$ contains a dotted edge.

\end{lemma}

\begin{lemma}
\label{order3}
Let $S\subset \Sigma$ be an elliptic subdiagram of order $3$
containing no component of the type $A_n$. Then $S$ has at most $2$
bad neighbors. 

\end{lemma}

\begin{proof}
Suppose that $S$ has $3$ or more bad neighbors.
Then $\Sigma_S$ is a Coxeter $5$-polytope with at most $6$ facets,
which is impossible.

\end{proof}

\begin{lemma}
\label{5-e2-3}
Suppose that $\Sigma$ contains a Lann\'er subdiagram $L_0$ of order
$5$ and that
the dotted edge is unique in $\Sigma$. Then $\Sigma$ contains a subdiagram
$S$ of the type $F_4$ or $H_4$ such that $P(S)$ is not a simplex.   
 
\end{lemma}

\begin{proof} 
Suppose that the lemma is broken. Then we may assume that for any diagram
$S'\subset \Sigma$ of the type $H_4$ or $F_4$ the face $P(S')$ is a Coxeter
simplex.

Denote by $S_0$ a subdiagram of $L_0$ of the type $H_4$ or $F_4$. In
particular, $P(S_0)$ is a Coxeter 4-simplex.   
The Coxeter diagram $\Sigma_{S_0}$ of $P(S_0)$ is a Lann\'er diagram
of order 5. Denote it by $L_1$. Denote by $S_1$ a subdiagram of $L_1$
of the type $H_4$ or $F_4$. $P(S_1)$ is a Coxeter 
simplex, and its Coxeter diagram
$\overline S_1= \Sigma_{S_1}$ is a Lann\'er diagram of order 5. 
Denote it by $L_2$. Notice
that $S_0\subset L_2$, so let $S_2=S_0$. Let $v_i=L_i\setminus S_i$,
$i=0,1,2$. 

Since $S_2=S_0$ does not attach to $L_1$, $v_2$ is joined with $L_1$. On
the other hand, $L_2$ does not attach to $S_1$. Thus, the subdiagram
$\[L_1, L_2\]$ consists of two Lann\'er diagrams $L_1$ and $L_2$ joined
by a unique edge $v_1v_2$, such that
$L_i\setminus v_i$ are of the type $H_4$ or $F_4$. By
Lemma~\ref{2lan}, $v_1v_2$ is a dotted edge.  

Now suppose the dotted edge $v_1v_2$ is unique.
Denote by $a$ and $b$ two nodes not contained in $L_1$ and $L_2$. Note
that if $v_2\ne v_0$ then either $a$ or $b$ coincides with $v_0$. 
Both $a$ and $b$ are neighbors of $S_1$ and $S_2$. The diagram
$\Sigma$ consists of the following parts:

\begin{center}
\psfrag{L1}{$L_1$}
\psfrag{L2}{$L_2$}
\psfrag{S1}{$S_1$}
\psfrag{S2}{$S_2=S_0$}
\psfrag{v1}{$v_1$}
\psfrag{v2}{$v_2$}
\psfrag{a}{$a$}
\psfrag{b}{$b$}
\epsfig{file=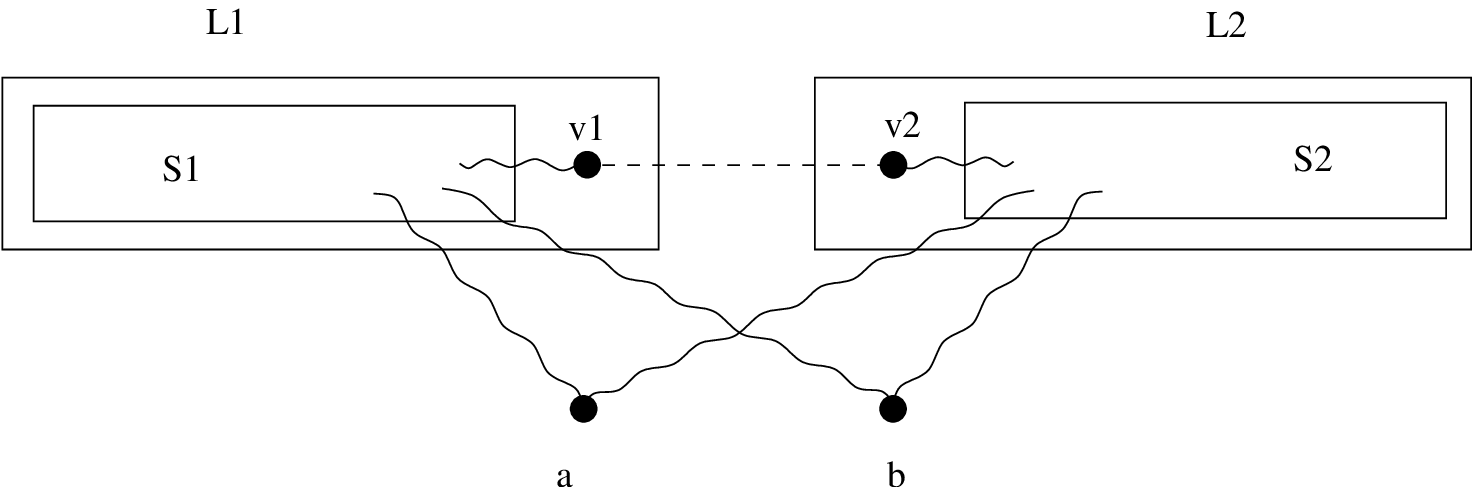,width=0.55\linewidth}
\end{center}
where a wave line means that the vertex is joined with the subdiagram
by some non-dotted edges, and some additional non-dotted edges may
appear in the subdiagram $\[a,b,v_1,v_2\]$.

By Lemma~\ref{no3-8}, vertices $v_1$ and $v_2$ attach to $S_1$ and
$S_2$ respectively by simple edges only. Thus, $L_1$ (as well as
$L_2$) is one of the diagrams $\L_1^5$, $\L_4^5$ or $\L_5^5$ 
(see Table~\ref{Lan}).

Suppose that $L_1=\L^5_5$. Let $w_1\in S_1$ 
be a vertex attached to $a$. If $w_1a$ is a triple edge, we obtain
either a Lann\'er diagram of order 3 or a diagram of the type $H_3$
with at least three bad neighbors, which is impossible by
Cor.~\ref{simplex}. If $w_1a$ is a double edge or a
simple edge, we obtain either a parabolic diagram or a Lann\'er
diagram of order 3. So, the multiplicity of the edge $w_1a$ should be
greater than 3, which contradicts  Lemma~\ref{nok8-2}. 

Thus, both $S_1$ and $S_2$ are of the type $H_4$. 
For each of 4 possible pairs of diagrams $L_1$ and $L_2$ there exists
a unique label of the dotted edge $v_1v_2$ such that the determinant
of the diagram $\Sigma\setminus\[a,b\]$ vanishes. The label equals
$(1+\sqrt{5})/2$ for $L_1=L_2=\L^5_1$, \ $4+2\sqrt{5}$ for
$L_1=L_2=\L^5_4$, and $(3+\sqrt{5})/\sqrt{2}$ for $L_1=\L^5_1$,
$L_2=\L^5_4$ (or $L_2=\L^5_1$, $L_1=\L^5_4$).

Denote by $u_1\in
S_1$ and $u_2\in S_2$ the leaves of $\Sigma\setminus\[a,b\]$ incident
to the triple edges. If we assume that each of $a$ and $b$ attaches to
$S_1\setminus u_1$ and to $S_2\setminus u_2$, all edges joining
them with $S_i\setminus u_i$ are  simple,
and $a$ is not joined with $b$,
 then the diagram $\Sigma\setminus\[v_1,v_2\]$ contains 
 a parabolic  subdiagram of the type $\widetilde A_m$ for some $m$,
$2\le m \le 7$.
If in addition there are some multiple edges joining $a$ and $b$ with
$S_1\setminus u_1$ and to $S_2\setminus u_2$, or $a$ is joined with
$b$, then $\Sigma\setminus v_1$ or  $\Sigma\setminus v_2$
contains a subdiagram of the type
$H_3$ with at least 3 bad neighbors, or a Lann\'er diagram of order 3,
or a parabolic diagram of the type $\widetilde B_3$. 
Therefore, at least one of $a$ and $b$, say
$a$, can attach to one of $S_i\setminus u_i$, say $S_1\setminus u_1$,
by multiple edges only. 
This means that we have got two cases: $a$
either attaches to $S_1\setminus u_1$ by a unique multiple edge
(Lemma~\ref{no3-8})
or does not attach to $S_1\setminus u_1$ (in the latter case $a$
attaches to $u_1$). 

\bigskip

\noindent
{\bf Case 1:} $a$ attaches to $S_1\setminus u_1$ by a unique multiple
edge $au$. It is easy to see that in this case $L_1=\L^5_4$, $u$ is
the leaf of $L_1$ different from $v_1$, and $a$ does not attach 
to $v_1$ 
(otherwise we obtain either a diagram of the type $H_3$ with at least
3 bad neighbors, or
a Lann\'er diagram of order 3, or a
parabolic subdiagram). If $au$ is a double edge then 
$\[a,L_1\setminus u_1\]$ is a parabolic diagram
$\widetilde B_4$. Suppose that $au$ is a triple edge.  The vertex $a$
is not joined with $u_1$, otherwise the diagram $\[u,a,u_1\]$ of the
type $H_3$ has at least 3 bad neighbors, which contradicts 
Cor.~\ref{simplex}. Now consider the diagram $\Sigma_1=\[L_1,a,v_2\]$. 
The vertex $a$ may attach (by simple edge only, see Lemma~\ref{no3-8}) 
to $v_2$ only. An elementary computation shows that $\Sigma_1$ is
superhyperbolic whenever $a$ attaches to $v_2$ or not. 

\bigskip

\noindent
{\bf Case 2:} $a$ does not attach to $L_1\setminus u_1$. Therefore,  
$a$ attaches to $u_1$ by a simple edge. Again, consider the diagram
$\Sigma_1=\[L_1,a,v_2\]$. The vertex $a$ may attach (by
simple edge only) to  $v_1$ and $v_2$ only. For each of 4 pairs of diagrams
$L_1$ and $L_2$ the diagram $\Sigma_1$ is superhyperbolic whenever 
$a$ attaches to $v_1$ and (or) $v_2$ or not. 

\end{proof}

\begin{lemma}
\label{nok8-1}
Any node of $\Sigma$ is incident to at most one dotted edge.

\end{lemma}

\begin{proof}
Suppose the contrary. Denote by $u$ a node of $\Sigma$ incident to at least
two dotted edges. Since the number of bad neighbors of $u$ cannot
exceed 3, $u$ is incident to either two or three dotted edges. If
there are three, the polytope $P(u)$ is a simplicial
7-face of $P$, which contradicts  Cor.~\ref{simplex} of Lemma~\ref{lift}.

Thus, $u$ is incident to exactly two dotted edges, and $P(u)$ is a 
7-polytope with 9 facets, i.e. a polytope whose  diagram of the missing
faces is $\lf k,9-k \rf_1$ for some $k<9$.
Since the missing faces
diagram of such a polytope contains no large missing faces,
it is of the type $\lf 4,5 \rf_1$.
By Lemma~\ref{lift54}, this
diagram of missing faces  
has a unique lifting, namely the diagram $\Theta_1$ 
shown in Fig.~\ref{10only5}. 
Denote the nodes of $\Theta_1$ as shown in Fig.~\ref{lift-v},
denote
by $v$ and $w$ the remaining two nodes of $\Sigma$. Each of them 
attaches to the node $u$ by a dotted edge and attaches to the
subdiagram $\[v_5,\dots,v_9\]\subset \Sigma$, where the nodes are indexed
as shown in Fig.~\ref{lift-v}.
\begin{figure}[htb]
\begin{center}
\smallskip
\psfrag{v}{$u$}
\psfrag{v1}{$v_1$}
\psfrag{v2}{$v_2$}
\psfrag{v3}{$v_3$}
\psfrag{v4}{$v_4$}
\psfrag{v5}{$v_5$}
\psfrag{v6}{$v_6$}
\psfrag{v7}{$v_7$}
\psfrag{v8}{$v_8$}
\psfrag{v9}{$v_9$}
\epsfig{file=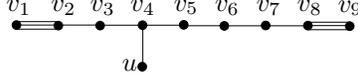,width=0.32\linewidth}
\caption{A unique lifting of $[4,5]_1$. }
\label{lift-v}
\end{center} 
\end{figure}

\noindent
{\bf Case 1:} no dotted edge joins $\[ v,w\]$ with $\Sigma\setminus
\[v,w,u\]=\[v_1,\dots,v_9\]$.   

Let $S_0=\[v_2,\dots,v_8\]$. The face $P(S_0)$ is one-dimensional, so
$S_0$ is contained in exactly two elliptic diagrams of order 8. This
means that one of $v$ and $w$, say $v$, either is not joined with $S_0$
or is a good neighbor of $S_0$. In particular, any double edge from
$v$ to $S_0$ may end in $v_2$ or $v_8$ only. Now Lemma~\ref{no3-8}
implies that no double edge joins $v$  with $\Sigma\setminus
\[v,w,u\]$. Thus, $v$ may be joined by 
simple edges only with $v_1$, $v_9$ and one node of $S_0$. A
straightforward calculation shows that none of obtained subdiagrams
$\Sigma\setminus \[w,u\]$ is admissible with positive inertia index at
most 8.

%
%

\bigskip

\noindent
{\bf Case 2:} at least one of $v$ and $w$, say $v$, is joined by a dotted edge
with some node $v_x\in \Sigma\setminus \[v,w,u\]$.

In this case $P(v)$ is a 7-face of $P$ with 9 facets, 
i.e. a face with a diagram of missing faces of the type $\lf 5,4 \rf_1$.
A unique lifting of its missing faces is the diagram $\Theta_1$ 
(Lemma~\ref{lift54}). In other words, 
diagram $\Sigma \setminus\[u,v_x\]$ looks as shown in 
Fig.~\ref{lift-v}, where
$v$ takes place of $u$, $w$ takes place of $v_x$,
and, possibly, the node $v$ is not joined with $v_4$, but joined with 
either $v_6$ (if $v_x\ne v_6$) or with $w$ (if $v_x=v_6$).

If $v_x$ and $w$ are not joined by a dotted edge, we obtain a
parabolic diagram $\[v_{x-1},v_x,v_{x+1},w\]$ (if $v_x\ne v_1$, $v_9$)
or a Lann\'er diagram  of order 3 (otherwise). Thus, $v_x$ and
$w$ are  joined by a dotted edge. Now consider $P(v_x)$. It is
again a 7-polytope with 9 facets, but in the lifting of its missing faces 
diagram the node $v_x$ is joined with at least two other nodes, which
contradicts Lemma~\ref{lift}.   

\end{proof}

\begin{lemma}
\label{2dot}
$\Sigma$ contains at least two dotted edges.

\end{lemma}

\begin{proof}
By Lemma~\ref{e2}, there exists at least one dotted
edge. Suppose it is unique. 
Then, by Lemmas~\ref{nok8-2} and \ref{no3-8},
$\Sigma$ satisfies  conditions (1)--(3) of Lemma~\ref{nok}. 
Below we show that this implies that no
Lann\'er diagram of order 5 is contained in $\Sigma$. In view of
Lemmas~\ref{no4-9},~\ref{no3-8} and~\ref{disjoint}, this proves the lemma.

We assume that there exists a Lann\'er diagram $L_0\subset \Sigma$ of order
5. Denote by $S_0$ a subdiagram of $L_0$ of the type $H_4$ or $F_4$.
The Coxeter diagrams of all Coxeter 4-polytopes with 7 facets
contain either a Lann\'er diagram of order 3 or at least two dotted
edges. The Coxeter diagrams of all Coxeter 4-polytopes with 
6 facets that are not prisms contain Lann\'er diagrams of order 3. Thus,
$P(S_0)$ is either a Coxeter 4-prism or a Coxeter
4-simplex. 

By Lemma~\ref{5-e2-3}, we can assume that $P(S_0)$ is a prism.
Consider two cases.

\bigskip

\noindent
{\bf Case 1:}
$P(S_0)$ is a prism whose Coxeter diagram $\Sigma_{S_0}$ is different from 
\smallskip
\begin{center}
\psfrag{2,3}{\scriptsize $2,3$}
\psfrag{v1}{}
\psfrag{v2}{}
\psfrag{v3}{}
\psfrag{v4}{}
\psfrag{v5}{}
\psfrag{v6}{}
\epsfig{file=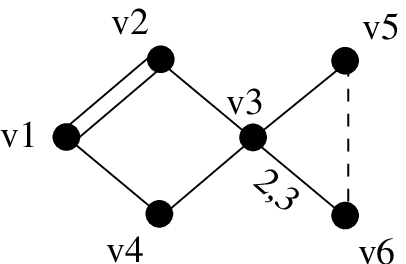,width=0.26\linewidth}
\end{center} 
where the label ``$2,3$'' means that the nodes are either not joined
or joined by a simple edge.

Looking at the list of Coxeter diagrams of 4-prisms~\cite{K}, we
see that $\Sigma_{S_0}=\overline S_0$ has a subdiagram 
$S_1$ of the type $H_4$ or $F_4$
with a dotted edge attached to $S_1$. 
Since $\Sigma$ contains a unique dotted edge,
the Coxeter diagram $\Sigma_{S_1}=\overline S_1$ of $P(S_1)$ 
does not contain dotted
edges, so $P(S_1)$ 
is a simplex. Notice that $S_0\subset \overline S_1$.
Denote by $u$ and $v$ the neighbors of $S_0$, denote by $a$
the end of the dotted edge not contained in $S_1$, and by $b$ the remaining
node of $\Sigma$ not contained in $\[S_0,S_1\]$.
So, the diagram $\Sigma$ consists of the following parts:

\begin{center}
\psfrag{S0}{$S_0=H_4$ or $F_4$}
\psfrag{S1}{$S_1=H_4$ or $F_4$}
\psfrag{s0}{}  
\psfrag{s1}{$\overline S_0=\Sigma_{S_0}=$ {\small a 4-prism}}
\psfrag{u}{$u$}
\psfrag{v}{$v$}
\psfrag{a}{$a$}
\psfrag{b}{$b$}
\epsfig{file=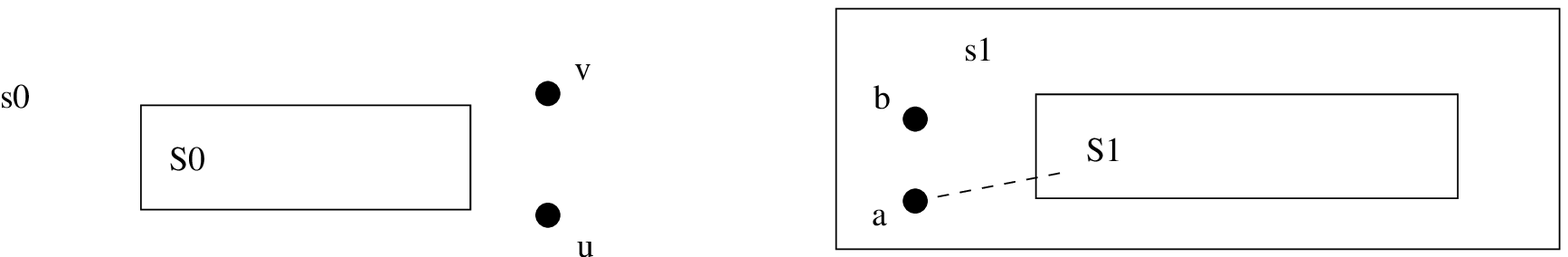,width=0.85\linewidth}
\end{center}

%
\noindent
where $u$ and $v$ attach to $S_0$, and $S_0$ does not attach to
$\overline S_0$. 
Since  $P(S_1)$
is a simplex, $b$ and one of $u$ and $v$, say $v$, attach to $S_1$,
and $\[S_0,u\]=\overline S_1$ is a Lann\'er diagram.

Consider the diagram $\Sigma\setminus \[v,a\]$. It consists of a Lann\'er
diagram  $\[S_0,u\]=\overline S_1$ and an indefinite diagram 
$\[S_1,b\]$. 
These two subdiagrams may be joined by a unique non-dotted
edge $ub$ only. It is easy to see that the diagram 
$\Sigma \setminus\[v,a\]$ is superhyperbolic.

\bigskip

\noindent
{\bf Case 2.}
$P(S_0)$ is a prism whose Coxeter diagram $\Sigma_{S_0}$ is  
\begin{center}
\smallskip
\psfrag{2,3}{\scriptsize $2, 3$}
\psfrag{v1}{$v_1$}
\psfrag{v2}{$v_2$}
\psfrag{v3}{$v_3$}
\psfrag{v4}{$v_4$}
\psfrag{v5}{$v_5$}
\psfrag{v6}{$v_6$}
\epsfig{file=./pic_new/prism4-i.eps,width=0.26\linewidth}
\end{center} 
Denote by $S_1$ a subdiagram of $\Sigma_{S_0}$ of the type $B_4$, say
$S_1=\[v_1,v_2,v_3,v_5\]$. 
So, $\Sigma$ consists of the following parts:
\begin{center}
\psfrag{S0}{$S_0=H_4$ or $F_4$}
\psfrag{S1}{$S_1$}
\psfrag{s0}{}  
\psfrag{s1}{$\overline S_0$}
\psfrag{u}{$u$}
\psfrag{v}{$v$}
\psfrag{v1}{$v_1$}
\psfrag{v2}{$v_2$}
\psfrag{v3}{$v_3$}
\psfrag{v4}{$v_4$}
\psfrag{v5}{$v_5$}
\psfrag{v6}{$v_6$}
\psfrag{2,3}{\scriptsize $2, 3$}
\psfrag{a}{$a$}
\psfrag{b}{$b$}
\epsfig{file=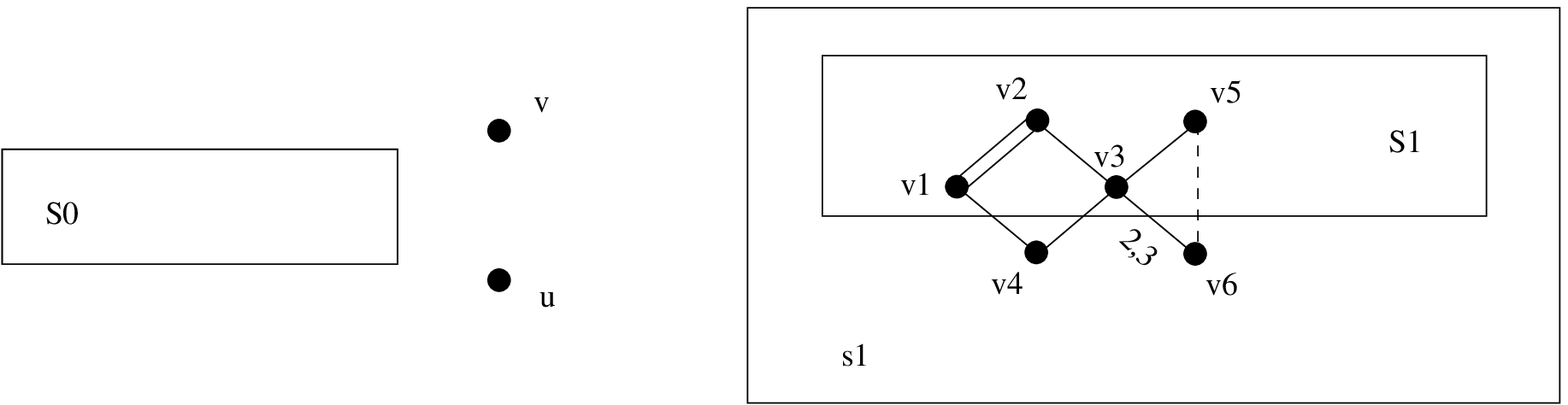,width=0.8\linewidth}
\end{center}
where $u$ and $v$ attach to $S_0$, and $S_0$ does not attach to
$\overline S_0$. Notice that
$S_1$ has at least two bad neighbors,
$v_4$ and $v_6$, so $P(S_1)$ is either an Esselmann polytope,
or a prism, or a simplex.
Consider these three cases.

\smallskip

\noindent
{\bf Case 2.1:} $P(S_1)$ is an Esselmann polytope.\\
Then $\Sigma_{S_1}$ contains two disjoint Lann\'er diagrams of order 3,
while $\o S_1$ contains no subdiagram of order 3. Therefore,
  $\Sigma_{S_1}$ contains two nodes that are good neighbors of 
the diagram $S_1$ in $\Sigma$. This is impossible
by Cor.~\ref{dif2} (we use the list of Esselmann
diagrams~\cite{Ess_eng} and the fact that
$S_1$ is a diagram of the type $B_4$).

\noindent
{\bf Case 2.1:} $P(S_1)$ is a prism.\\
As above, $u$ and $v$ are  neighbors of $S_0$. 
By the assumption, $\Sigma$ contains a unique dotted edge,
hence 
at least one of $u$ and $v$,
say $v$, should be a good neighbor of $S_1$. Then $\[S_1,v\]$ is
of the type $B_5$. If $u$ is also a good neighbor of $S_1$, we
obtain either a parabolic subdiagram of $\Sigma$ or a Lann\'er diagram of
order 3. Thus, $v$ is a unique good neighbor of $S_1$.

Consider the diagram $\overline S_1=\[S_0,u,v\]$. 
The only node of $\overline S_1$ that attaches to $S_1$
is $v$. Thus, 
$\overline S_1$ may differ from $\Sigma_{S_1}$ by multiplicities of edges incident to $v$
only. By Prop.~\ref{al}, any simple edge of $\overline S_1$ incident to $v$ becomes
a double edge in $\Sigma_{S_1}$, and any other edge of $\overline S_1$
incident to $v$ becomes a
dotted edge. Since  $P(S_1)$ is a prism, $\Sigma_{S_1}$ contains a 
unique dotted
edge. At the same time, no Coxeter diagram of a Coxeter 4-prism contains a
node incident to a dotted edge and to a multiple edge
simultaneously. Thus, $v$ is 
joined with a unique node of $\overline S_1 \setminus v=
\[S_0,u\]$, say $w$, and
$\overline S_1$ may be obtained from $\Sigma_{S_1}$ by replacing  
the dotted edge by
a double or a triple edge.   
  
Recall that $v$ is a unique neighbor of $S_1$ contained in $\overline
S_1$, and it is
joined with $v_5$ by a simple edge. Thus, we obtain either a parabolic
subdiagram $\[v_1,v_2,v_3,v_5,v,w\]$ or a subdiagram $\[v_3,v_5,v,w\]$
of the type $H_4$ with at least 4 neighbors.

\smallskip

\noindent
{\bf Case 2.3:} $P(S_1)$ is a simplex.\\
We preserve the notation from Case 1. The difference is that now
one of $u$ and $v$, say $v$, is a bad neighbor of $S_1$. The only way
to attach a bad neighbor to $S_1$ such that no parabolic or Lann\'er
diagram of order 3 appears is to join $v$ by a triple edge with
$v_5$. But in this case we obtain again a subdiagram
$\[v_2,v_3,v_5,v\]$ of the type $H_4$ with at least 4 neighbors.

\bigskip

Hence, no Lann\'er subdiagram of $\Sigma$ of order five exists, and the
lemma is proved.     
   
\end{proof}

\begin{theorem}
\label{no8}
There are no compact hyperbolic Coxeter $8$-polytopes with $12$ facets.

\end{theorem}

\begin{proof}
Suppose there exists a polytope $P$ with a Coxeter diagram $\Sigma$. 
By Lemma~\ref{e2}, $\Sigma$ contains a dotted edge $vw$.
Lemma~\ref{nok8-1} implies that $P(v)$ is a 7-polytope with
10 facets. By Lemma~\ref{2dot}, $P(v)$ has a pair of disjoint
facets. 

There are only four two-dimensional Gale diagrams of order 10
containing a missing face of order 2 and not containing large
missing faces, all of them are listed below. 
\smallskip
\begin{center}
\psfrag{1}{\small $1$}
\psfrag{2}{\small $2$}
\psfrag{3}{\small $3$}
\psfrag{4}{\small $4$}
\psfrag{5}{\small $5$}
\psfrag{532}{}
\psfrag{442}{}
\psfrag{14132}{}
\psfrag{32311}{}
\epsfig{file=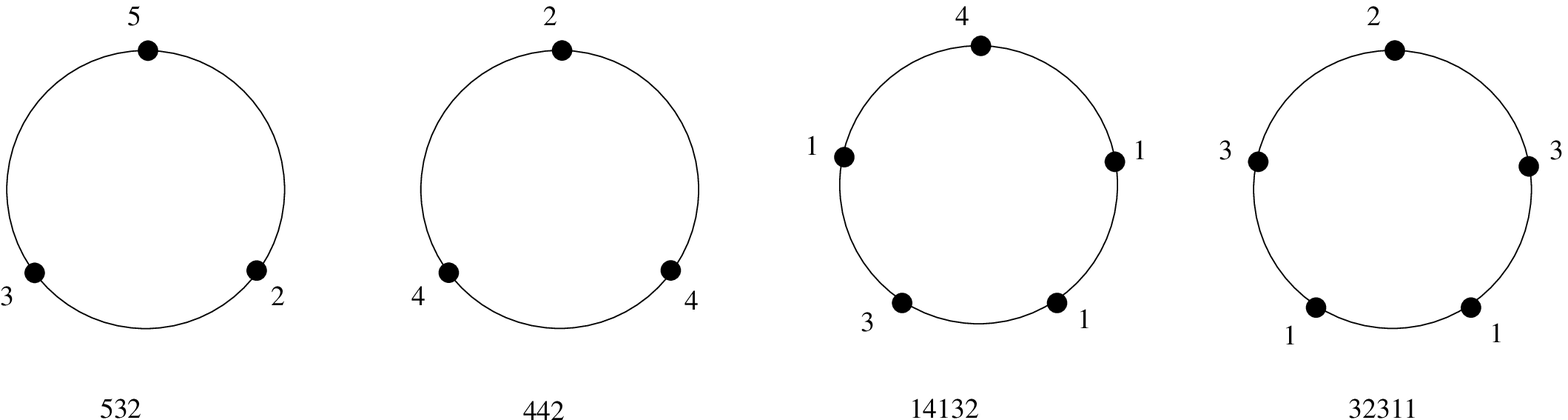,width=0.8\linewidth}
\end{center} 
%


By Lemma~\ref{lift}, $\Sigma$ contains either 0- or 1-lifting of one 
of these diagrams. By Lemma~\ref{no3-8}, $\Sigma$ contains no Lann\'er
diagram of order 3. Since $P$ is an 8-dimensional polytope, the
positive inertia index of $\Sigma$ is at most 8.
Hence, by Lemmas~\ref{442},~\ref{14131},~\ref{32311} and~\ref{1lift53},
$\Sigma$ contains $\Sigma(P_8)$,
which contradicts Lemma~\ref{same}.

\end{proof}

\section{Polytopes in dimension 7}
\label{dim7}
\begin{figure}[htb]
\begin{center}
\epsfig{file=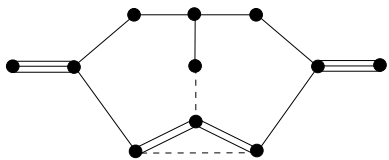,width=0.35\linewidth}
\caption{ A unique compact Coxeter 7-polytope with 11 facets.}
\label{p7}
\end{center}
\end{figure}

In this section we assume that $\Sigma$ is the Coxeter diagram of 
a compact Coxeter 7-polytope with $11$ facets 
and prove that $\Sigma$ coincides with
$\Sigma_{P_7}$, where $\Sigma_{P_7}$ is a diagram 
found by  V.~Bugaenko in~\cite{Bu1} and shown in Fig.~\ref{p7}

\begin{theorem}
\label{th7}
If  $\Sigma$ is the Coxeter diagram of 
a compact Coxeter $7$-polytope with $11$ facets 
then $\Sigma=\Sigma_{P_7}$.

\end{theorem}

The proof is similar to the proof in larger dimensions,
however it is rather long and routine, with many cases to consider.
First, we prove that $\Sigma$ contains a subdiagram of the type
$F_4$ or $H_4$, and then finish the proof by really boring
lemmas~\ref{f4,2sos}~--~\ref{h4,3sos} 
concerning the subdiagrams $F_4$ and $H_4$.

We say that a Coxeter diagram satisfies the {\it signature condition} 
if it is admissible and its positive
inertia index does not exceed 7. 

Recall that if $u,v\in \Sigma$ then $\[u,v\]=m$ ($\infty$ or
2) means that $u$ and $v$ are joined by an $(m-2)$-fold edge
(dotted or empty edge, respectively).

\subsection*{Existence of a subdiagram of the type $F_4$ or $H_4$}

In this subsection we show the following properties of $\Sigma$:

\begin{itemize}
\item
Any node of $\Sigma$ belongs to at most one dotted edge 
(Lemma~\ref{2dotted});
\item
$\Sigma$ contains no subdiagram of the type $G_2^{(k)}$ for $k>5$
(Lemma~\ref{multi-mult});
\item
$\Sigma$ contains at least one subdiagram of the type $F_4$ or $H_4$
(Lemma~\ref{f,h}),
and any such subdiagram has at least 2 bad neighbors
(Lemma~\ref{H_4,1sos}).
\end{itemize}

Recall from Lemma~\ref{bad} that any elliptic subdiagram of
$\Sigma$ has at most 3 bad neighbors.

\begin{lemma}
\label{g2,3badsos}
A subdiagram of the type $G_2^{(k)}$, $k>3$, has at most $2$ bad neighbors.

\end{lemma}

\begin{proof}
Suppose that $S_0\subset \Sigma$ is a 
subdiagram of the type $G_2^{(k)}$, $k>3$, with 3 bad neighbors.
Then $P(S_0)$ is a Coxeter 5-simplex, which is impossible.

\end{proof}

\begin{lemma}
\label{2dotted}
Any node of $\Sigma$ is incident to at most one dotted edge.

\end{lemma}

\begin{proof}
Suppose that a node $v$ is incident to at least two dotted edges.
Then the facet $f$ of $P$ corresponding to $v$ is a (possibly non-Coxeter)
6-polytope with at most $6+2$ facets.
By Cor.~\ref{simplex}, $f$ is not a simplex. 
Hence (by Prop.~\ref{prod}), it is a product of two
simplices, i.e. either $\Delta^5\times \Delta^1$, or  
$\Delta^4\times \Delta^2$, or  $\Delta^3\times \Delta^3$.
The first case is impossible, since $\Sigma$ contains no large missing
faces. The second and the third cases are impossible since the diagrams
$\lf 5,3 \rf_1$ and $\lf 4,4 \rf_1$ have no 0- and 1-liftings with positive
inertia index smaller than 8 (see Lemmas~\ref{lift44} and~\ref{1lift53}).

\end{proof}

\begin{lemma}
\label{H_4,1sos}
Any subdiagram of $\Sigma$ of the type $H_4$ or $F_4$ 
has at least two neighbors. 

\end{lemma}

\begin{proof}
Suppose that $S_0$ is a diagram of the type $H_4$ or $F_4$.
Since $\Sigma$ is a connected diagram,  $S_0$ has at least one neighbor. 
Suppose that $S_0$  has a unique neighbor, $a$. 
Then $P(S_0)$ is a Coxeter 3-polytope with
$3+3$ facets.
There are two simple combinatorial types of 3-polytopes with 6 facets,
namely,
a cube and a doubly truncated tetrahedron, i.e. a polytope with 2 
pentagonal, 2 quadrilateral and 2
triangular facets. The later case is impossible for $P(S_0)$, since
any of its
triangular facets does not meet two other facets in contradiction
to lemma~\ref{2dotted} (here we use that 
$\Sigma_{S_0}=\o S_0\subset \Sigma$, since $S_0$ is a diagram of the
type $F_4 $ or $H_4$). Hence, $P(S_0)$ is a cube.
Denote by $b_1$ and $b_2$, $c_1$ and $c_2$, $d_1$ and $d_2$ the
ends of the dotted edges in $\overline S_0=\Sigma_{S_0}$.
By Lemma~\ref{bad&Lanner}, $a$ is joined with each of the dotted 
edges $b_1b_2$, 
$c_1c_2$ and
$d_1d_2$. We assume that $a$ is joined with $b_1$, $c_1$ and $d_1$.

Suppose that $[b_1,c_1] \ge 4$. Then the subdiagram $\[b_1,c_1\]$
has at least 3 bad neighbors ($b_2$, $c_2$ and $a$), which is
impossible by Lemma~\ref{g2,3badsos}. 
Furthermore $[b_1,c_1] \ne \infty$, since $P(S_0)$ is a cube.
Hence,   $[b_1,c_1]=2$ or 3.
Similarly, $[b_1,d_1]\le 3$ and  $[c_1,d_1]\le 3$.
 
Suppose that  $[a,b_1]\ge 4$. Let $S_1=\[a,b_1\]$.
If $[a,b_1]\ge 6$, then $S_1$ has at least 3 bad neighbors
$b_2,c_1$ and $d_1$ in contradiction to Lemma~\ref{g2,3badsos}.
Hence,  $[a,b_1]= 4$ or $5$.
By Lemma~\ref{g2,3badsos}, the diagram
$S_1$ has at most 2 bad neighbors, thus, at least one of $c_1$ and
$d_1$ is a good neighbor of $S_1$ (recall that both  $c_1$ and
$d_1$ are joined with $a$). We may assume that $c_1$ is a good neighbor 
of $S_1$ and consider the diagram $S_2=\[S_1,c_1\]$ of the type
$H_3$ or $B_3$.
The diagram $S_2$ has at least four bad neighbors, 
namely $b_2,c_2,d_1$ and   
one  of the nodes of $S_0$ (since $a$ is a neighbor of $S_0$),
which is impossible by Lemma~\ref{bad}.
The contradiction shows that $[a,b_1]=3$ ($[a,b_1]\ne \infty$ 
by Lemma~\ref{2dotted}).
Similarly, $[a,c_1]=[a,d_1]=3$. 

Since  $\Sigma $ contains no parabolic subdiagrams and
$[b_1,c_1]\le 3$,
we obtain that  $[b_1,c_1]=2$. Similarly, $[b_1,d_1]=[c_1,d_1]=2$,  
 so the diagram $\[a,b_1,c_1,d_1\]$ is a diagram of the type
$D_4$. This diagram has at least four bad neighbors, $b_2,c_2,d_2$  and   
 one of $x_i$, $1\le i \le 4$. 
We come to a contradiction, and the lemma is proved.
 
\end{proof}

\begin{lemma}
\label{multi-mult}
$\Sigma$ contains no subdiagram $G_2^{(k)}$ for $k>5$.

\end{lemma}

\begin{proof}
Suppose that $\Sigma$ contains a subdiagram $S_0=\[x_1,x_2\]$ of the
type  $G_2^{(k)}$, $k>5$. Without loss of generality 
we may assume that the edge of $S_0$ has the maximal multiplicity
amongst all edges in $\Sigma$. Since $\Sigma$ is connected, 
$S_0$ has at least one (evidently bad) neighbor; 
by Cor.~\ref{bad},  
$S_0$ has at most 2 neighbors. 
Hence, $S_0$ has either 1 or 2 neighbors.
We consider these two cases.

\bigskip
\noindent
{\bf Case 1.}
Suppose that $S_0$ has a unique neighbor $a$.
Then $P(S_0)$ is a 5-polytope with $5+3$ facets.
Corollary~\ref{dif} implies that $\Sigma_{S_0}=\overline S_0$.
Hence, by Lemma~\ref{2dotted}, any node of $\Sigma_{S_0}$ is incident
to at most one dotted edge. The list of 5-polytopes with 8 facets
contains a unique entry satisfying this condition
(see Fig.~\ref{multi12}(a) for this diagram and notation for its nodes).
By Lemma~\ref{bad&Lanner},
$a$ is joined either with $z_1$ or with $z_2$,
say with $z_1$.
Denote $S_1=\[b_1,z_1\]$. If $a$ is a bad neighbor of $S_1$, then
$S_1$ has 3 bad neighbors ($y_1,z_2,a$) in contradiction 
to Lemma~\ref{g2,3badsos}.
Therefore, $a$ is a good neighbor of $S_1$,
$[a,z_1]=3$, 
and the diagram $S_2=\[a,S_1\]$ is a diagram of the type $B_3$.
Recall that $a$ is a neighbor of $S_0$, and we may assume that
$x_1$ is joined with $a$. If $x_1$ is a good neighbor of $S_2$, then
the diagram $\[x_1,S_2\]$ has more than 3 bad neighbors
(namely, $x_2,z_2,y_1$ and some node of the Lann\'er diagram
$\[b_2,y_4,y_3,y_2\]$ attached to $a$; 
the latter neighbor is a bad one since $a$ is not a leaf of  
$\[x_1,S_2\]$). This is impossible by  Lemma~\ref{g2,3badsos}, 
so $x_1$ is a bad
neighbor of $S_2$ and $S_2$ has 3 bad neighbors ($x_1,y_1,z_2$). 
Thus, 
$P(S_2)$ is a 4-dimensional simplex and 
$\Sigma_{S_2}$ is a Lann\'er diagram of order 5.
By Cor.~\ref{dif2}, this implies that $\o S_2$ is a Lann\'er diagram, too. 
However,
$\o S_2=\[x_2,b_2,y_4,y_3,y_2\]$ cannot be a Lann\'er diagram,
since it contains a Lann\'er diagram $\[b_2,y_4,y_3,y_2\]$ 
of order 4.
The contradiction shows that $S_0$ has 2 neighbors.   

\bigskip
\noindent
{\bf Case 2.}
Suppose that the diagram $S_0=\[x_1,x_2\]$ has 2 neighbors $a_1$ and $a_2$.
Then $P(S_0)$ is a 5-prism. We have two possibilities for this prism
shown in Fig.~\ref{multi12}(b),(c). We denote the nodes of $\overline S_0$ as
it is shown in the figure.
Denote by $S_1$ the subdiagram $\[y_1,y_2,y_3,y_4\]$ of the type $F_4$
or $H_4$. 
Since the subdiagram
$\overline S_1=\Sigma_{S_1}$ contains the diagram
$S_0$ of the type $G_2^{(k)}$, $k>5$, the diagram $\o S_1$
is not a Lann\'er diagram of order 4, so the face $P(S_1)$ is not a
3-simplex. Hence, $S_1$ has at most two neighbors, and at least one of 
$a_1$ and $a_2$ is not a neighbor of $S_1$. 
We may assume that $a_1$
is not a neighbor of $S_1$. By Lemma~\ref{bad&Lanner}, 
this means that $a_1$ is joined with $y_5$
(since $y_5$ belongs to a Lann\'er diagram 
$\[y_1,y_2,y_3,y_4,y_5\]\subset \o S_0$). 
We consider cases $S_1=F_4$ and $S_1=H_4$ separately.

\begin{figure}[!h]
\begin{center}
\psfrag{a}{$\t a$}
\psfrag{y1}{$\t y_1$}
\psfrag{y2}{$\t y_2$}
\psfrag{y3}{$\t y_3$}
\psfrag{y4}{$\t y_4$}
\psfrag{y5}{$\t y_5$}
\psfrag{z1}{$\t z_1$}
\psfrag{z2}{$\t z_2$}
\psfrag{2,3}{{\scriptsize  $2,3$}}
\psfrag{3,4}{{\scriptsize  $3,4$}}
\psfrag{2,3,4}{{\scriptsize  $2,3,4$}}
\psfrag{(a)}{(a)}
\psfrag{(b)}{(b)}
2\psfrag{(c)}{(c)}
\psfrag{b1}{$\t z_1$}
\psfrag{c1}{$\t b_1$}
\psfrag{d1}{$\t y_1$}
\psfrag{e1}{$\t y_2$}
\psfrag{b2}{$\t z_2$}
\psfrag{c2}{$\t b_2$}
\psfrag{d2}{$\t y_4$}
\psfrag{e2}{$\t y_3$}
\epsfig{file=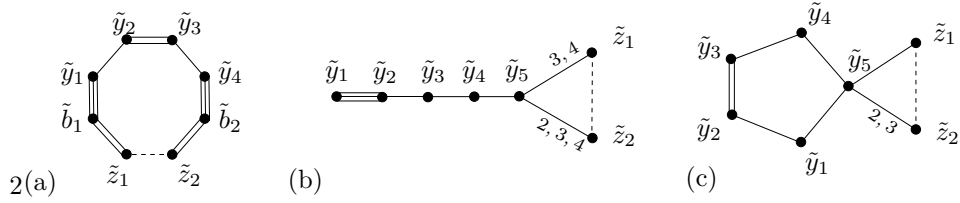,width=0.85\linewidth}
\caption{To the proof of Lemmas~\ref{multi-mult} and~\ref{f4,2sos}}
\label{multi12}
\end{center}
\end{figure}

In the case $S_1=F_4$ consider two subdiagrams
$\[y_2,y_3,y_4,y_5,z_1\]$ and $\[y_3,y_2,y_1,y_5,z_1\]$ of the type
$B_5$. Each of these diagrams has 3 bad neighbors ($a_1,y_1,z_2$
and $a_1,y_4,z_2$ respectively), so $a_2$ is not a bad neighbor for
these diagrams. Therefore, $a_2$ is not joined with 
a Lann\'er diagram $\[S_1,y_5\]$,
which contradicts  Lemma~\ref{bad&Lanner}. 

Now we are left with the case $S_1=H_4$.
Since $a_1$ is not a neighbor of $S_1$, 
the node $a_2$ is a neighbor of $S_1$, otherwise $S_1$ has a
unique neighbor in contradiction to Lemma~\ref{H_4,1sos}.
Therefore, $S_1$ has 2 bad neighbors, $y_5$ and $a_2$, 
and $P(S_1)$ is a 3-prism.
This means that the diagram $\o S_1=\Sigma_{S_1}$ consists of a dotted
edge $z_1z_2$ and a Lann\'er diagram of order 3 $\[S_0,a_1\]$.
Thus, the subdiagram
$X= \Sigma\setminus \[z_1,z_2\]=\[S_0,a_1,S_1,y_5\]$ consists of two   
Lann\'er diagrams $\[S_0,a_1\]$ and $\[ S_1,y_5\]$ 
joined by a unique edge $a_1y_5$. 
If this edge is not dotted one, then the    
subdiagram $X$ is superhyperbolic, unless $[a_1,y_5]=5$. 
We consider two cases,  $[a_1,y_5]=5$ and  $[a_1,y_5]=\infty$.

\medskip
\noindent
{\bf Case 2.1.}
Suppose that $[a_1,y_5]=5$.
Then the subdiagram $\[a_1,y_5,y_4,y_3\]$ of the type $H_4$ has 
3 bad neighbors, namely $z_1,y_2$ and one of $x_1$ and $x_2$, say $x_1$.
Hence, $z_2$ is not a neighbor of  $\[a_1,y_5,y_4,y_3\]$.
In particular, $[z_2,a_1]=2$, which means (Lemma~\ref{bad&Lanner}) that
$[a_1,z_1]\ne 2$. Therefore, $z_1$ is a bad neighbor of
$S_2=\[a_1,y_5\]$, and $\Sigma_{S_2}$ is a diagram of a 5-polytope
with at most 8 facets. Consider  $\Sigma_{S_2}$. 
By Prop.~\ref{al}, the subdiagram 
$S_1 =\[y_1,y_2,y_3,y_4\]$ turns into a linear diagram of order 4 
with a triple edge $\t y_1 \t y_2$, simple edge $\t y_2 \t y_3$ 
and an edge $\t y_3 \t y_4$ labeled by 10 in  $\Sigma_{S_2}$.
However, no  diagram of 5-polytope
with at most 8 facets contains such a subdiagram.
Thus, the case  $[a_1,y_5]=5$ is impossible.

\medskip
\noindent
{\bf Case 2.2.}
Suppose that  $[a_1,y_5]=\infty$.
Consider the subdiagram $S_2=\[y_1,y_2,y_3\]$ of the type $H_3$.
If $a_2$ is a bad neighbor of $S_2$ then $P(S_2)$ is a 4-polytope
with $4+3=7$ facets, and $\Sigma_{S_2}$ contains a subdiagram
$G_2^{(k)}$ for $k>5$ and at least 2 dotted edges
($\t z_1 \t z_2$ and $\t y_4 \t y_5$). The list of
4-polytopes with $7$ facets contains no entry with these properties.
Hence, $a_2$ is not a bad neighbor of $S_2$.
If $a_2$ is a good neighbor of $S_2$ then the diagram
$\[y_1,y_2,y_3,a_2\]$ of the type $H_4$ has at least 4 bad neighbors
(namely, $y_4$ and at least one of $x_1$ and $x_2$, one of $z_1$ and
$z_2$, and one of $a_1$ and $y_5$).
Therefore, $a_2$ is not a neighbor of $S_2$, and hence, 
by Lemma~\ref{H_4,1sos}, $a_2$ is joined with $y_4$. 
Consider the subdiagram $S_3=\[y_2,y_3,y_4,y_5,z_1\]$
of the type $A_5$ or $B_5$. $S_3$ has three bad neighbors ($y_1,a_1,z_2$),
so $a_2$ is a good neighbor of $S_3$ and $\[a_2,S_3\]$ is a
diagram of the type $E_6$ having at least four 
bad neighbors ($y_1,z_2,a_1$ and one of $x_1$ and $x_2$). 
This contradiction proves the lemma.

\end{proof}

\begin{lemma}
\label{f,h}
$\Sigma$ contains at least one subdiagram of the type $F_4$ or $H_4$. 

\end{lemma}

\begin{proof}
Suppose that the lemma is broken, i.e. 
$\Sigma$ contains no subdiagram of the types 
$F_4$ and $H_4$. 

Suppose that $\Sigma$ contains 
 a subdiagram 
$S_0=\[x_1,x_2\]$ of the type $G_2^{(4)}$ or $G_2^{(5)}$, having a bad
neighbor. Then, $P(S_0)$ is a 5-polytope with at most $5+3$ facets.
Therefore, $\Sigma_{S_0}$ contains a subdiagram
of the type  either $F_4$ or $H_4$. Cor.~\ref{dif2} implies that
$\o S_0$ contains a subdiagram of this type,
which contradicts  the assumption.

In particular, we conclude that $\Sigma$ contains no Lann\'er diagram of
order 3 (we use also Lemma~\ref{multi-mult}). 
Since any Lann\'er diagram of order 5 contains a subdiagram of type $F_4$
or $H_4$,
$\Sigma$ contains no Lann\'er diagram of order 5 either.
By \cite[Satz~6.9]{Ess2},  any simple $d$-polytope ($d>4$) with $d+4$ facets
contains at least one missing face of order greater than 2.
Thus,
$\Sigma$ contains a Lann\'er subdiagram $L$ of order 4. 
Let $S_0\subset L$ be a subdiagram of the type $H_3$ or $B_3$.
Denote $a_1=L\setminus S_0$ and $S_0=\[x_1,x_2,x_3\]$. 

Suppose that $S_0$ has three bad neighbors.
Then $P(S_0)$ is a 4-simplex, so $\Sigma_{S_0}$ (and, hence, $\o S_0$)
contains a subdiagram of the type either $F_4$ or $H_4$,
which is impossible by the assumptions.

Suppose that $S_0$ has a unique bad neighbor, $a_1$.
Then $P(S_0)$ is a 4-polytope with $4+3$ facets. 
It follows from the assumption of the lemma that
$\Sigma_{S_0}$ is a diagram of  a 4-polytope with $4+3$
facets containing neither a subdiagram of type $F_4$ nor a subdiagram of
the type $H_4$, and containing no subdiagram of the type  $G_2^{(k)}$,
$k>5$ (Lemma~\ref{multi-mult}).
Each of these diagrams contains a subdiagram 
$S_1=\[y_1,y_2\]$ of the type $G_2^{(4)}$ or $G_2^{(5)}$ having a bad
neighbor. 
If $S_0$ is a diagram of the type $H_3$, then it has no good neighbor
(by the assumption $\Sigma$ contains no subdiagram of the type $H_4$), 
and $\o S_0=\Sigma_{S_0}$   
contains a subdiagram 
$\[\t y_1,\t y_2\]$ of the type $G_2^{(4)}$ or $G_2^{(5)}$ having a bad
neighbor. As it is shown above, this is impossible.
Thus, $S_0$ is a diagram of the type $B_3$. By Cor.~\ref{dif},
either $S_1=\[y_1,y_2\]\subset \o S_0$ is a multiple edge 
with a bad neighbor (which is impossible)
or we have one of the following possibilities for the diagram
$\Sigma_{S_0}$:

\begin{center}
\psfrag{a1}{$\t a_1$}
\psfrag{x1}{$\t x_1$}
\psfrag{x2}{$\t x_2$}
\psfrag{x3}{$\t x_3$}
\psfrag{y1}{$\t y_4$}
\psfrag{y2}{$\t y_3$}
\psfrag{y3}{$\t y_2$}
\psfrag{y4}{$\t y_1$}
\psfrag{y5}{$\t y_5$}
\psfrag{y6}{$\t y_6$}
\epsfig{file=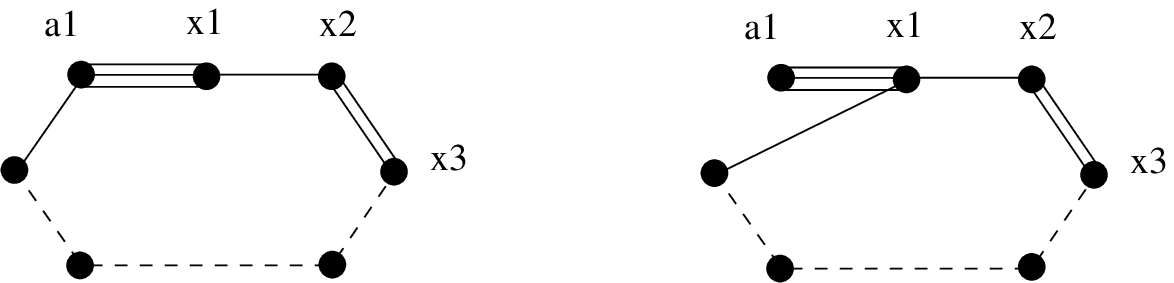,width=0.4399\linewidth}

\end{center}

In this case the double edge $\t y_1 \t y_2$ may turn into a simple edge
 $y_1y_2$ in $\o S_0$, 
that leads to a subdiagram $\[y_1,y_2,y_3,y_4\]\subset \o S_0$ of the
type $H_4$, which contradicts  the assumption.
Therefore, the multiple edge $\t y_1\t y_2$ remains multiple 
edge $y_1y_2$ in $\o
S_1$, and Corollary~\ref{dif} implies that the bad neighbor of $S_1$
remains bad in $\o S_0$, which is impossible. 
So, $S_0$ cannot have 3 bad neighbors.
 
Therefore, $S_0$ has exactly 2 bad neighbors, $a_1$ and $a_2$, 
and $P(S_0)$ is either an Esselmann polytope or a 4-prism.
Since $\Sigma$ contains no diagram of the types $F_4$ and $H_4$,
using Cor.~\ref{dif2} we obtain that
$\Sigma_{S_0}$ is the following diagram: 
%
%
\begin{center}
\psfrag{y1}{$\t y_1$}
\psfrag{y2}{$\t y_2$}
\psfrag{y3}{$\t y_3$}
\psfrag{y4}{$\t y_4$}
\psfrag{y5}{$\t y_5$}
\psfrag{y6}{$\t y_6$}
\psfrag{2,3}{{\scriptsize $2,3$}}
\epsfig{file=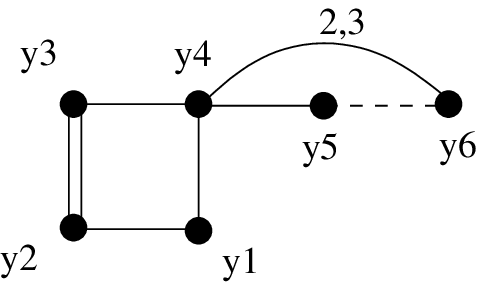,width=0.2399\linewidth}
\end{center}
%
%
%
%
By Cor.~\ref{dif2},
the nodes of $\Sigma_{S_0}$
(with possible exclusion for $y_6$ in case of $[y_4,y_6]=2$) cannot be good neighbors
of $S_0$. In particular, $\o S_0$ contains a cyclic Lann\'er diagram
of order 4 with a unique double edge, and contains no subdiagram of
the type $H_3$.
Furthermore, it is easy to check that a dotted edge $\t y_5\t y_6$ 
of $\Sigma_{S_0}$ corresponds to a dotted edge $y_5y_6$ in $\o S_0$ 
(indeed, otherwise $[y_4,y_6]=2$, and if the edge 
$y_5y_6$ is triple one, then
 $\[y_3,y_4,y_5,y_6\]$ is a subdiagram of 
the type $H_4$, which contradicts  the assumption; if $y_5y_6$
is a double edge, then $\[y_2,y_3,y_4,y_5,y_6\]$ is a
parabolic subdiagram of $\Sigma$ of the type $\w C_4$, which is also 
impossible; if  $y_5y_6$ is a simple edge,
then $\[y_5,y_6,S_0\]$ is a diagram of the type $B_5$
and Prop.~\ref{al} shows that $\t y_5 \t y_6$ must not be a dotted
but a  double edge of $\Sigma_{S_0}$).
So $\o S_0$ consists of a dotted edge and of a cyclic Lann\'er 
diagram of order 4, and $\o S_0=\Sigma_{S_0}$. 

Consider the diagram $S_1=\[y_1,y_2,y_3\]$. It is a subdiagram of the
type $B_3$ in the Lann\'er diagram   $\[y_1,y_2,y_3,y_4\]$ of order 4.
The same reasoning as for $\o S_0$ shows that
the diagram $S_1$ has exactly 2 bad neighbors ($y_4$ and one of $a_1$
and $a_2$), and
the diagram $\o S_1$ consists of a dotted edge $y_5y_6$
and a cyclic Lann\'er diagram (either $\[S_0,a_2\]$ or 
$\[S_0,a_1\]$ respectively).
Without loss of generality we may assume that
$a_2$ is a bad neighbor of $S_1$ (and $\[a_1,S_0\]$ is a cyclic
Lann\'er diagram). By Cor.~\ref{dif2}, $a_1$ is not a good neighbor of $S_1$,
so $a_1$ is joined with $y_4$ (Lemma~\ref{bad&Lanner}). 
If $[a_1,y_4]= 3$ or 4 then
the subdiagram $\[x_1,x_2,a_1,y_2,y_3,y_4\]$ contains a parabolic
subdiagram of the type $\w B_5$  or $\w B_3$. If  $[a_1,y_4]=5$ then
$\[x_2,x_3,a_1,y_4\]$ is a  subdiagram of the type $H_4$.
Therefore, $[a_1,y_4]=\infty$.

Consider the diagram $S_2=\[y_2,y_3,y_4\]$. It has two bad neighbors,
$a_1$ and $y_1$, so $a_2$ is not a bad neighbor of $S_2$.
Therefore, $a_2$ is joined with $y_1$, and repeating the reasoning
above we obtain that $[a_2,y_1]=\infty$.
Therefore, the diagram $\[y_1,y_3,y_4,y_5\]$ of the type $D_4$ has 4
bad neighbors $a_1,a_2,y_2,y_6$, which is impossible.

\end{proof}

\subsection*{Treatment of the subdiagrams $F_4$ and $H_4$}

We have proved that $\Sigma$ contains at least one subdiagram of the
type $F_4$ or $H_4$. 
In Lemmas~\ref{f4,2sos} and~\ref{h4,2sos}
we prove that such a subdiagram always has 3 bad neighbors.
Next, in Lemma~\ref{f4,3sos} we show that
  $\Sigma$ contains no subdiagram of the type $F_4$.
Lemma~\ref{h4,3sos} finishes the proof of 
Theorem~\ref{th7}. 
Notice, that the proofs of Lemmas~\ref{h4,2sos} and~\ref{h4,3sos}
(concerning the subdiagrams of the type $H_4$) turn out to be
much more complicated than the proofs of the similar 
Lemmas~\ref{f4,2sos} and~\ref{f4,3sos}
(concerning the subdiagrams of the type $F_4$). 
The possible reason is that a diagram of the type $H_4$
appears in diagrams of $d$-polytopes with at most $d+3$
facets more often than  a diagram of the type $F_4$ does.

\begin{lemma}
\label{f4,2sos}
Any subdiagram of $\Sigma$ of the type $F_4$ 
has three neighbors.

\end{lemma}

\begin{proof}
Suppose the contrary. By Lemma~\ref{H_4,1sos},
this implies that $\Sigma$ contains a subdiagram $S_0$ of the type
$F_4$ with 2 neighbors. Then $P(S_0)$ is a 3-prism. The diagram
$\Sigma_{S_0}=\overline S_0$ consists of a dotted edge $z_1z_2$ and
 a Lann\'er diagram $\[x_1,x_2,x_3\]$ of order
3, which in its
turn contains some multiple edge. Choose in $L=\[x_1,x_2,x_3\]$
an edge  $S_1=\[x_1,x_2\]$ of maximal possible multiplicity in $L$. 
The diagram $S_1$ has at least one bad neighbor, $x_3$.
By Lemma~\ref{g2,3badsos}, $S_1$ has either 1 or 2 bad neighbors.

Suppose that $S_1$ has a unique bad neighbor $x_3$.
Then $P(S_1)$ is a 5-polytope with $5+3=8$ facets. 
The diagram $S_0$ of the type $F_4$ is not joined with $S_1$.
Hence, the diagram $\Sigma_{S_1}$
contains a subdiagram of the type $F_4$, and $\Sigma_{S_1}$
is the diagram shown in Fig.~\ref{multi12}(a). 
We denote the nodes of $S_0\subset \o S_1$ by $y_1,y_2,y_3,y_4$ and 
denote the neighbors of $S_0$ by $b_1$ and $b_2$.
By Cor.~\ref{dif2},  
the nodes $b_1$ and $b_2$ are not
neighbors of $S_1$. A Lann\'er diagram $\[S_1,x_3\]$ should be
connected with a Lann\'er diagram $\[b_1,y_1,y_2,y_3\]$, thus,
$x_3$ is joined with $b_1$. Similarly, $x_3$ is joined with $b_2$.
Consider the diagram $S_2=\[b_1,y_1,y_2\]$ of the type $H_3$.
It has 3 bad neighbors ($y_3,x_3$ and $z_1$) 
and has no good neighbor. Therefore, $P(S_2)$
is a 4-simplex and $\o S_2=\Sigma_{S_2}=\[y_4,b_2,z_2,x_1,x_2\]$ is a
Lann\'er diagram of order 5.
Since $y_4 b_2$ and $x_1x_2$ are
two unjoined multiple edges, the diagram $\o S_2$ is a linear Lann\'er
diagram (with edges labeled by 5,3,3,4 or by 5,3,3,5).
Since $[y_4,b_2]=5$,
the diagram $\[y_4,b_2,z_2,x_1\]$
(or $\[y_4,b_2,z_2,x_2\]$) is a diagram of the type $H_4$.
In particular, $[b_2,z_2]=3$ in $\Sigma$. 
So, a simple edge $b_2z_2$ of $\Sigma$ turns into a double edge
$\t b_2 \t z_2$ in $\Sigma_{S_1}$.
In view of Cor.~\ref{dif},
this implies that $[x_1,x_2]=4$. Since the multiplicity of
$x_1x_2$ is maximal in the Lann\'er diagram $\[x_1,x_2,x_3\]$,  
we obtain that $x_3$ is joined with both $x_1$ and $x_2$.
Therefore, the diagram  $\[y_4,b_2,z_2,x_1\]$ (or
$\[y_4,b_2,z_2,x_2\]$)
of the type $H_4$ 
has at least 4 neighbors: $y_3,x_2$ (or $x_1$),
$x_3$  and $z_1$, which is impossible. 

The contradiction shows that $S_1$ has 2 bad neighbors, $x_3$ and
some node $a_1$. 
So, $P(S_1)$ is a 5-prism containing a subdiagram of the type $F_4$ 
(see Fig.~\ref{multi12}(c) for the notation).
Recall that $a_1\notin \o S_1$, and hence, $a_1$ is a bad neighbor of 
$S_0=\[y_1,y_2,y_3,y_4\]$.  
Consider the diagrams $\[y_2,y_3,y_4,y_5,z_1\]$ and
$\[y_3,y_2,y_1,y_5,z_1\]$ of the type $B_5$. 
The node $a_1$ is a bad neighbor for at
least one of these diagrams, say  $\[y_2,y_3,y_4,y_5,z_1\]$. Then
this diagram  $\[y_2,y_3,y_4,y_5,z_1\]$  has at least 4 bad neighbors
($y_1,a_1,z_2$ and one of $x_1,x_2,x_3$).
The contradiction proves the lemma.
 
\end{proof}

\begin{lemma}
\label{h4,2sos}
Any subdiagram of $\Sigma$ of the type $H_4$ 
has three neighbors. 

\end{lemma}

\begin{proof}
Suppose the contrary. By Lemma~\ref{H_4,1sos},
this implies that $\Sigma$ contains a subdiagram $S_0$ of the type
$H_4$ with 2 neighbors. Then $P(S_0)$ is a 3-prism. The diagram
$\Sigma_{S_0}=\overline S_0$ consists of a dotted edge $z_1z_2$ and
 a Lann\'er diagram $L=\[y_1,y_2,y_3\]$ of order
3, which in its
turn contains some multiple edge. We suppose that $y_1y_2$ 
is an edge of  maximal possible multiplicity in $L$, 
 and let $S_1=\[y_1,y_2\]$. 
By Lemma~\ref{g2,3badsos}, $S_1$ has either 1 or 2 neighbors.

\bigskip
\noindent
{\bf Case 1.}
Suppose that $S_1$ has a unique bad neighbor $y_3$.
Then $P(S_1)$ is a 5-polytope with $5+3=8$ facets. Since $\Sigma_{S_1}$
contains a subdiagram $S_0$ of the type $H_4$, the diagram $\Sigma_{S_1}$
is the diagram of one of three types (a), (b) and (c)
shown in Fig.~\ref{h4,2sos,1}. 
We denote the nodes of $S_0\subset \o S_1$ by $x_1,x_2,x_3,x_4$ and 
denote the neighbors of $S_0$ by $a_1$ and $a_2$ 
(see  Fig.~\ref{h4,2sos,1}). We consider cases (a), (b) and (c)
separately.

\begin{figure}[!h]
\begin{center}
\psfrag{x1}{$\t x_1$}
\psfrag{x2}{$\t x_2$}
\psfrag{x3}{$\t x_3$}
\psfrag{x4}{$\t x_4$}
\psfrag{z1}{$\t z_1$}
\psfrag{z2}{$\t z_2$}
\psfrag{a1}{$\t a_1$}
\psfrag{a2}{$\t a_2$}
\psfrag{2,3,4}{{\scriptsize $2,3,4$}}
\psfrag{(a)}{(a)}
\psfrag{(b)}{(b)}
\psfrag{(c)}{(c)}
\epsfig{file=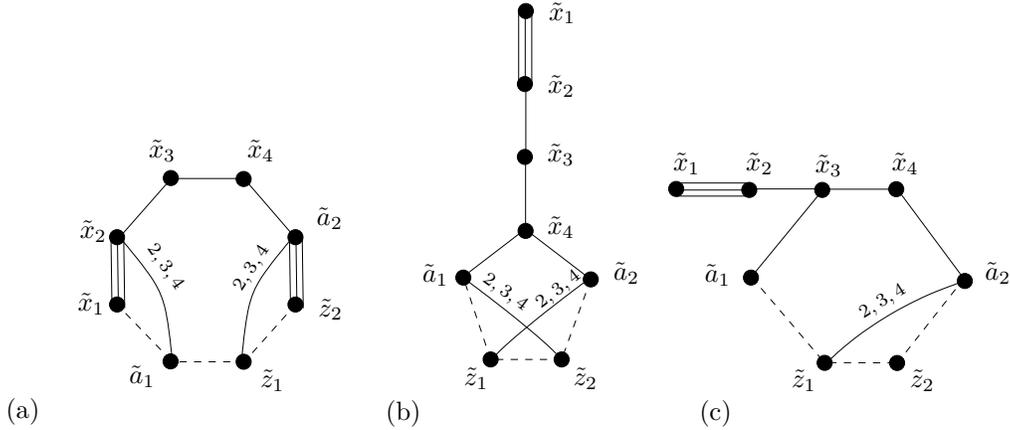,width=0.92\linewidth}
\caption{Possibilities for $\Sigma_{S_1}$ (Case 1).}
\label{h4,2sos,1}
\end{center}
\end{figure}

\medskip
\noindent
{\bf Case 1.1.}
Suppose that  $\Sigma_{S_1}$ is the diagram shown in  
Fig.~\ref{h4,2sos,1}(a).
Recall that $z_1z_2$ is a dotted edge in $\Sigma$,
so Lemma~\ref{2dotted} implies that $a_1$ is a good neighbor of $S_1$.
Without loss of generality we may assume that $[a_1,y_1]=3$ and 
$[a_1,y_2]=2$.
On the other hand, $a_2$ is not a good neighbor of $S_1=\[y_1,y_2\]$,
since $[a_2,z_2]=5$.
Hence, $a_2$ is joined with $y_3$ 
(otherwise a Lann\'er diagram $\[x_1,x_2,x_3,x_4,a_2\]$ is not joined with
a Lann\'er diagram $L=\[S_1,y_3\]$). 
Notice that the diagram $S_2=\[z_2,a_2,x_4,x_3\]$ of the type $H_4$
has 3 neighbors ($x_2,z_1,y_3$), hence, $P(S_2)$ is a 3-simplex,
so the diagram $\o S_2=\[x_1,a_1,y_1,y_2\]$ 
is a Lann\'er diagram.
Furthermore, $\o S_2$ is a linear Lann\'er diagram ($x_1$ is not joined
with $\[y_1,y_2\]$, and $a_1$ is not joined with $y_2$ and joined with
$y_1$ by a simple edge). 
Therefore, $[x_1,a_1]=4$ or $5$. Consider the diagram
$S_3=\[x_1,a_1,y_1\]$ of the type $H_3$ or $B_3$. 
It has 3 bad neighbors
($x_2,z_1,y_2$), hence $P(S_3)$ is a 4-simplex, and
$\Sigma_{S_3}$ is a Lann\'er diagram of order 5.
By Cor.~\ref{dif2}, this implies that $\o S_0$ is a Lann\'er diagram, too.
However, in $\o S_3=\[x_3,x_4,a_2,z_2,y_3\]$ one end ($a_2$)
of the triple edge has valency at least three,
which is impossible in a Lann\'er diagram of order 5.
    
\medskip
\noindent
{\bf Case 1.2.}
Suppose that  $\Sigma_{S_1}$ is one of the diagrams shown in  
Fig.~\ref{h4,2sos,1}(b) and (c).
In these cases neither $a_1$ nor $a_2$ can be a good neighbor of
$S_1$, and Lemma~\ref{2dotted} implies that both $z_1$ and $z_2$ are
good neighbors of $S_1$ (recall that $z_1z_2$ is a dotted edge in
$\Sigma$). 
By Lemma~\ref{multi-mult}, $[a_i,z_i]=3$, 4 or 5 ($i=1,2$). If $[a_1,z_1]=5$, then the diagram
$\[z_1,a_1,x_4,x_3\]$ of the type $H_4$
has at least four neighbors ($x_2,z_2,a_2$ and at least one of $y_1$
and $y_2$). 
 If $[a_1,z_1]=4$, then the diagram
$\[z_1,a_1,x_4,x_3,x_2\]$
of the type $B_5$
has at least four bad neighbors ($x_1,z_2,a_2$ and at least one of
$y_1$ and $y_2$). 
Finally, if  $[a_1,z_1]=3$, then the edge $\t a_1\t z_1$ of 
 $\Sigma_{S_1}$ should be either a double edge  
or an edge marked by $10$, but not a dotted edge.  

\bigskip
\noindent
{\bf Case 2.}
Suppose that 
$S_1$ has 2 bad neighbors: $y_3$ and one of $a_1$ and $a_2$, say
$a_1$. Then, $P(S_1)$ is a 5-prism containing a subdiagram of the type
$H_4$, so, the diagram $\Sigma_{S_1}$ is the following one:

\begin{center}
\psfrag{y1}{$\t x_1$}
\psfrag{y2}{$\t x_2$}
\psfrag{y3}{$\t x_3$}
\psfrag{y4}{$\t x_4$}
\psfrag{z2}{$\t z_1$}
\psfrag{z1}{$\t z_2$}
\psfrag{y5}{$\t a_2$}
\psfrag{3,4}{{\scriptsize $3,4$}}
\psfrag{2,3,4}{{\scriptsize $2,3,4$}}
\epsfig{file=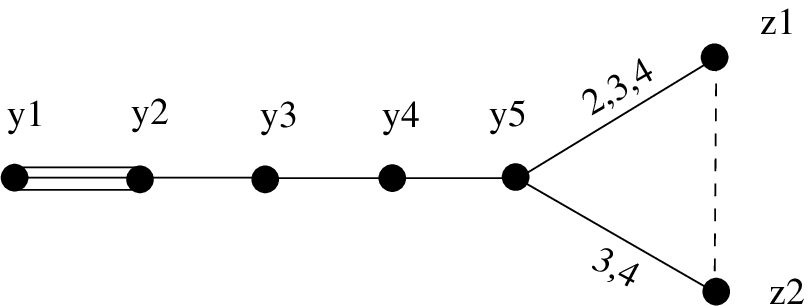,width=0.32\linewidth}
\end{center}

By Cor.~\ref{dif2}, $a_2$ 
is not a good neighbor of $S_1$, hence, $a_2$ is a neighbor of
$y_3$ (see Lemma~\ref{bad&Lanner}). 
If $[a_2,y_3]\ne \infty$, then the subdiagram 
$\[S_0,a_2,y_3,S_1\]$ is superhyperbolic (in assumption
that $\Sigma$ contains no subdiagrams of the type $G_2^{(k)}$ for
$k>5$
and
the multiplicity of the edge $y_1y_2$ is maximal possible in
$\[y_1,y_2,y_3\]$).
Hence, we assume  $[a_2,y_3]= \infty$.
Let $S_4=\[x_1,x_2,x_3\]$. Consider 3 cases:
$a_1$ is either a bad neighbor of $S_2$, or a good neighbor, or a
non-neighbor.

\medskip
\noindent
{\bf Case 2.1.}
Suppose that $a_1$ is a bad neighbor of $S_4$. 
Then $P(S_4)$ is a 4-polytope
with $4+3$ facets, and $\Sigma_{S_4}$ contains at least 3 dotted edges
$\t z_1 \t z_2$, $\t a_2\t y_3$ and $\t x_4\t a_2$. 
However, no Coxeter diagram of 4-polytope with $4+3$ facets
contains more than 3 dotted edges and
if a Coxeter
diagram of a 4-polytope with $4+3$ facets have exactly 3 dotted edges 
then any dotted edge is incident to some other dotted edge.
The edge $\t z_1 \t z_2$ has no common vertices with
$\t a_2\t y_3$ and $\t x_4\t a_2$, so we obtain a contradiction.

\medskip
\noindent
{\bf Case 2.2.}
Suppose that $a_1$ is a good neighbor of $S_4=\[x_1,x_2,x_3\]$. 
Consider the diagram
$S_5=\[x_1,x_2,x_3,a_2\]$ of the type $H_4$.
$S_5$ has 3 bad neighbors: $x_4$, one of $y_1$ and $y_2$ (say $y_1$), 
and one of
$z_1$ and $z_2$. This implies that $y_2$ is not a neighbor of $S_5$,
in particular, $[a_1,y_2]=2$ and $[a_1,x_1]=2$.
On the other hand, $a_1$ is a  bad neighbor of $\[y_1,y_2\]$.
Therefore, $[a_1,y_1]\in \{4,5,\infty \}$.

Consider the subdiagram $S_6=\[x_2,x_3,x_4,a_2,z_1\]$ of the type
either $A_5$ or $B_5$. Since $a_1$ is a neighbor of 
$S_0=\[x_1,x_2,x_3,x_4\]$ and as we have shown above, $[a_1,x_1]=2$,
we obtain that $a_1$ is a neighbor of $S_6$.
The diagram $S_6$ has 3 bad neighbors: $x_1,z_2,y_3$.
Hence, $a_1$ is a good neighbor of $S_6$. Therefore, 
$\[a_1,S_6\]$ is a diagram of the type $D_6$ with 3 bad neighbors 
($x_1,z_2,y_3$).
However, since $[a_1,y_1]=4$, 5 or $\infty$, $y_1$ is also a bad
neighbor of this diagram, which is impossible. 

\medskip
\noindent
{\bf Case 2.3.}
Suppose that $a_1$ is not a neighbor of $S_4=\[x_1,x_2,x_3\]$. 
Then $a_1$ is joined with $x_4$ (since $a_1$ is a neighbor of 
$S_0=\[x_1,x_2,x_3,x_4\]$).
Consider the subdiagram $S_7=\[x_2,x_3,x_4,a_2,z_1\]$ of the type
$A_5$ or $B_5$. $S_7$ has 3 bad neighbors: $x_1,z_2,y_3$. 
Thus, $a_1$ is a good neighbor of $S_7$,
and $\[S_7,a_1\]$ is a diagram of the type $E_6$.
However, this diagram $\[S_7,a_1\]$ has 4 bad neighbors:
$x_1,z_2,y_3$ and one of $y_1$ and $y_2$ (since $a_1$ is a bad
neighbor of $S_1$ by assumption). 
  
\medskip
\noindent
This shows that the assumption that
$S_1$ has 2 bad neighbors is impossible, which completes the proof of
the lemma. 

\end{proof}

\begin{lemma}
\label{f4,3sos}
$\Sigma$ contains no subdiagram of the type $F_4$.

\end{lemma}

\begin{proof}
Suppose that $S_0=\[x_1,x_2,x_3,x_4\]$ 
is a subdiagram of $\Sigma$ of the type $F_4$.
By Lemma~\ref{f4,2sos}, $S_0$ has 3 neighbors
$a_1,a_2$ and $a_3$. Then $P(S_0)$ is a 3-simplex, and 
$\o S_0=\Sigma_{S_0}=\[y_1,y_2,y_3,y_4\]$ is a Lann\'er diagram of order 4.
Let $S_1=\[y_1,y_2,y_3\]$ be a subdiagram of $\o S_0$ 
of the type $H_3 $ or $B_3$.
$S_1$ has at least one bad neighbor, $y_4$.
Let $S_2$ be a subdiagram of $\o S_0$ of the type $H_3 $ or $B_3$
containing $y_4$ (it does exist by Lemma~\ref{2subd_of_L4}),
we may assume that $S_2=\[y_2,y_3,y_4\]$.

The diagram $S_1$ has at least one bad neighbor, $y_4$.
We consider three cases where $S_1$ has 1, 2 or 3 bad neighbors
respectively.

\bigskip
\noindent
{\bf Case 1.} 
Suppose that $S_1$ has a unique bad neighbor $y_4$.
Then $P(S_1)$ is a 4-polytope with $4+3$ facets.

Suppose in addition that the diagram  $\[a_i,S_0 ,\o S_0\]$
contains a dotted edge for any $i= 1,2,3$. 
Then each of $a_i$
 ($i=1,2,3$) is an end of some dotted edge. Consider three other
ends.
If $S_0$ contains three ends of dotted edges,
then $\Sigma_{S_1}$ contains three dotted edges with mutually distinct 
ends (since $S_0\subset \Sigma_{S_1}$ and $a_i\in \Sigma_{S_1}$ for
$i=1,2,3$). 
This is impossible for a diagram of  a 4-polytope with $4+3$
facets.
Therefore, at least one of the ends of the dotted edges belongs to $\o
S_0$. Since $a_i\in \Sigma_{S_1}$ ($i=1,2,3$), 
$S_1$ contains no end of a dotted edge. Hence,
$y_4$ is the only end of a dotted edge in $\o S_0$. 
Then $S_2$ has at least 2 bad neighbors
(namely, $y_1$ and some $a_j$ such that 
$[a_j,y_4]=\infty$). 
If $S_2$ has exactly 2 bad neighbors, then $\Sigma_{S_2}$ is a
diagram of either an Esselmann polytope or of a 4-prism. 
However, it contains two dotted edges, which is
impossible. If $S_2$ has 3 bad neighbors, then  $\Sigma_{S_2}$ is a
diagram of a 4-simplex. At the same time, it contains at least one
dotted edge. 

Therefore, we may assume that the diagram  $\[a_1,S_0,\o S_0\]$
contains no dotted edges. In view of Lemma~\ref{multi-mult},
we are left with finitely many possibilities for the diagram
$\[a_1,S_0,\o S_0\]$.
The only two of these diagrams satisfying the signature condition and the
condition that $S_1$ has a unique bad neighbor are the diagrams shown
in Fig.~\ref{F4_1}.

For the diagram shown in  Fig.~\ref{F4_1}(a) consider the subdiagram 
$S_3=\[a_1,y_4\]$ having 2 bad neighbors. Then $\Sigma_{S_3}$ is a
diagram of a 5-prism, but $\Sigma_{S_1}$ contains a diagram 
$\[\t x_1,\t x_2,\t x_3, \t y_1,\t y_2\]$ of the type $B_3+B_2$ which is
impossible for a diagram of a 5-prism. 

Now, consider the diagram 
 shown in  Fig.~\ref{F4_1}(b). The subdiagram $S_4=\[x_4,a_1,y_3,y_2\]$
of the type $H_4$ has 3 bad neighbors $x_3, y_4$ and $y_1$.
Therefore, $a_2$ and $a_3$ are not neighbors of $S_4$.
In particular, $[a_i,y_3]=[a_i,y_2]=2$ for $i=2,3$.

Suppose that $[a_2,y_4]\ne 0$. Then the diagram $S_5=\[a_1,y_3,y_4\]$
of the type $H_3$ has 3 bad neighbors ($x_4,y_2$ and $a_2$).
Hence, $[a_3,y_4]=2$, and $[a_3,y_1]\ne 0$ (see Lemma~\ref{bad&Lanner}
applied to a Lann\'er diagram $L=\[y_1,y_2,y_3,y_4\]$). 
Furthermore, $P(S_5)$ is a simplex and $\Sigma_{S_5}$ is a Lann\'er
diagram of order 5.
By Cor.~\ref{dif2}, $\o S_5$ is also a Lann\'er diagram. 
Since $\o S_5$ contains the diagram $S_0$ of the type $F_4$,
$\o S_0$ is a cyclic Lann\'er diagram and $[a_3,x_1]=[a_3,x_4]=3$.
However, the diagram $\[x_1,a_3,x_4,a_1\]$ of the type $H_4$ has
at least 4 neighbors ($x_2,x_3,y_3$ and $y_1$,
 which is a neighbor of $a_3$).
The contradiction shows that $[a_2,y_4]=0$.
Similarly, $[a_3,y_4]=0$.

By Lemma~\ref{bad&Lanner}, $a_2$ and $a_3$ should
be joined with a Lann\'er diagram $\[y_1,y_2,y_3,y_4\]$, so,
both $a_2$ and $a_3$ are joined with $y_1$. Therefore, the diagram 
$\[a_1,y_3,y_2,y_1\]$ of the type $B_4$ has 4 bad neighbors
$a_2,a_3,x_4,y_4$,
which is impossible.

\begin{figure}[!h]
\begin{center}
\psfrag{a1}{$a_1$}
\psfrag{x1}{$x_1$}
\psfrag{x2}{$x_2$}
\psfrag{x3}{$x_3$}
\psfrag{x4}{$x_4$}
\psfrag{y1}{$y_1$}
\psfrag{y2}{$y_2$}
\psfrag{y3}{$y_3$}
\psfrag{y4}{$y_4$}
\psfrag{a}{(a)}
\psfrag{b}{(b)}
\epsfig{file=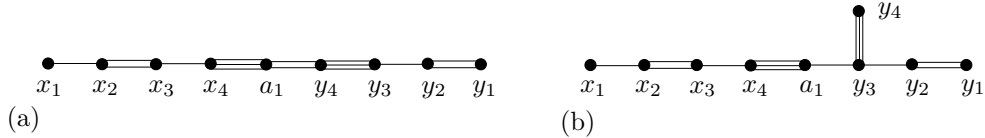,width=0.892\linewidth}
\caption{Case 1: two possibilities for  the diagram  
$\[a_1,S_0,\o S_0\]$}
\label{F4_1}
\end{center}
\end{figure}

\bigskip
\noindent
{\bf Case 2.} 
Suppose that $S_1$ has two bad neighbors, $y_4$ and $a_1$.
Then $P(S_1)$ is either an Esselmann polytope or a 4-prism.
Notice that  $\Sigma_{S_1}$ contains
a subdiagram $S_0$ of the type $F_4$ and contains no $G_2^{(k)}$ for $k>5$
(see Cor.~\ref{dif}).
There are two prisms and one Esselmann polytope
satisfying these conditions (see
Fig.~\ref{4-prism,f4}).
We consider these polytopes separately.

\begin{figure}[!h]
\begin{center}
\psfrag{a2}{$\t a_2$}
\psfrag{a3}{$\t a_3$}
\psfrag{x1}{$\t x_1$}
\psfrag{x2}{$\t x_2$}
\psfrag{x3}{$\t x_3$}
\psfrag{x4}{$\t x_4$}
\psfrag{a}{(a)}
\psfrag{b}{(b)}
\psfrag{c}{(c)}
\epsfig{file=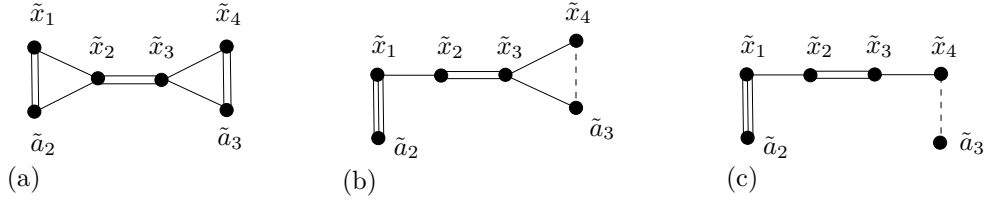,width=0.892\linewidth}
\caption{Diagrams of 4-polytopes with 6 facets containing $F_4$ 
and containing no $G_2^{(k)}$ for $k>5$.}
\label{4-prism,f4}
\end{center}
\end{figure}

\medskip
\noindent
{\bf Case 2.1.} 
Suppose that  $\Sigma_{S_1}$ is the diagram shown in
Fig.~\ref{4-prism,f4}(a). Then $\Sigma_{S_1}=\o S_1$. 
The node $a_1$ is a neighbor of $S_0=\[x_1,x_2,x_3,x_4\]$. 
Without loss of generality we may assume that $a_1$ is a neighbor of 
$\[x_1,x_2,x_3\]$. Then the subdiagram $\[x_1,x_2,x_3,a_3\]$ of the
type $F_4$ has 4 bad
neighbors ($x_4,a_1, a_2$ and some node of $\o S_0$ joined 
with $a_3$), which is impossible.

\medskip
\noindent
{\bf Case 2.2.} 
Suppose that  $\Sigma_{S_1}$ is the diagram shown in
Fig.~\ref{4-prism,f4}(b).  Then again $\Sigma_{S_1}=\o S_1$.
In particular, $a_2$ and $a_3$ are non-neighbors of 
$S_1$, and  $a_2$ and $a_3$ are joined with $y_4$.
Furthermore, 
the diagram $\[x_1,x_2,x_3,a_3\]$ of the type $F_4$ has 3 neighbors:
$x_4,a_2$ and $y_4$. Therefore, $a_1$ is not a neighbor of this
diagram,
hence $[a_1,x_4]\ne 2$.
Consider the diagram $S_6=\[a_2,x_1,x_2\]$ of the type $H_3$.
Clearly, it has no good neighbors, and $\o S_6=\Sigma_{S_6}$.
Furthermore,
$\o S_6$ contains a dotted edge $a_3x_4$, so, $S_6$ has at most 2
bad neighbors. Thus, it has exactly 2 bad neighbors, $x_3$ and
$y_4$.
Therefore, the subdiagram $\[a_3,x_4,a_1\]$ belongs to the subdiagram 
$\o S_6=\Sigma_{S_6}$, which is a diagram of a 4-prism.
Hence, the edge $\t x_4\t a_1$ is adjacent to a dotted edge in the diagram
of a 4-prism. This implies that $[x_4,a_1]=3$, and $\[x_1,x_2,x_3,x_4,a_1\]$
is a parabolic diagram of the type $\widetilde F_4$, which is impossible.

\medskip
\noindent
{\bf Case 2.3.} 
Suppose that  $\Sigma_{S_1}$ is the diagram shown in
Fig.~\ref{4-prism,f4}(c). 
We consider 2 cases: either $[a_3,x_4]=\infty$  or  $[a_3,x_4]\ne
\infty$ in $\Sigma$.

\smallskip
\noindent
{\bf Case 2.3.1.} 
Suppose that  $[a_3,x_4]=\infty$.
Since $a_2$ is not a good neighbor of $S_1$, $y_4$ is the
only node of $\o S_0$ joined with $a_2$.

Consider the diagram $S_7=\[a_2,x_1\]$ of the type $G_2^{(5)}$.
If $S_7$ has a bad neighbor, then $P(S_7)$ is a 5-polytope with
at most $5+3$ facets. At the same time, $\Sigma_{S_7}$ contains
a subdiagram
$\[\t x_2,\t x_3,\t x_4,\t a_3\]$, and $\t x_2\t x_3$ is a dotted edge in
$\Sigma_{S_7}$. However, no diagram of a 5-polytope with at most $5+3$ 
facets contains two dotted edges ($\t x_2\t x_3$ and $\t x_4\t a_3$)  
joined by a unique simple edge ($\t x_3\t x_4$).
Therefore, $S_7$ has no bad neighbors. 
In particular, $[a_2,y_4]=3$.

Consider the diagram $\[a_2,S_0,\o S_0\]$. It contains no
dotted edges and we completely know this diagram modulo 
finitely many possibilities for a Lann\'er
subdiagram $\o S_0$. 
However, for any of these possibilities 
the diagram  $\[a_2,S_0,\o S_0\]$ does not satisfy the
signature condition.
Hence, the case  $[a_3,x_4]=\infty$ is impossible.

\smallskip
\noindent
{\bf Case 2.3.2} 
Suppose that  $[a_3,x_4]\ne \infty$.
Then $a_3$ is a good neighbor of $S_1$ and
 $[a_3,x_4]=5$
(otherwise we have either $[a_3,x_4]=3$ and $\[a_3,S_0\]$
is a parabolic diagram of the type $\widetilde F_4$,
or $[a_3,x_4]=4$ and $\[a_3,x_4,x_3,x_2\]$
is a parabolic diagram of the type $\widetilde C_4$).

Consider the subdiagram $X=\[S_0,\o S_0,a_3\]$.
We know this diagram modulo finitely many possibilities for 
$\o S_0$  
($S_0$ is a diagram of the type $F_4$, $\o S_0$ is a Lann\'er diagram of
order 4, $S_1\subset \o S_0$ is a diagram of the type $H_3$ or $B_3$,
$a_3$ is a good neighbor of $S_1$,
and $a_3$ is joined with $S_0$ as shown in  Fig.~\ref{4-prism,f4}(c)).
The only diagram satisfying these conditions and the signature condition 
is the following:
%
\begin{center}
\psfrag{a3}{$a_3$}
\psfrag{x1}{$x_1$}
\psfrag{x2}{$x_2$}
\psfrag{x3}{$x_3$}
\psfrag{x4}{$x_4$}
\psfrag{y1}{$y_4$}
\psfrag{y2}{$y_3$}
\psfrag{y3}{$y_2$}
\psfrag{y4}{$y_1$}
\epsfig{file=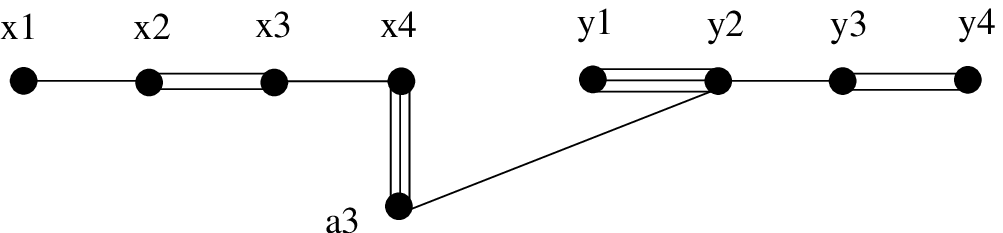,width=0.4592\linewidth}
\end{center}

Clearly, $a_2$ is not a good neighbor of $S_1$.
Hence, it is not a neighbor of $S_1=\[y_1,y_2,y_3\]$ and
$[a_2,y_4]\ne 2$ (Lemma~\ref{bad&Lanner}). 
Therefore, the diagram $S_2=\[y_2,y_3,y_4\]$ has 3 bad neighbors 
$a_2,a_3,y_1$. Thus, $\Sigma_{S_2}$ is a Lann\'er diagram of order 5
containing a subdiagram of the type $F_4$. This implies that 
$a_1$ is joined with $x_1$ and $x_4$. However, in this case the diagram
$\[x_4,a_3,y_3,y_2\]$ of the type $H_4$ 
has 4 bad neighbors $a_1,x_3,y_1,y_4$.

\medskip 

\noindent
We considered two prisms and an Esselmann polytope  
and no possibilities were found, so the proof of Case 2 is completed. 

\bigskip
\noindent
{\bf Case 3.} 
Suppose that $S_1$ has three bad neighbors, $y_4$, $a_1$
and $a_2$. Then $\o S_1$ is a cyclic Lann\'er diagram of order 5,
i.e. $\o S_1=\L_5^5$, and
$[a_3,x_1]=[a_3,x_4]=3$, $[a_3,x_2]=[a_3,x_3]=2$.
By Cor.~\ref{dif2}, $a_3$ is not a good neighbor of $S_1$, 
and $y_4$ is the only
node of $\o S_0$ joined with $a_3$.
 
As it is shown in Case~1 and Case~2, 
we may assume that the diagram $S_2=\[y_2,y_3,y_4\]$ has also 3 bad
neighbors. The reasoning above shows that each of $a_i$, $i=1,2,3$, 
is either a bad neighbor of $S_2$ or a non-neighbor, and exactly one
of $a_i$ is a non-neighbor.
Thus, $a_3$ and one of $a_1$ and $a_2$ (say $a_1$) are bad
neighbors of $S_2$
and $a_2$ is a non-neighbor.
Moreover, $\[a_1,S_0\]$ is a cyclic Lann\'er diagram,
and $y_1$ is the only node of $\o S_0$ joined with $a_1$.
Without loss of generality we may assume that $a_2$ is a neighbor of 
$\[x_1,x_2,x_3\]$. Then the diagram $\[a_1,x_1,x_2,x_3\]$ of the type
$B_4$ has 3 bad neighbors $x_1,a_2$ and $a_3$.
Hence, $y_1$ is a good neighbor of $\[a_1,x_1,x_2,x_3\]$.
So, $[y_1,a_1]=3$. Recall, that $a_1$ is a bad neighbor of $S_1$.
This implies that either $[y_1,y_2]=4$, $5$ or
$[y_1,y_3]=4$, 5.
Therefore, either $y_2$ or $y_3$ is a bad neighbor of 
$\[y_1,a_1,x_1,x_2,x_3\]$, which is impossible. 

\end{proof}

\begin{lemma}
\label{h4,3sos}
If $\Sigma$ contains a subdiagram of the type $H_4$
then $\Sigma$ is a diagram $\Sigma_{P_7}$.

\end{lemma}

\begin{proof}
Suppose that $S_0=\[x_1,x_2,x_3,x_4\]$ 
is a subdiagram of $\Sigma$ of the type $H_4$.
By Lemmas~\ref{H_4,1sos} and~\ref{h4,2sos}, $S_0$ has 3 neighbors
$a_1,a_2$ and $a_3$. Then $P(S_0)$ is a 3-simplex, and 
$\o S_0=\Sigma_{S_0}=\[y_1,y_2,y_3,y_4\]$ is a Lann\'er diagram of order 4.
Let $S_1=\[y_1,y_2,y_3\]$ be a subdiagram of $\o S_0$ 
of the type $H_3 $ or $B_3$.
$S_1$ has at least one bad neighbor, $y_4$.
We consider three cases in which $S_1$ has 1, 2 and 3 bad neighbors
respectively. 

\bigskip
\noindent
{\bf Case 1.} 
Suppose that $S_1$ has a unique bad neighbor $y_4$.
Then $P(S_1)$ is a 4-polytope with $4+3$ facets.
The same reasoning as in Case 1 of Lemma~\ref{f4,3sos} shows that
for some $i\in \{1,2,3\}$ the diagram  $\[a_i,S_0,\o S_0\]$
contains no dotted edges. 
We
assume that the diagram  $\[a_1,S_0,\o S_0\]$
contains no dotted edges. In view of Lemma~\ref{multi-mult}, 
we are left with finitely many possibilities for the diagram
 $\[a_1,S_0,\o S_0\]$.
Recall that $S_1\subset \o S_0$ is a subdiagram of the type $H_3$ or
$B_3$ with a unique bad neighbor. There are only three possibilities for
 the diagram $\[a_1,S_0,\o S_0\]$ satisfying this condition
together with the signature condition, namely, 
the diagrams shown in Fig.~\ref{2shemi}.
We consider these diagrams separately.

\begin{figure}[!h]
\begin{center}
\psfrag{a1}{$a_1$}
\psfrag{a2}{$a_2$}
\psfrag{a3}{$a_3$}
\psfrag{x1}{$x_1$}
\psfrag{x2}{$x_2$}
\psfrag{x3}{$x_3$}
\psfrag{x4}{$x_4$}
\psfrag{y1}{$y_1$}
\psfrag{y2}{$y_2$}
\psfrag{y3}{$y_3$}
\psfrag{y4}{$y_4$}
\psfrag{a}{(a)}
\psfrag{c}{(b)}
\psfrag{b}{(c)}
\epsfig{file=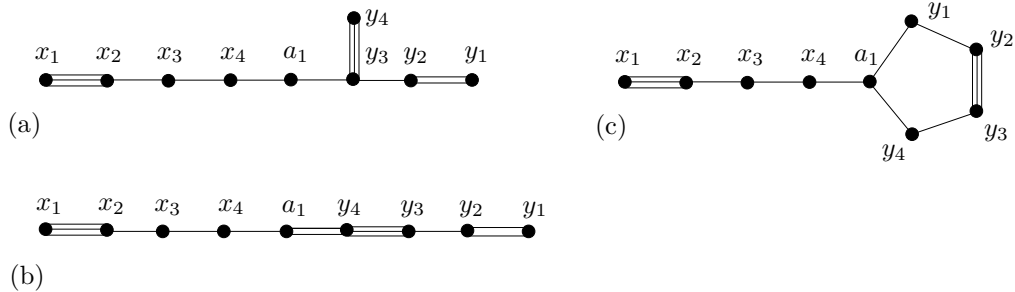,width=0.92\linewidth}
\caption{Three possibilities for
 the diagram $\[a_1,S_0,\o S_0\]$.}
\label{2shemi}
\end{center}
\end{figure}

Consider the diagrams shown in  Fig.~\ref{2shemi}(a) and~\ref{2shemi}(b).
Clearly, $S_1=\[y_2,y_3,y_4\]$. 
Suppose that $a_2$ is joined with $y_3$ or $y_4$.
Since $a_2$ is not a bad neighbor of $S_1$, we obtain a subdiagram
$\[a_2,S_1\]$ of the type $F_4$, which contradicts 
Lemma~\ref{f4,3sos}.
Hence, both $a_2$ and $a_3$ are joined with either $y_1$ or $y_2$ 
(Lemma~\ref{bad&Lanner}),
and the subdiagram $\[y_1,y_2,y_3\]$ of the type $H_3$ has 4 bad neighbors 
($y_4,a_1,a_2,a_3$), which is impossible.

Consider the diagram shown in  Fig.~\ref{2shemi}(c).
Without loss of generality we assume that $S_1=\[y_1,y_2,y_3\]$.
Each of the diagrams $S_2=\[a_1,y_1,y_2,y_3\]$ and 
$S_3=\[a_1,y_4,y_3,y_2\]$ of the type $H_4$ has 
bad neighbors $x_4$ and $y_1$ (or $y_4$ respectively).
Thus, each of these diagrams has at most one extra neighbor
(Lemma~\ref{h4,2sos}). On the other hand, by Lemma~\ref{bad&Lanner},
each of the nodes $a_2$ and $a_3$ is joined with the diagram
$\[a_1,y_1,y_2,y_3,y_4\]$.  
Thus, we may assume that 
$a_2$ is not joined with $S_2$ and joined with  $y_4$, and
$a_3$ is not joined with $S_3$ and joined with  $y_1$.
Since $S_1$ has no bad neighbor besides $y_4$, $[a_4,y_1]=3$.
Furthermore, $[a_3,a_1]=2$ (otherwise the diagram $\[S_1,a_3\]$
of the type $H_4$ has 4 bad neighbors $y_4,a_2,a_1$ and some $x_i$,
$i\in \{1,2,3,4\}$).
If $[a_1,x_4]=2$, then $\[x_4,a_1,y_4,y_1,y_2,a_3\]$ is a parabolic
diagram of the type $\widetilde D_5$, hence,  $[a_3,x_4]\ne 2$.
Moreover,  $[a_1,x_4]\ne 3$ (otherwise, $\[x_4,a_1,y_1,a_3\]$ is a
parabolic diagram $\widetilde A_3$). Therefore, the diagram 
$S_4=\[x_2,x_3,x_4,a_1,y_4,y_3\]$ of the type $A_6$ has three bad
neighbors $x_1,y_2$ and $a_3$. Hence, $a_2$ is a good neighbor of
$S_4$, and $[a_2,y_4]=3$, $[a_2,x_4]=[a_2,x_3]=[a_2,x_2]=2$.
Therefore, $\[x_3,x_4,a_1,y_4,a_2,y_1,y_2\]$ is a parabolic diagram of the
type $\widetilde E_6$, which is impossible.

\bigskip
\noindent
{\bf Case 2.} 
Suppose that $S_1$ has two bad neighbors, $y_4$ and $a_1$.
Then $P(S_1)$ is either an Esselmann polytope or a 4-prism, 
$\Sigma_{S_1}=\[\t x_1, \t x_2, \t x_3, \t x_4,\t a_2,\t a_3\]$.

A Coxeter diagram of an Esselmann polytope containing no subdiagram of
the types $F_4$ and $G_2^{(k)}$, $k\ge 6$, is one of the two diagrams
shown in  Fig.~\ref{ess}. 
A Coxeter diagram of a 4-prism containing no subdiagram of the type
$F_4$ is one of the diagrams shown in Fig.~\ref{4-prisms}.
The diagram (a) is drown for
3 times, since there are 3 different possibilities for the location of
$S_0$  in this diagram
(in any other diagram containing two subdiagrams of the type $H_4$
these subdiagrams are permuted by an automorphism of the diagram).

\begin{figure}[!h]
\begin{center}
\psfrag{a2}{$\t a_2$}
\psfrag{a3}{$\t a_3$}
\psfrag{x1}{$\t x_1$}
\psfrag{x2}{$\t x_2$}
\psfrag{x3}{$\t x_3$}
\psfrag{x4}{$\t x_4$}
\psfrag{b}{(b)}
\psfrag{a}{(a)}
\epsfig{file=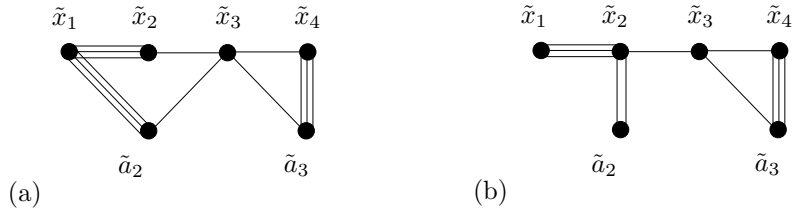,width=0.72\linewidth}
\caption{Diagrams of Esselmann polytopes containing 
$H_4$ and containing no $G_2^{(k)}$ for $k>5$.}
\label{ess}
\end{center}
\end{figure}

We assume that $\[\t x_1,\t x_2,\t x_3\]$ is
the subdiagram of $\Sigma_{S_1}$ 
of the type $H_3$. It turns out that in any possible case for the
diagram $\o S_1$ the subdiagram  $\[\t x_1,\t x_2,\t x_3\]$ has a unique bad
neighbor in $\Sigma_{S_1}$.
We assume that $\t a_2$ is the bad neighbor of the diagram 
$\[\t x_1,\t x_2,\t x_3\]$,
then $\t a_3$ is the remaining node of $\Sigma_{S_1}$.

First we consider the Esselmann polytopes and then prisms.

\medskip
\noindent
{\bf Case 2.1.} 
Suppose that  $\Sigma_{S_1}$ is an Esselmann diagram.
The subdiagram $S_5=\[x_1,x_2,x_3,a_3\]$ of the type $H_4$ 
has 3 bad neighbors in $\Sigma$ ($a_2,x_4$
and some node of $\o S_0$ joined with $a_3$).
Hence, $a_1$ is not a neighbor of $S_5=\[x_1,x_2,x_3,a_3\]$,
so $a_1$ is joined with $x_4$ (as a neighbor of $S_0$).
Furthermore, the subdiagram $\[S_0,a_3,\o S_0\]$ 
consists of two diagrams  $\[S_0,a_3\]$ and $\o S_0$ 
joined by the edge $a_3y_4$ only (and this edge does exist
by Lemma~\ref{bad&Lanner}).
Moreover, this edge is a dotted one, otherwise the diagram 
 $\[S_0,a_3,\o S_0\]$  is superhyperbolic.
In particular, three bad neighbors of the diagram
$S_5=\[x_1,x_2,x_3,a_3\]$ are $a_2,x_4,y_4$, so
  $\[S_1,a_1\]=\o S_5 = \Sigma_{S_5}$ is a Lann\'er diagram of order 4. 

Consider the diagram  $S_6=\[a_3,x_4\]$  of the type $G_2^{(5)}$. 
$S_6$ has at least two bad neighbors, 
$y_4$ and $x_4$, and $P(S_6)$ is a 5-dimensional polytope with at most
$5+2$ facets. This implies that   $P(S_6)$ is a 5-prism.
Notice that the subdiagram  $X=\[x_1,x_2,a_2,y_1,y_2,y_3\]$  is not
joined with $S_6$ in $\Sigma$, so it does not differ from
the subdiagram $\[\t x_1,\t x_2,\t a_2,\t y_1,\t y_2,\t y_3\]$  
of $\Sigma_{S_6}$. It is clear, that $X$ does not contain dotted edges
($a_2$ cannot be joined with $S_1$ by a dotted edge, since
$a_2$ is not a bad neighbor of $S_1$).
Moreover, as a diagram of a 5-prism should not contain 
Lann\'er diagrams of order 3, we obtain that $a_2$ is a good neighbor of $S_1$,
the diagram $\Sigma_{S_1}$ is the one shown in  Fig.~\ref{ess}(b),
$[a_2,x_2]=3$, and $S_1$ is a diagram of the type $B_3$
(if it is of the type $H_3$, then $\t a_2\t x_2$ is a dotted edge
in $\Sigma_{S_6}$, so it
cannot be adjacent to a triple edge  $\t x_1 \t x_2$ 
in a diagram of a 5-prism).
Without loss of generality we may assume that $[y_1,y_2]=3$,
$[y_2,y_3]=4$.
Then we obtain a linear diagram $\[x_1,x_2,a_2,y_1,y_2,y_3\]$
with the following labels on its edges: 5,3,3,3,4.
This implies that $[a_1,y_3]=\infty $ in $\Sigma_{S_6}$,
thus, $[a_1,y_3]\ne 2,3$ in $\Sigma$ (if  $[a_1,y_3]=3$,
the edge $\t a_1\t y_3$ of $\Sigma_{S_6}$ 
is labeled by 10,  and it is not a 
dotted edge as it should be).
However, the multiple
edges $a_1y_3$ and $y_3y_2$ cannot be adjacent
 in Lann\'er diagram  $\[S_1,a_1\]$ of order 4.

The contradiction shows that $\Sigma_{S_1}$ is not a diagram 
of an Esselmann polytope.

\medskip
\noindent
{\bf Case 2.2.} 
Suppose that $\Sigma_{S_1}$ is a 4-prism.
Notice, that in all diagrams of 4-prisms all edges incident to $a_2$
are simple, and by Cor.~\ref{dif2},  
$a_2$ cannot be a good neighbor of $S_1$.
Hence, $a_2$ is a neighbor of $y_4$ (Lemma~\ref{bad&Lanner}).

Furthermore, suppose that $[a_2,y_4]\ne \infty$. Then the diagram
$\[S_0,a_2,\o S_0\]$ 
satisfies the following properties:
it contains no dotted edges;  
$y_4$ is the only node of $\o S_0$ joined with  $a_2$;
  $a_2$ is a
bad neighbor of $H_3\subset S_0$ and $a_2$ is joined with $S_0$ by
simple edges only.
It is easy to check that no diagram satisfies these properties
together with signature condition.
We obtain that  $[a_2,y_4]= \infty$. 

Consider the diagrams of 4-prisms case-by-case.

\begin{figure}[!h]
\begin{center}
\psfrag{a2}{$\t a_2$}
\psfrag{a3}{$\t a_3$}
\psfrag{x1}{$\t x_1$}
\psfrag{x2}{$\t x_2$}
\psfrag{x3}{$\t x_3$}
\psfrag{x4}{$\t x_4$}
\psfrag{a_1}{$(a1)$}
\psfrag{a_2}{$(a2)$}
\psfrag{a_3}{$(a3)$}
\psfrag{b}{(b)}
\psfrag{c}{(c)}
\psfrag{d}{(d)}
\psfrag{e}{(e)}
\psfrag{f}{(f)}
\psfrag{g}{(g)}
\psfrag{2,3}{{\scriptsize $2,3$}}
\epsfig{file=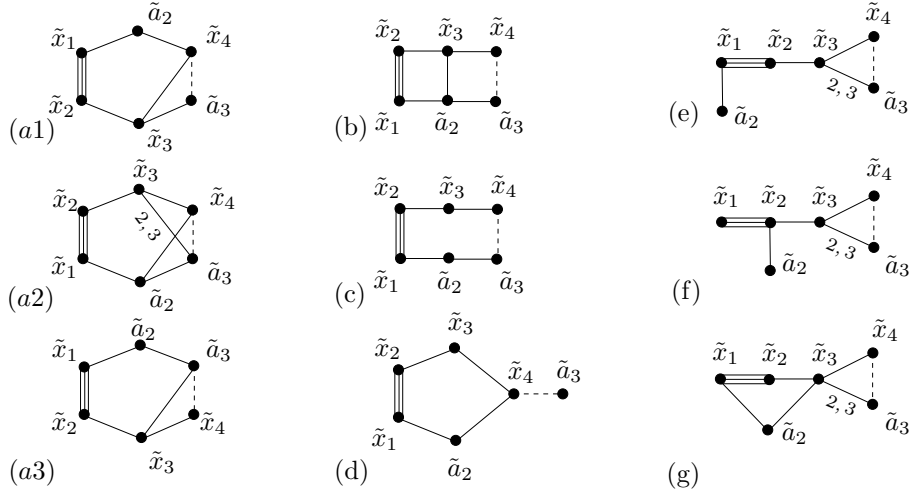,width=0.82\linewidth}
\caption{Diagrams of 4-prisms containing $H_4$.}
\label{4-prisms}
\end{center}
\end{figure}

\smallskip
\noindent
{\bf Case 2.2.1.} 
Suppose that  $\Sigma_{S_1}$ is the diagram shown in
Fig.~\ref{4-prisms}(a1).
If $a_1$ is a neighbor of $\[x_1,x_2,x_4\]$ 
then $\[x_2,x_1,a_2,x_4\]$ is a diagram of the type $H_4$ with 4
neighbors $x_3,a_3,a_1,y_4$ (recall that $[a_2,y_4]=\infty$). If  
$a_1$ is not a neighbor of $\[x_1,x_2,x_4\]$, then $[a_1,x_3]\ne 2$ 
and  $\[x_1,x_2,x_3,a_3\]$ is a diagram of the type $H_4$ with 4
neighbors ($x_3,a_2,a_1$ and some $y_i$, $i\in \{1,2,3,4\}$, 
see Lemma~\ref{bad&Lanner}).  

\smallskip
\noindent
{\bf Case 2.2.2.} 
Suppose that  $\Sigma_{S_1}$ is one of the diagrams shown in
Fig.~\ref{4-prisms}(a2)-(c).
Since the diagram $\[x_2,x_1,a_2,a_3\]$ of the type $H_4$
has 3 bad neighbors $x_3, x_4$ and $y_4$ (since
$[a_2,y_4]=\infty$), the diagram $\[S_1,a_1\]$ is
a Lann\'er diagram. Clearly, the nodes $x_1,x_2, a_2$ are not attached
to this Lann\'er diagram
(otherwise $\[x_2,x_1,a_2,a_3\] $ has 4 bad neighbors, $x_3,x_4,y_4$
and some node of $\[S_1,a_1\]$).
Since
a Lann\'er diagram $\[a_1,S_1\]$ is joined with a Lann\'er
diagram $\[x_1,x_2,x_3,a_2\]$, we obtain that $[x_3,a_1]\ne 2$.
The signature condition applied to the  
diagram $\[S_0,a_1,a_2,S_1\]$ implies that
$[x_3,a_1]=\infty$.
However, this means that the diagram
$\[x_2,x_3,x_4,a_2\]$ of the type $A_4$ (in case (a2)), or $A_3+A_1$ 
(in cases (a3) and (c)), or $D_4$ (in case (b)) has 4 bad neighbors,
$x_1,a_3,y_4,a_1$ (again, we use that $[a_2,y_4]=\infty$).

\smallskip
\noindent
{\bf Case 2.2.3.} 
Suppose that  $\Sigma_{S_1}$ is the diagram shown in
Fig.~\ref{4-prisms}(d).
Since the diagram $\[x2,x_1,a_2,x_4\]$ of the type $H_4$
has 3 bad neighbors, $x_3,a_3$ and
$y_4$,
the node $a_1$ is not a neighbor of  $\[x2,x_1,a_2,x_4\]$ and 
$\[a_1,S_1\]$ is a Lann\'er diagram. Hence,
$[a_1,x_3]\ne 2$  ($a_1$ is a neighbor of $S_0=\[x_1,x_2,x_3,x_4\]$
and non-neighbor of $\[x_1,x_2,x_4\]$).
The signature condition applied to the diagram
 $\[S_0,a_1,a_2,S_1\]$ implies that
$[x_3,a_1]=\infty$.

Suppose that $a_3$ is not a good neighbor of $S_1$.
Then the dotted edge $\t x_4\t a_3\subset \Sigma_{S_1}$ remains a
dotted edge  $x_4a_3\subset \o S_1$ in $\Sigma$,   
and the diagram $\[x_2,x_3,x_4,a_2\]$ of the type $A_4$ has 4 bad
neighbors, $x_1,a_1,a_3,y_4$ (we use that $[x_3,a_1]=\infty$ and 
$[a_2,y_4]=\infty$), which is impossible.

Suppose that $a_3$ is a good neighbor of $S_1$.
If the edge  $x_4a_3\subset \o S_1$ is dotted or multiple one,
then $a_3$ remains a bad neighbor of
 $\[x_2,x_3,x_4,a_2\]$  and this diagram still has 4 bad neighbors.
Therefore, we may assume that   $x_4a_3$ is a simple
edge. 
Hence, $S_1$ is a diagram of the type $H_3$ (otherwise 
$\[a_2,x_3,x_4,a_3,S_1\]$ is a parabolic diagram of the type
$\widetilde B_6$).  So, we have the following  diagram:

\begin{center}
\psfrag{a1}{$a_1$}
\psfrag{a2}{$a_2$}
\psfrag{a3}{$a_3$}
\psfrag{x1}{$x_1$}
\psfrag{x2}{$x_2$}
\psfrag{x3}{$x_3$}
\psfrag{x4}{$x_4$}
\psfrag{y1}{$y_1$}
\psfrag{y2}{$y_2$}
\psfrag{y3}{$y_3$}
\psfrag{y4}{$y_4$}
\epsfig{file=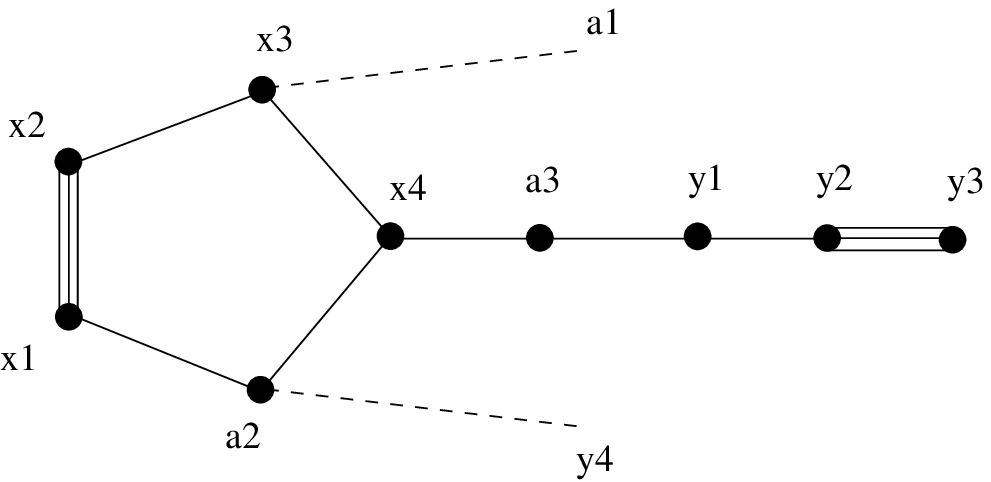,width=0.35\linewidth}
\end{center}

Furthermore, $[y_4,a_3]\ne \infty$ by Lemma~\ref{2dotted}.
If $[y_4,a_3]=5$ then  $\[y_4,a_3,x_4,x_3\]$ 
  is a diagram of the type $H_4$ with 4 bad neighbors, $a_1,a_2,y_1,x_2$.
If $[y_4,a_3]=4$ then  $\[y_4,a_3,x_4,x_3,x_2\]$ 
  is a diagram of the type $B_5$ with 4 bad neighbors, $a_1,a_2,y_1,x_1$.
Hence,  $[y_4,a_3]=2$ or 3.
Similarly,  $[a_1,a_3]=2$ or 3.

Consider the diagram $\[S_0,a_3,S_1\]$.
The node $y_4$ is not joined with $S_0$.
Since the diagram $\[x_2,x_3,x_4,a_3,y_1,y_2\]$
has 3 bad neighbors ($x_1,a_1,y_3$), $y_4$ may be joined
with $a_3$ by a simple edge only. Recall also 
that $\[S_1,y_4\]=\o S_0$ is a Lann\'er diagram,
so by Lemma~\ref{bad&Lanner} we obtain that $[y_4,a_3]\ne 2$. 
Therefore, $[y_4,a_3]=3$.
It turns out that
for any Lann\'er diagram $\[y_4,S_1\]$ (where $S_1$ is of the
type $H_3$) the diagram 
$\[S_0,a_3,y_4,S_1\]$ does not satisfy the signature condition
unless $[y_4,y_3]=3$ and 
$[y_4,y_1]=[y_4,y_2]=2$.
By the similar consideration of the diagram
$\[x_2,x_1,a_2,x_4,a_3,a_1,S_1\]$ we obtain 
$[a_1,a_3]=[a_1,y_3]=3$ and 
$[a_1,y_1]=[a_1,y_2]=2$.
However, in this case the diagram $\[a_3,y_4,y_3,y_2\]$ of the type $H_4$ 
has 4 bad neighbors ($a_2,x_4,y_1,a_1$).

\smallskip
\noindent
{\bf Case 2.2.4.} 
Suppose that  $\Sigma_{S_1}$ is one of the diagrams shown in
Fig.~\ref{4-prisms}(e)-(g).
In these cases $a_2$ still cannot be a good neighbor of $S_1$,
therefore, it is not a neighbor and $[a_2,y_4]=\infty$ as it was shown above.
At the same time, if $[x_3,a_3]=2$, 
then $a_3$ may be a good neighbor of $S_1$
and it may occur that $[x_4,a_3]\ne \infty$ in $\Sigma$.

Suppose that  $[x_4,a_3]\ne \infty$ in $\Sigma$, and hence,  $[x_3,a_3]=2$.
Consider the subdiagram $Y=\[S_0,a_2,a_3,S_1\]$.
If $[x_4,a_3]=3$ or 5, then $Y$ does not satisfy the signature
condition (namely, at least one of the subdiagrams $Y\setminus a_2$ and
$Y\setminus y_4$ does not satisfy the condition). 
If $[x_4,a_3]=4$  then $Y$ contains a parabolic subdiagram 
$\[x_2,x_3,x_4,a_3,y_1\]$ of the type $\w F_4$
(here we assume that the nodes of $S_1$ are numbered in such a way
that $[y_1,y_2]=3$, $[y_2,y_3]=4$ or 5).

Therefore,  $[x_4,a_3]=\infty$.
Consider the subdiagram $S_7=\[x_2,x_1,a_2\]$ of the type $H_3$.
It has at least 2 bad neighbors $x_3$ and $y_4$. 
Furthermore, the dotted edge
$x_4a_3$ is not joined with $S_7$, so $\[x_4,a_3\]\subset \o S_7$,
which implies that $S_7$ has no extra bad neighbors. 
In particular, $a_1$ is not a bad neighbor of $S_7$. 
Suppose that $a_1$ is a good neighbor of $S_7$. Then 
$S_8=\[a_1,S_7\]$ is a diagram of the type $H_4$ 
with 3 bad neighbors ($x_3,y_4$
and some node of $S_1$ joined with $a_1$). 
However, $\o S_8$ contains a dotted edge 
$x_4a_3$. The contradiction shows that $a_1$ is not a neighbor
of $S_7$. Thus, the diagram 
$\o S_7 =\Sigma_{S_7}=\[x_4,a_3,a_1,S_1\]$ 
is a diagram of a 4-prism, which implies that $\[a_1,S_1\]$ is a
Lann\'er diagram.
A Lann\'er diagram $\[a_2,x_1,x_2,x_3\]$ should be joined with
a Lann\'er diagram $\[a_1,S_1\]$.
Since $a_1$ is not a neighbor of $S_7$, this implies that 
$[a_1,x_3]\ne 2$. 
Consider the subdiagram $S_9=\[a_2,x_2,x_3,x_4\]$ (of the type $A_1+A_3$ or
$A_4$ or $D_4$ in cases (e), (f) and (g) respectively). 
It has 3 bad neighbors $x_1,a_3,y_4$, so
$a_1$ is a good neighbor of $S_9$, 
which implies $[a_1,x_3]=3$, $[a_1,x_4]=2$.

Finally, consider the diagram $Z=\[S_0,a_2,a_1,S_1\]$.
This diagram contains no dotted edges, so we have finitely many
possibilities for $Z$. Moreover, $Z$ satisfies the following
conditions:
$\[S_0,S_1\]$ is a diagram of the type $H_4+H_3$ or $H_4+B_3$;
$a_2$ is not joined with $S_1$ and is joined with $S_0$ in one of the
ways shown in Fig.~\ref{4-prisms}(e)-(g);
$[a_1,x_3]=3$ and $a_1x_3$ is a unique edge connecting $a_1$ with 
$\[a_2,S_0\]$; $\[a_1,S_1\]$ is a Lann\'er diagram.
However, no diagram satisfies these conditions together with the
signature condition.

\bigskip
\noindent
{\bf Case 3.} 
Suppose that $S_1$ has 3 bad neighbors.
Let $S_1'\subset \o S_0$ be a subdiagram of the type $H_3$ or $B_3$,
$S_1'\ne S_1$
 (see Lemma~\ref{2subd_of_L4}). 
We denote $\[y_1,y_2,y_3\]=S_1$, $\[y_2,y_3,y_4\]=S_1'$.
By Cases 1 and 2 we may assume that 
both $S_1$ and $S_1'$ have 3 bad neighbors.
We have two possibilities (up to permutation of the nodes $a_1, a_2,
a_3$): either $a_1$ and $a_2$ are bad neighbors of $S_1$ as well as
of $S_1'$ (in addition to $y_4$ and $y_1$ respectively), or 
$y_4$, $a_1$, $a_2$ are bad neighbors of $S_1$
and   $y_1$, $a_2$, $a_3$ are bad neighbors of $S_1'$.

Suppose that $a_1$ and $a_2$ are bad neighbors for each of $S_1$ and 
$S_1'$. Then the node $a_3$ is not a
bad neighbor for both $S_1$ and $S_1'$.
The diagram $\Sigma_{S_1}$ is
a Lann\'er diagram of order 5. By Cor.~\ref{dif}, $\o S_1=\[S_0,a_3\]$
is a Lann\'er diagram, too.
So, the diagram $\[S_0,a_3,\o S_0\]$
consists of a Lann\'er diagram $\[S_0,a_3\]$ (where $S_0$ is of the
type $H_4$), and a Lann\'er diagram $\o S_0$, the node $a_3$ is not a
bad neighbor for both $S_1$ and $S_1'$. The only such a diagram
satisfying the signature condition is the diagram shown in 
 Fig.~\ref{2shemi}(b).
However, by Proposition~\ref{al}, the edge $\t x_4\t a_3$ of
 $\Sigma_{S_1}$ is a dotted edge, which contradicts  the fact that
 $\Sigma_{S_1}$ is a Lann\'er diagram of order 5.

Therefore, we may assume that 
$y_4$, $a_1$, $a_2$ are bad neighbors of $S_1$,
and   $y_1$, $a_2$, $a_3$ are the bad neighbors of $S_1'$.
By consideration of $\o S_1$ and $\o S_1'$ we conclude that
$\[S_0,a_1\]$ and $\[S_0,a_3\]$ are Lann\'er diagrams.

Now, we consider two cases: either both $a_1$ and $a_3$ are joined with
$\o S_0$ by dotted edges, or at least one of $a_1$ and $a_3$ 
(say $a_1$)
is joined with $\o S_0$ by non-dotted edges only.

\medskip
\noindent
{\bf Case 3.1.} 
Suppose that  both $a_1$ and $a_3$ are joined with
$\o S_0$ by dotted edges. Since $a_1$ is a bad neighbor of $S_1$
and it is not a bad neighbor of $S_1'$, the dotted edge joining $a_1$ with 
$\o S_0$ is $a_1y_1$. Similarly, $[a_3,y_4]=\infty$.
By Lemma~\ref{2dotted},  $[a_1,a_3]\ne \infty$.
 Hence, we have finitely many possibilities for the diagram
$\[S_0,a_1,a_3\]$ (4 possibilities for each of Lann\'er diagrams
$\[S_0,a_1\]$ and $\[S_0,a_3\]$ and 3 possibilities for 
$[a_1,a_3]\in \{3,4,5\}$).

Consider the subdiagram $S_{10}=\[x_2,x_3,x_4,a_1\]$ of the type
$D_4,A_4,B_4$ or $H_4$. If $a_3$ is a bad neighbor of $S_{10}$ then
$S_{10}$ has 3 bad neighbors $x_1,y_1$ and $a_3$, thus, $a_2$ is not a
bad neighbor of $S_{10}$. 
If $a_3$ is not a bad neighbor of $S_{10}$, then the diagram
$\[S_{10},a_3\]$ has 3 bad neighbors $x_1,y_1,y_4$, and $a_2$ is not a
bad neighbor of $\[S_{10},a_3\]$.       
In any case, $a_2$ is not a bad neighbor of $S_{10}$.
Similarly, $a_2$ is not a bad neighbor of $S_{11}=\[x_2,x_3,x_4,a_3\]$.

Suppose that $a_2$ is a bad neighbor of the diagram 
$S_{12}=\[x_1,x_2,x_3\]$ of the type $H_3$. Then $P(S_{12})$ is a 
4-polytope with at most $4+3$ facets. 
However, $\Sigma_{S_{12}}\ne \o S_{12}$, and
$\Sigma_{S_{12}}$ contains 4 dotted edges
$\t y_4\t a_3,\t a_3\t x_4,\t x_4\t a_1,\t a_1\t y_1$,
which is impossible for a diagram of a 4-polytope with at most 7 facets.
Therefore, $a_2$ is not a bad neighbor of $S_{12}$, which implies that
$a_2$ is not a neighbor of $x_1x_2$. Hence, $a_2$ is a neighbor
of $x_3x_4$.
 Thus, each of $S_3$ and $S_{11}$ is a 
diagram of the type $D_4$ or $A_4$.

Suppose that $S_{10}$ and $S_{11}$ are of the same type.
Then $[a_1,a_3]=2$ (otherwise either $[a_1,a_3]=\infty$ in
contradiction to Lemma~\ref{2dotted}, 
or $\[a_1,a_3\]$ has 3 bad neighbors in 
contradiction to Lemma~\ref{g2,3badsos},
or one of the diagrams $\[x_3,a_1,a_3\]$ and $\[x_4,a_1,a_3\]$ is a
parabolic diagram of the type $\w A_2$). 
This implies that $\[S_{10},a_3,a_2\]$
contains a parabolic subdiagram.
Hence, $S_{10}$ and $S_{11}$ are of different types.
If $[a_1,a_3]=2$ or 3 then  $\[S_{10},a_3,a_2\]$
still contains a parabolic diagram.
If $[a_1,a_3]=4$ or 5 then the diagram $S_{13}=\[a_1,a_3\]$ has 2 bad
neighbors, so $\Sigma_{S_{13}}$ is a diagram of a 5-prism. However, 
in $\Sigma_{S_{13}}$
the subdiagram $\[\t x_1,\t x_2,\t x_3,\t x_4\]$ is a linear diagram
in which one triple edge and one simple edge are joined by the dotted
edge $\t x_2\t x_3$. No diagram of a 5-prism contains such
a subdiagram, which shows that the case 3.1 is impossible.

\medskip
\noindent
{\bf Case 3.2.} 
Suppose that $a_1$ is joined with $S_0$ by non-dotted edges only.
Then the subdiagram $T=\[S_0,a_1,\o S_0\]$ contains no dotted
edges, so we have finitely many possibilities for the diagram $T$.
Notice also that $T$ satisfies the following properties:\\
\phantom{\qquad}
$S_0$ is a diagram of the type $H_4$, and $\[S_0,a_1\]$ 
is a Lann\'er diagram;\\
\phantom{\qquad}
$\o S_0$ is a Lann\'er diagram of order 4, $\o S_0$ is not joined with $S_0$;\\
\phantom{\qquad}
$S_1$ and $S_1'$ are subdiagrams of $\o S_0$ or the types $H_3$ or
$B_3$;\\
\phantom{\qquad}
$a_1$ is a bad neighbor of $S_1$ and not a bad neighbor of $S_1'$.\\
However, there are only two diagrams satisfying these properties
together with signature condition, namely the diagrams shown in
Fig.~\ref{9points}(a) and~\ref{9points}(b).

Consider the diagram shown in  Fig.~\ref{9points}(a).
The subdiagram $\[x_2,x_3,x_4,a_1,y_1\]$ of the type $B_5$ has 3 bad
neighbors in $\Sigma$: $x_1,y_2$ and $a_3$. 
Hence, $a_2$ is not a bad neighbor of
 $\[x_2,x_3,x_4,a_1,y_4\]$, which implies that $a_2$ is joined with 
$\[x_1,x_2\]$. Therefore, $a_2$ is a bad neighbor of
$S_{14}=\[x_1,x_2,x_3\]$, and $P(S_{14})$ is a 4-polytope with $4+3$ facets
(it is easy to see that $S_{14}$ has no other bad neighbors). 
However, $\Sigma_{S_{14}}$ contains $\[\t a_1,\t y_1,\t y_2\]$, 
the subdiagram
composed of a double edge $\t a_1\t y_1$ adjacent to a 
triple edge $\t y_1\t y_2$,
which is impossible for a diagram of a 4-polytope with $4+3$ facets. 

In particular, we obtain that for any diagram $S\subset \Sigma$ of the
type $H_4$ the diagram $Q=\[S,\o S,a_1\]$ is the diagram shown
in  Fig.~\ref{9points}(b) and
$\o S$ is a linear Lann\'er diagram with one
double, one simple and one triple edge.

\begin{figure}[!h]
\begin{center}
\psfrag{a1}{$a_1$}
\psfrag{a3}{$a_3$}
\psfrag{y1}{$y_1$}
\psfrag{y2}{$y_2$}
\psfrag{y3}{$y_3$}
\psfrag{y4}{$y_4$}
\psfrag{x1}{$x_1$}
\psfrag{x2}{$x_2$}
\psfrag{x3}{$x_3$}
\psfrag{x4}{$x_4$}
\psfrag{a}{  (a)}
\psfrag{b}{(b)}
\psfrag{c}{(c)}
\epsfig{file=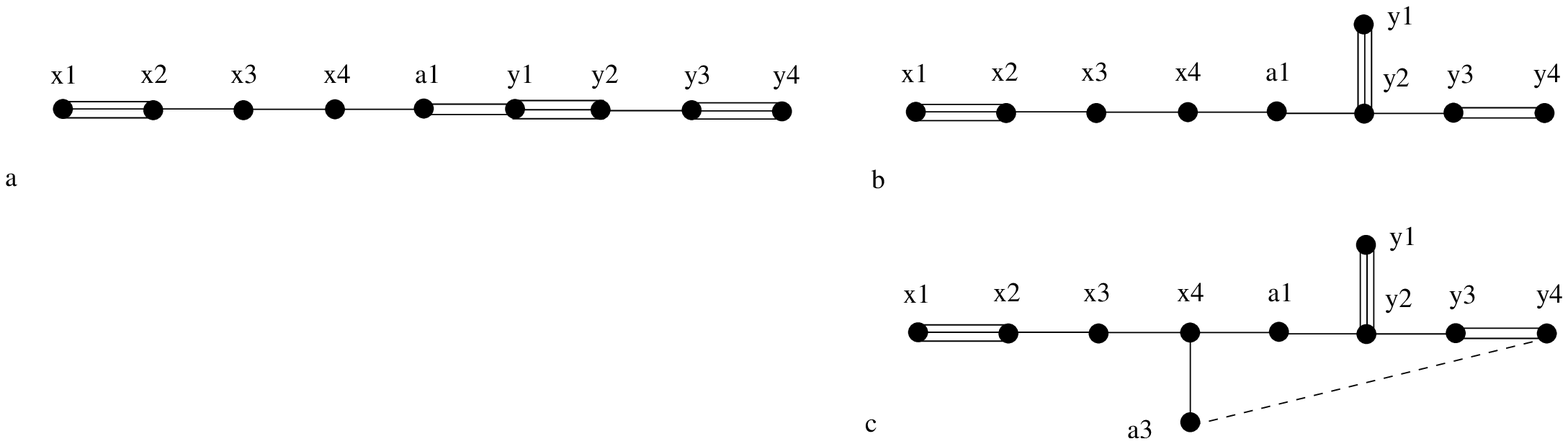,width=0.699\linewidth}
\caption{}
\label{9points}
\end{center}
\end{figure}

Consider the diagram $\[Q,a_3\]$. It contains 
at most 1 dotted edge (Lemma~\ref{2dotted}), which cannot have
endpoint neither in $S_0$ ($\[S_0,a_3\]$ is a Lann\'er diagram)
nor in $S_1$ ($a_3$ is not a bad neighbor of $S_1$).
So, either $\[S_0,a_3,y_4,S_1\]$ or
 $\[S_0,a_3,a_1,S_1\]$  is a diagram containing no
dotted edges. A direct check shows that the former diagram never
satisfies signature condition, while the latter satisfies it in a
unique case shown in Fig.~\ref{9points}(c). Moreover,
$[a_3,y_4]=\infty$, otherwise $\det(\[S_0,a_1,a_3,\o S_0\])\ne 0$.
 
Now, we are left to determine how it is possible to attach $a_2$ to the
diagram $\[S_0,a_1,a_3,\o S_0\]$.
Consider the subdiagram $S_{15}=\[y_1,y_2,a_1,x_4\]$ of the type $H_4$ 
with its bad neighbors $x_3,a_3,y_2$. 
As it was shown above, the diagram 
$\o S_{15}=\[x_1,x_2,a_2,y_4\]$ is a linear Lann\'er diagram with one  
double, one simple and one triple edge.
Hence, either
$[a_2,x_2]=2, [a_2,x_1]=3, [a_2,y_4]=4$
or
$[a_2,x_2]=3, [a_2,x_1]=2, [a_2,y_4]=4$.
Furthermore, $a_2$ is not joined with $\[x_3,x_4,a_3,a_1\]$ 
since the diagram
$\[x_2,x_3,x_4,a_3,a_1,y_2\]$ of the type $E_6$ already has
3 bad neighbors $x_1,y_1,y_4$. Also, $a_2$ is not joined with
$\[y_1,y_2\]$ since the diagram $\[y_1,y_2,a_1,x_4\]$ of the type $H_4$ has 3
bad neighbors $x_3,a_3,y_3$. Since $a_2$ is a bad neighbor of $S_1$,
$[a_2,y_3]\ne 2$. If $[a_2,y_3]= 3$, then $S_{16}=\[x_2,x_1,a_2,y_3\]$ is a
diagram of the type $H_4$ such that $\o S_{16}=\[y_1,a_1,x_4,a_3\]$ 
is a non-connected diagram, which is impossible.
If $[a_2,y_3]=5 $ (or 4)
then $\[a_2,y_3,y_2,a_1\]$ (or  $\[a_2,y_3,y_2,a_1,x_4,a_3\]$)
is a diagram of the type $H_4$ ($B_6$) with 4 bad neighbors
$y_4,y_1,x_4$ and one of $x_1$ and  $x_2$.
Hence, $ [a_2,y_3]=\infty $.

Therefore, $a_2$ is an end of three edges: one dotted edge
$a_2y_3$, one double edge $a_2y_4$ and one simple
edge (either $a_2x_1$ or $a_2x_2$). 
If $[a_2,x_1]=3$ then the subdiagram $\[S_0,a_1,a_2,a_3,y_1,y_2\]$
does not satisfy the signature condition (it contains no dotted edges,
so we check the signature directly).
Thus,  $[a_2,x_1]=2$,  $[a_2,x_2]=3$ and we arrive with the diagram
$\Sigma_{P_7}$.  

\bigskip
\noindent
All cases are considered, the only polytope $P_7$ is found.
So, the lemma is proved.

\end{proof}


\vspace{35pt}
\noindent
Independent University of Moscow\\
e-mail: 
\phantom{ow} felikson@mccme.ru\\
\smallskip
\phantom{owowwww} pasha@mccme.ru

\end{document}